\documentclass[final,5p,times,twocolumn]{elsarticle}

\usepackage{amssymb}
\usepackage{amsmath}
\usepackage{subcaption}
\usepackage{float}
\usepackage{float}

\journal{Elsevier}

\begin{document}

\begin{frontmatter}


%% \ead{yaji@mech.eng.osaka-u.ac.jp}

\title{Homogenization-based optimization of wall thickness distribution for TPMS two-fluid heat exchangers}

%% 所属
\affiliation[1]{organization={Department of Mechanical Engineering, The University of Osaka},
                addressline={2-1 Yamadaoka},
                city={Suita},
                postcode={565-0871},
                state={Osaka},
                country={Japan}}
                
\affiliation[2]{organization={Department of Mechanical Engineering, University of Maryland},
                city={College Park},
                postcode={20742},
                state={MD},
                country={USA}}

\author[1]{Kaito Ohtani}
\author[1]{Hiroki Kawabe}
\author[1]{Kentaro Yaji\corref{cor1}}
\cortext[cor1]{Corresponding author}
\ead{yaji@mech.eng.osaka-u.ac.jp}
\author[1]{Kikuo Fujita}
\author[2]{Vikrant Aute}

%% Abstract
\begin{abstract}
Triply Periodic Minimal Surface (TPMS) structures are attracting growing attention as promising geometries for next-generation high-performance heat exchangers (HXs), due to their continuous flow paths and high surface-area-to-volume ratio that enhance heat transfer performance.
Among these, graded TPMS structures with spatially varying thickness have emerged as a potential means to further improve performance.
This study proposes an optimization method of wall thickness distribution for TPMS HXs based on an effective porous media model, which allows accurate performance prediction while significantly reducing computational cost.
The proposed method is applied to a gyroid two-fluid HX aiming to improve the thermal-hydraulic performance.
Furthermore, full-scale numerical simulations of the optimized-thickness design show a 12.2\% improvement in the performance evaluation criterion (PEC) compared to the uniform-thickness design.
The improvement is primarily attributed to the optimized non-uniform wall thickness, which directs more flow toward the core ends and enhances velocity uniformity.
As a result, heat transfer is enhanced at the core ends, leading to more effective use of the entire HX core and improved overall thermal performance.
\end{abstract}


% \begin{highlights}
% \item An effective porous model and an optimization method are proposed for TPMS HXs.
% \item The validity of the effective model is confirmed by comparison with a full-scale model.
% \item Full-scale analysis reveals the mechanism of enhanced heat transfer in the optimized-thickness HX.
% \item Optimized design shows higher PEC than uniform-thickness, improving performance.
% \end{highlights}

%% Keywords

\end{frontmatter}



\section{Introduction}
\label{sec1}

Rapid advances in science and technology have substantially increased the demand for effective thermal management solutions across various industrial sectors.
Heat exchangers (HXs) are one of the most important components in thermal management, playing vital roles in a wide range of applications such as air conditioning, refrigeration, automotive cooling, and industrial process systems.
Additionally, advances in additive manufacturing (AM) technology have enabled the fabrication of products with complex internal geometries that were once difficult to produce \cite{frazier2014,Tofail2018}.
These advances in manufacturing capabilities have facilitated the practical implementation of lattice structures and triply periodic minimal surface (TPMS) structures, which offer superior mechanical properties \cite{Yuan2019,Feng2022}.

TPMS structures have attracted significant attention for HX applications due to their high surface-area-to-volume ratio and continuous flow channels, promoting excellent heat and fluid transport performance \cite{Kaur2021,Yeranee2022,Dutkowski2022,Gado2024}.
Li et al.~\cite{Li2020} designed HXs based on TPMS structures and demonstrated through numerical analysis that these HXs outperform conventional printed circuit heat exchangers (PCHEs).
In addition, Iyer et al.~\cite{Iyer2022}, Li et al.~\cite{Li2022}, and Yan et al.~\cite{Yan2023} confirmed through simulations that TPMS HXs exhibit superior heat transfer performance compared to conventional designs.
Experimental validations using AM-fabricated prototypes further confirm these advantages, supporting their application in compact and efficient thermal management systems \cite{Dixit2022,Liang2023,Reynolds2023,Yan2024,Qian2024}.
These findings underscore the potential of TPMS HXs as promising candidates for next-generation thermal management solutions.

Recent studies have developed various design approaches to enhance TPMS HX performance.
Some emerging research has focused on applying conventional heat transfer enhancement mechanisms to TPMS-based HXs. 
Geometry manipulation, such as stretching, was investigated by several researchers \cite{Mahmoud2023,Yan2024b,Coe2025}. 
These studies demonstrated that stretching the TPMS structure improves both thermal and hydraulic performance.
Similar to microfins or grooves used to enhance in-tube performance, Qin et al.~\cite{Qin2025} investigated the use of wrinkles, which are generated by a modification of the gyroid equation using a control factor term. They reported a maximum increase in Nu of 3.5\% and the corresponding increase in heat transfer coefficient of 8\%.
In addition to conventional enhancement techniques, researchers have also explored strategies that leverage the unique characteristics of TPMS structures to further improve heat transfer and hydraulic performance.
Tang et al.~\cite{Tang2023} and Barakat et al.~\cite{Barakat2024} proposed deformation control parameters that alter the TPMS structure's shape to enhance heat transfer performance.
Barakat et al.~proposed novel TPMS lattice structures that outperformed Gyroid and Diamond based HX designs \cite{Barakat2024new}.
Dharmalingam et al.~developed a rapid and practical optimization framework for TPMS HXs by constructing CFD-based correlations using the smallest repeatable section \cite{Dharmalingam2025}.
The results showed that reducing the hydraulic diameter enhanced convective heat transfer and decreased thermal resistance, achieving up to ten times higher volumetric power density compared to the initial design.
Although these approaches have successfully improved heat transfer and hydraulic performance, they have focused on TPMS HXs with uniform wall thickness distribution.
Therefore, there remains scope for further performance improvement by considering non-uniform wall thickness distribution.
Investigating these spatial variations has motivated numerous recent studies to improve TPMS HX performance.

Functionally graded TPMS structures have attracted growing interest due to their potential to enhance HX performance \cite{Yu2019,AlKetan2019}. 
Zhang et al.~developed a conformal design approach to improve flow uniformity \cite{Zhang2025conformal}.
Another study by Zhang et al.~proposed a gradient design in which the unit cell size varies radially \cite{Zhang2025gradient}.
Oh et al.~introduced a mathematical gradient framework to enable localized morphological changes \cite{Oh2023} and implemented graded variations in cell size while maintaining constant wall thickness \cite{Oh2025}.
Moreover, Wang et al.~\cite{Wang2025} and Qian et al.~\cite{Qian2025} proposed hybrid graded TPMS HXs combining Gyroid and Diamond structures, achieving superior overall thermal-hydraulic performance.
These studies demonstrate that graded TPMS designs offer a promising approach to both enhance heat transfer and hydraulic performance. 
However, most design methods for graded TPMS HXs rely on empirical or experimental approaches, and optimization strategies are still underdeveloped.
Previous studies have shown that wall thickness strongly influences thermal and hydraulic characteristics \cite{Samson2023,Attarzadeh2021,Wang2023}.
Therefore, determining the optimal wall thickness distributions is essential for further improving graded TPMS HX performance.
From a macroscopic perspective, TPMS cells with different wall thicknesses can be treated as different materials.
Consequently, the optimization of wall thickness distributions for TPMS HXs can be formulated as a material distribution problem.

Topology optimization (TO) \cite{Bendsoe1988} is a well-established method for optimizing material distribution.
It offers a high degree of design freedom by representing the material distribution using design variables.
TO has been actively studied for thermal-fluid systems \cite{Dbouk2017,Alexandersen2020} and has recently been applied to two-fluid HXs \cite{Fawaz2022}.
For example, density-based \cite{Hoghoj2020,Kobayashi2021} and level-set-based \cite{Feppon2021} TO methods have been developed for two-fluid HXs, both aiming to maximize heat transfer under pressure constraints while preventing fluid mixing.
Galanos et al.~\cite{Galanos2022} proposed an optimization method that combines topology and shape optimization.
Pimanov et al.~\cite{Pimanov2025} proposed a sparse narrow-band TO framework. 
This framework restricts design updates to the fluid-solid interface and employs an efficient solver, enabling effective optimization of large-scale two-fluid HXs.
However, two-fluid HXs optimized using conventional TO are often unmanufacturable due to their complex geometries.
For industrial applications, it is essential to either incorporate manufacturing constraints or develop TO methods based on manufacturable structures.

In recent years, optimization of TPMS structures has gained significant attention due to its manufacturability via AM.
However, directly optimizing the detailed microscale features of TPMS geometries is computationally prohibitive due to their complexity.
To address this challenge, researchers have employed a homogenization approach.
Due to the periodic and spatially regular nature of TPMS structures, they can be represented as homogenized materials with equivalent effective properties.
This method enables the prediction of macroscopic behavior without resolving detailed microscale geometries, substantially reducing computational cost.
In addition, the effective properties can be linked to design variables of TO.
Therefore, by treating TPMS structures with different wall thicknesses, cell sizes, or cell types as different homogenized materials and incorporating them into TO frameworks, we can effectively optimize graded TPMS structures.
These approaches have been pioneered mainly in structural problems, including studies on wall thickness optimization \cite{Li2019,Stromberg2021}, sequential optimization of cell size followed by wall thickness \cite{Hu2022} and its extension to thermal conduction problems \cite{Wang2022}, optimization considering multiple TPMS cell types \cite{Feng2022topology,Zhao2023}, and optimization including void phases \cite{Stromberg2024}.
Taken together, these efforts highlight the potential of homogenization-based TO as a powerful tool for improving the functional performance of graded TPMS structures.

Recently, researchers have extended homogenization-based TO to thermal-fluid problems.
For large-scale thermal-fluid systems, directly resolving detailed geometries requires an extremely large number of mesh elements, making full-scale simulations often computationally prohibitive.
To overcome this limitation, porous media models have been widely adopted \cite{Missirlis2010,Hooman2010,Yang2014,An2018,Juan2018,Kim2019,Li2021,Chen2023,Wang2024}.
They simplify the heat exchange regions by treating them as porous media, reducing computational cost while maintaining sufficient accuracy in predicting pressure drop and heat transfer.
These models also enable the application of homogenization-based TO to thermal-fluid systems with periodic microstructures.
In previous studies, lattice structures of varying diameters and TPMS structures with different wall thicknesses have been modeled as homogenized materials with distinct effective properties, and their lattice diameters and wall thickness distributions have been optimized based on a TO framework \cite{Takezawa2019,Men2025}.
Takezawa et al.~modeled lattice structures as porous media and optimized diameter distributions for water-cooled heat sinks \cite{Takezawa2019}, while Men et al.~optimized wall thickness distributions in TPMS heat sinks \cite{Men2025}.
However, these studies are limited to single-fluid systems and cannot consider heat exchange between two fluids separated by partition walls.
Given the high surface area and favorable thermal-fluid transport properties of TPMS structures, extending homogenization-based TO to two-fluid systems is an important engineering challenge.
To enable homogenization-based TO for such systems, it is essential to develop a porous media model that accurately captures heat exchange between each fluid and the wall.

In this paper, we propose an optimization method of wall thickness distribution for TPMS HXs based on an effective porous media model, which replaces the detailed microscale geometry with a homogenized porous model to reduce computational cost.
To account for heat exchange between the two fluids in the effective model, we introduce effective heat transfer coefficients to characterize the thermal interaction between each fluid and the wall.
The effective properties are obtained from finite element simulations of the gyroid unit cell, and their relationship with the design variables is approximated using polynomials.
The optimization objective is defined as a weighted sum of the heat transfer rate and the pressure drop, and TO is conducted for multiple sets of weighting factors.
The optimized designs are reconstructed into full-scale gyroid structures via a dehomogenization process, and their performance is rigorously assessed through high-fidelity simulations.
This framework extends the application of homogenization-based TO using porous media models to two-fluid systems, which had previously been limited to single-fluid systems.
The effectiveness and limitations of the proposed approach are examined through numerical case studies.

The remainder of this paper is organized as follows.
In Section \ref{sec2}, we describe the effective porous media model for the TPMS two-fluid HX, the calculation methods for the effective material properties, and the full-scale analysis model.
In Section \ref{sec3}, we present the optimization model and the formulation of the optimization problem.
In Section \ref{sec4}, we validate the full-scale analysis model, verify the effective porous media model, and discuss the results of the optimization and full-scale analysis.
Finally, Section \ref{sec5} concludes this study and discusses future work.

\section{Numerical Modeling}
\label{sec2}

This section presents the mathematical foundations of the effective porous media model for TPMS HXs.
The core concept of the proposed model is the incorporation of heat exchange between each fluid and the wall into the effective thermal formulation.
The implementation of this mechanism and the overall construction of the model are presented.
In addition, the methodology for computing the effective properties and the metrics for evaluating full-scale analysis results are also presented.

\subsection{Gyroid Design with Controllable Wall Thickness}
\label{sec2.1}

Triply periodic minimal surfaces (TPMS) are periodic surfaces without self-intersections that extend infinitely along three principal axes \cite{Han2018}. 
The first example of TPMS was discovered by Schwarz in 1865 \cite{Schwarz1890}, and the gyroid structure, which is the focus of the present work, was discovered by Alan Schoen in 1970 \cite{Schoen1970}.
The gyroid surface is represented by the following equation:
\begin{equation}
\sin X \cos Y + \sin Z \cos X + \sin Y \cos Z = c/L_\mathrm{cell},
\label{gyroid_eq}
\end{equation}
where $c$ is the level-set parameter, $L_{\mathrm{cell}}$ is the unit cell size, and
$ X = 2\pi x / L_{\mathrm{cell}}, \quad Y = 2\pi y / L_{\mathrm{cell}}, \quad Z = 2\pi z / L_{\mathrm{cell}}$, where $\mathbf{x} = (x, y, z)$ represents a spatial position.
When $c$ is set to zero, the surface defined by Eq.~\eqref{gyroid_eq} partitions the space into two regions of equal volume.

%============================ 
\begin{figure}[bht]
    \centering
    \begin{minipage}[t]{0.45\textwidth}
        \centering
        \includegraphics[width=0.95\linewidth]{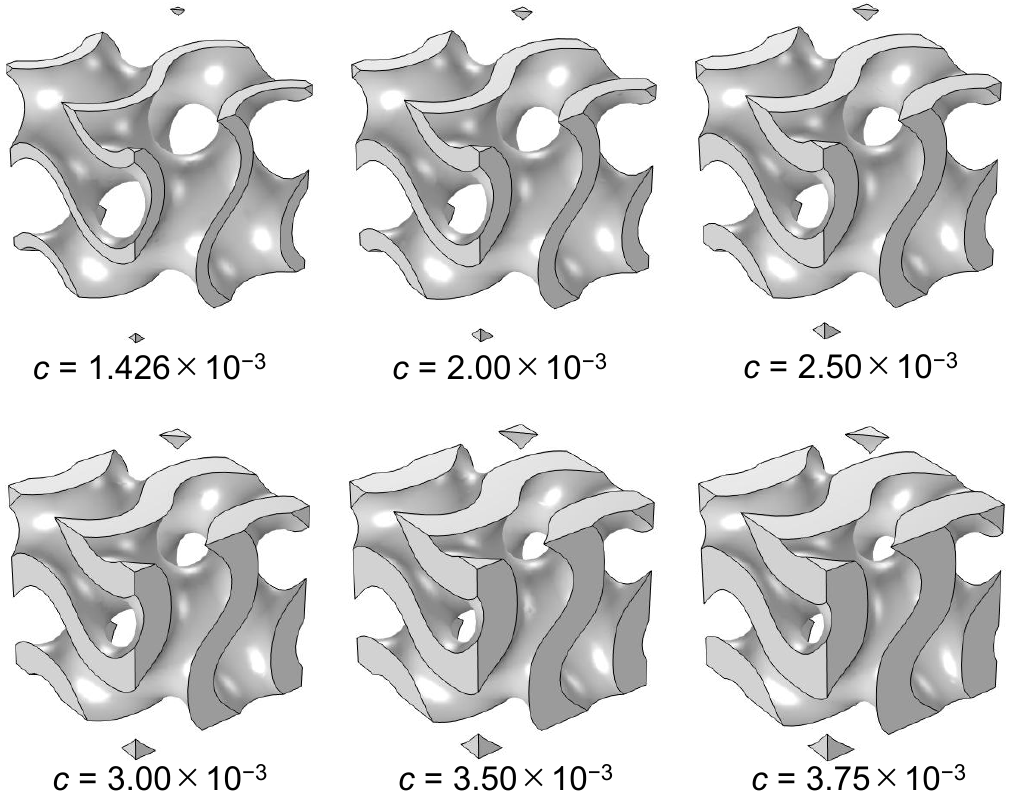}
        \subcaption{Gyroid unit cell with constant parameter $ c $ }
        \label{fig1a}
    \end{minipage} \\
    \vspace{0.2cm}
    \begin{minipage}[t]{0.45\textwidth}
        \centering
        \includegraphics[width=0.95\linewidth]{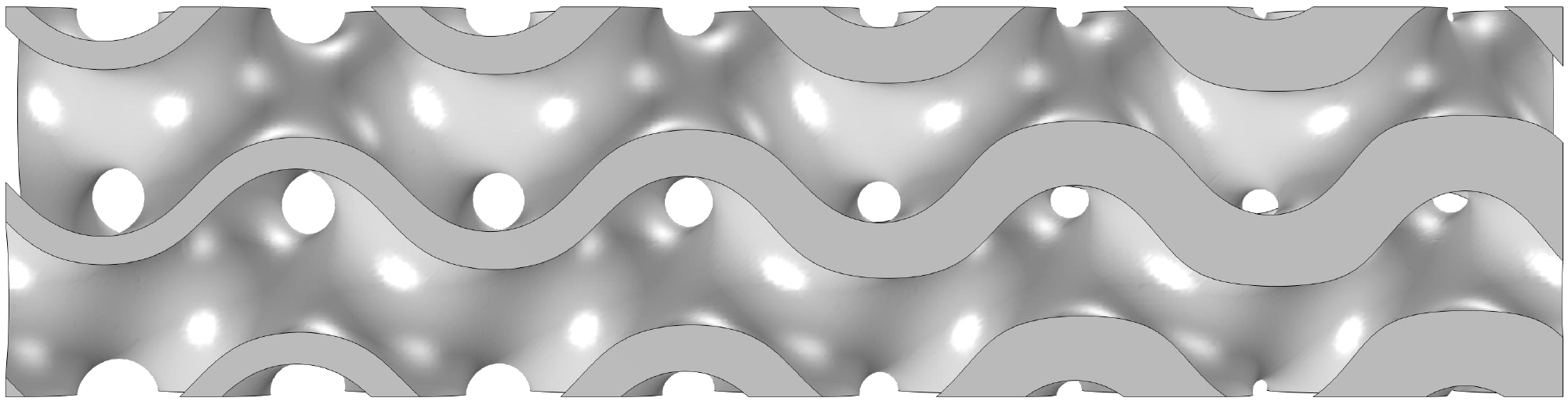}
        \subcaption{Gyroid structure with spatially varying $ c $}
        \label{fig1b}
    \end{minipage}
    %\vspace{-0.2cm}
    \caption{Illustration of gyroid structures showing the effect of the parameter $c$ on wall thickness for $L_{\mathrm{cell}} = 4.60 \times 10^{-3}$ [m]}
    \vspace{-0.2cm}
    \label{fig1}
\end{figure}
%============================
To apply gyroid structures to HXs, it is necessary to define solid domains.
In this study, the solid domain of the gyroid structure, denoted as $\Omega_s$, is defined using two level-set functions as follows \cite{Hu2022}:
\begin{equation} \label{eq:gyroid}
\begin{aligned}
   G_1(\mathbf{x}) &= c/L_\mathrm{cell} +  \left(\sin X \cos Y + \sin Z \cos X + \sin Y \cos Z \right), \\
   G_2(\mathbf{x}) &= c/L_\mathrm{cell} - \left( \sin X \cos Y + \sin Z \cos X + \sin Y \cos Z \right), \\
   \Omega_s &= \{ \mathbf{x} \mid G_1(\mathbf{x}) \geq 0, \, G_2(\mathbf{x}) \geq 0 \}.
\end{aligned}
\end{equation}
These level-set surfaces are generated by translating the reference surface defined by $c = 0$ an equal distance in the positive and negative directions.
These solid domains divide the space into two fluid regions of equal volume.

Fig.~\ref{fig1} illustrates the effect of the parameter $c$ on the gyroid structure.
Fig.~\ref{fig1}(a) shows a gyroid unit cell with a constant value of $c$, where it can be seen that the wall thickness increases as $c$ increases.
Fig.~\ref{fig1}(b) shows a gyroid structure generated by continuously varying $c$ spatially.
As demonstrated in Fig.~\ref{fig1}(b), smoothly varying $c$ results in a corresponding continuous change in the wall thickness of the gyroid structure.

In this study, functionally graded gyroid two-fluid HXs are designed by spatially varying $c$ within the gyroid HX core domain.
This approach enables the optimization of wall thickness distribution, thereby improving the thermal and hydraulic performance of the gyroid HXs.

\subsection{Effective Porous Media Model for TPMS two-fluid HXs}
\label{sec2.2}

Direct optimization of TPMS HX structures is challenging due to the high computational cost associated with their complex internal geometry.  
To address this issue, a porous media approximation is considered effective, as it can significantly reduce the computational cost.  
Toward this goal, previous studies have utilized the periodicity of TPMS structures to construct effective porous media models by treating them as homogeneous porous media \cite{Men2025}. Such models enable the prediction of physical behavior without explicitly resolving the complex internal geometry, thereby substantially reducing computational cost. However, an effective model applicable to the optimization of two-fluid TPMS HXs has not yet been established.
In this study, we propose an effective porous media model that incorporates heat exchange between each fluid and the solid wall to overcome this limitation.

This study employs two types of models: an effective porous media model for optimization calculations, and a Navier-Stokes (NS) model for full-scale analysis and the computation of effective properties. 
The following assumptions apply to both models:
\renewcommand{\labelenumi}{(\arabic{enumi})}  
\begin{enumerate}  
  \item The flow is steady and incompressible.  
  \item The flow is laminar.  
  \item The material properties of both the fluid and the solid are constant and temperature-independent.  
   \item The hot and cold fluids are assumed to be the same fluid.  
  \item Thermal radiation effects are not considered.
\end{enumerate}  
Additionally, the effective porous media model incorporates the following additional assumption:
\begin{enumerate}  
    \setcounter{enumi}{5}
  \item The effective properties are assumed to be isotropic.  
\end{enumerate}  
Based on these assumptions, each model is developed and analyzed.

Next, we describe the effective porous media model.
Flow through porous media is typically described by Darcy's law, which accounts only for viscous losses.
However, when the flow velocity is relatively high or when the medium has a complex internal geometry, such as in TPMS structures, inertial effects cannot be neglected.  
In such cases, the Brinkman-Forchheimer model \cite{Vafai1981, Nield2006}, which extends the Brinkman equation \cite{Brinkman1947} by incorporating the Forchheimer term \cite{Forchheimer1901}, is commonly used.  
This model includes the Darcy resistance, viscous diffusion, and inertial loss terms, making it suitable as an effective porous media model for TPMS structures.  
In this study, the porous flow model for the TPMS two-fluid HXs is expressed as follows:  
\begin{equation} \label{eq:BF_eq}
\begin{aligned}
    & \nabla \cdot \mathbf{U}_i = 0, \\
    &\frac{\rho}{\varepsilon_i} \left( \mathbf{U}_i \cdot \nabla \right) \frac{\mathbf{U}_i}{\varepsilon_i} = \nabla \cdot \left[ -p_i \mathbf{I} + \mathbf{K} \right] - \alpha_i \mathbf{U}_i - \beta_i |\mathbf{U}_i| \mathbf{U}_i,   \\
    & \mathbf{K} = \frac{\mu}{\varepsilon_i} \left( \nabla \mathbf{U}_i + (\nabla \mathbf{U}_i)^{\mathrm{T}} \right) - \frac{2}{3} \frac{\mu}{\varepsilon_i} (\nabla \cdot \mathbf{U}_i) \mathbf{I},
\end{aligned}
\end{equation}
where $\mathbf{U}_i$ is the Darcy velocity vector of fluid $i$ ($i = 1, 2$), defined as the volumetric flow rate per unit cross-sectional area within the porous medium. 
The parameters $\varepsilon_i$, $p_i$, $\alpha_i$, and $\beta_i$ denote the porosity, pressure, viscous resistance coefficient, and inertial resistance coefficient of fluid $i$, respectively. 
The viscosity and fluid density are represented by $\mu$ and $\rho$, respectively.
Using this model enables efficient and reasonable prediction of the macroscopic flow behavior through TPMS structures, without directly resolving the complex geometry.

%%%%%%%%%%%%%%%%%%%%%%%%%%%%%%%%%%%%%%%%%
\begin{figure}[htb]
    \begin{center}
        \includegraphics[width=0.95\linewidth]{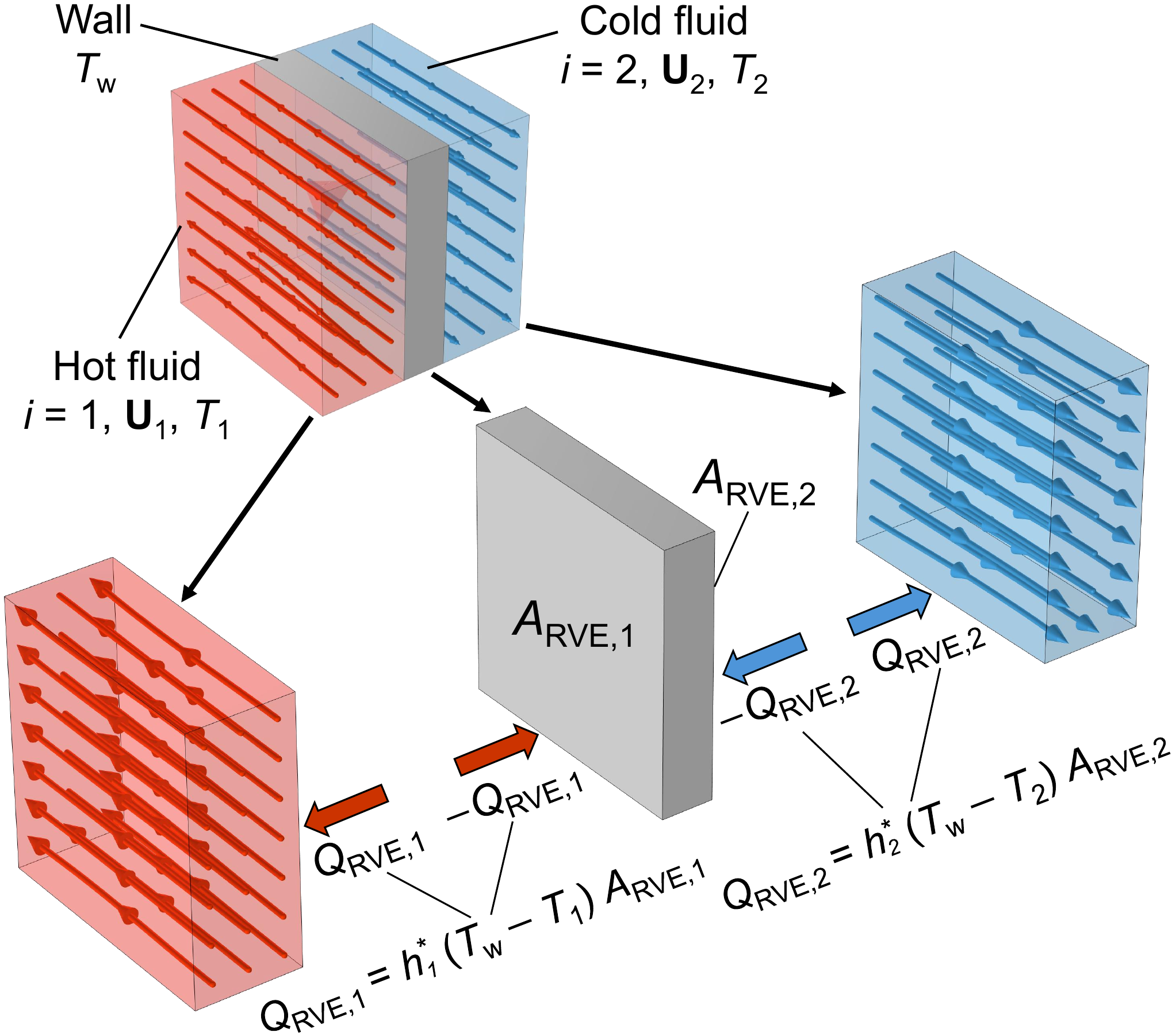}
        \caption{Schematic of heat exchange between fluid $i$ and the wall within a representative volume element (RVE)}
        %\vspace{-0.3cm}
        \label{fig2}
    \end{center}%
    %\vspace{-0.2cm}
\end{figure}
%%%%%%%%%%%%%%%%%%%%%%%%%%%%%%%%%%
Next, we discuss the effective heat transfer model.
It is essential to incorporate heat exchange between fluids and solids into the effective thermal model when constructing a porous media model for TPMS two-fluid HXs. 
We consider the representative volume element (RVE) of a two-fluid HX and propose incorporating the effect of heat exchange into the effective thermal model by expressing the heat exchange within the RVE in terms of its representative values.
Fig.~\ref{fig2} shows a schematic illustration of the heat exchange within the RVE.
As shown in Fig.~\ref{fig2}, both the hot and cold fluids exchange heat with the wall.
The total heat exchange of fluid $i$ within the RVE, denoted as $Q_{\mathrm{RVE},i}$, can be microscopically expressed as follows:
\begin{align}
    Q_{\mathrm{RVE},i} = \int_{\Gamma_i} h_i(\mathbf{x}) \left( \tilde{T}_{\mathrm{w}}(\mathbf{x}) - \tilde{T}_i(\mathbf{x}) \right) dA,
    \label{micro_heat_exchange}
\end{align}
where $\Gamma_i$ denotes the wall boundary in contact with fluid $i$, and $h_i(\mathbf{x})$ represents the local heat transfer coefficient between the wall and fluid $i$.
$\tilde{T}_{\mathrm{w}}(\mathbf{x})$ and $\tilde{T}_i(\mathbf{x})$ denote the wall and fluid temperatures, respectively.
By appropriately selecting the RVE to capture the average flow and temperature behavior, Eq.~\eqref{micro_heat_exchange} can be approximated using representative values within the RVE as follows:
\begin{align}
    Q_{i,\mathrm{RVE}} = h_i^* \left( {T}_{\mathrm{w}} - T_i \right) A_{\mathrm{RVE},i},
    \label{macro_heat_exchange}
\end{align}
where $h_i^*$ is the effective heat transfer coefficient of fluid $i$, ${T}_{\mathrm{w}}$ is volume-averaged wall temperature, $T_i$ is volume-averaged temperature of fluid $i$, and $A_{\mathrm{RVE},i}$ is the wall-fluid $i$ interfacial area in the RVE.
The effective heat transfer coefficient $h_i^*$ can be defined as:
\begin{align}
    h_i^* = \frac{\int_{\Gamma_i} \mathbf{q}_i \cdot \mathbf{n}_i \ dA}{A_{\mathrm{RVE},i} \left( {T}_{\mathrm{w}} - T_i \right)} ,
    \label{effective_heat_transfer_coefficient}
\end{align}
where $\mathbf{q}_i$ is the heat flux vector on the wall surface in contact with fluid $i$, and $\mathbf{n}_i$ is the outward unit normal vector on the wall surface $\Gamma_i$ with respect to the fluid domain.
Furthermore, dividing Eq.~\eqref{macro_heat_exchange} by the volume of the RVE, $V_{\mathrm{RVE}}$, yields the volumetric heat exchange rate:
\begin{align}
    Q_i^* = h_i^* \left( {T}_{\mathrm{w}} - T_i \right) A_{\mathrm{RVE},i}/V_{\mathrm{RVE}}.
\end{align}
By incorporating this volumetric heat source term, the effective thermal models for both fluid $i$ and the solid, accounting for the heat exchange effects, are formulated as follows:
\begin{equation}
    \rho c_p  \mathbf{U}_i \cdot \nabla  T_i + \nabla \cdot (-{k_{\text{f}}^*}_{,i} \nabla T_i)=Q_i^*,
    \label{eq:energy_fluid}
\end{equation}
\begin{equation}
    \nabla \cdot (-k_{\mathrm{s}}^* \nabla T_{\text{w}}) = -Q_1^* - Q_2^*,
    \label{eq:energy_solid}
\end{equation}
where $c_p$ is the specific heat capacity of the fluid, ${k_{\mathrm{f}}^*}_{,i} $ is the effective thermal conductivity of fluid $i$, and $k_{\mathrm{s}}$ is the effective thermal conductivity of the solid.

Porous media modeling of microchannel structures, in which heat exchange between fluid and solid is represented as a volumetric heat source term, has also been proposed in previous studies focusing on compact tube HXs \cite{Chen2023} and TPMS HXs \cite{Wang2024}.
In these studies, one of the two fluids was treated as a solid domain to simplify the porous media model.
Moreover, the heat transfer coefficient for fluid-solid surface was prescribed based on existing correlations.
In contrast, the porous media model developed in the present study adopts an RVE approach, wherein both fluids are explicitly modeled.
Furthermore, the effective heat transfer coefficient used in the model is analytically derived from microscale transport phenomena within the RVE, enabling a more accurate evaluation of the heat exchange between each fluid and the solid.

By solving the coupled equations, Eq.~\eqref{eq:BF_eq}, Eq.~\eqref{eq:energy_fluid}, and Eq.~\eqref{eq:energy_solid}, the flow and heat transfer behavior in TPMS two-fluid HXs can be predicted at a reasonable computational cost.
This effective porous media model is not limited to TPMS structures but can also be applied to porous media modeling of two-fluid HXs with other periodic architectures.

\subsection{Computing Effective Properties}
\label{sec2.3}

To solve the effective porous media model constructed in Section \ref{sec2.2}, it is necessary to compute the effective properties.  
The effective properties include porosity $\varepsilon_i$, surface area within the RVE $A_{\mathrm{RVE},i}$, viscous resistance coefficient $\alpha_i$, inertial resistance coefficient $\beta_i$, effective thermal conductivity of fluid $i$ $k^*_i$, effective thermal conductivity of solid $k_{\mathrm{s}}$, and effective heat transfer coefficient $h^*_i$.  
These effective properties are calculated based on the finite element method (FEM).
In simulations that do not employ the effective porous medium model, the following Navier-Stokes and energy equations are used:
\begin{equation} \label{eq:NS}
\begin{aligned}
    & \nabla \cdot \mathbf{u} = 0,\\
    & \rho \left( \mathbf{u} \cdot \nabla \right) \mathbf{u}
    = - \nabla p + \mu \nabla \cdot \left[ \nabla \mathbf{u} + \left( \nabla \mathbf{u} \right)^\mathsf{T} \right], \\
    & \rho c_p \mathbf{u} \cdot \nabla T 
    = \nabla \cdot \left( k \nabla T \right),
\end{aligned}
\end{equation}
where $\mathbf{u}$ is the velocity vector, $p$ is the pressure, $k$ is the thermal conductivity, and $T$ is the temperature.
Eq.~\eqref{eq:NS} is used for the calculation of the effective properties.

%%%%%%%%%%%%%%%%%%%%%%%%%%%%%%%%%%%%%%%%%
\begin{figure}[htb]
    \begin{center}
        \includegraphics[width=0.985\linewidth]{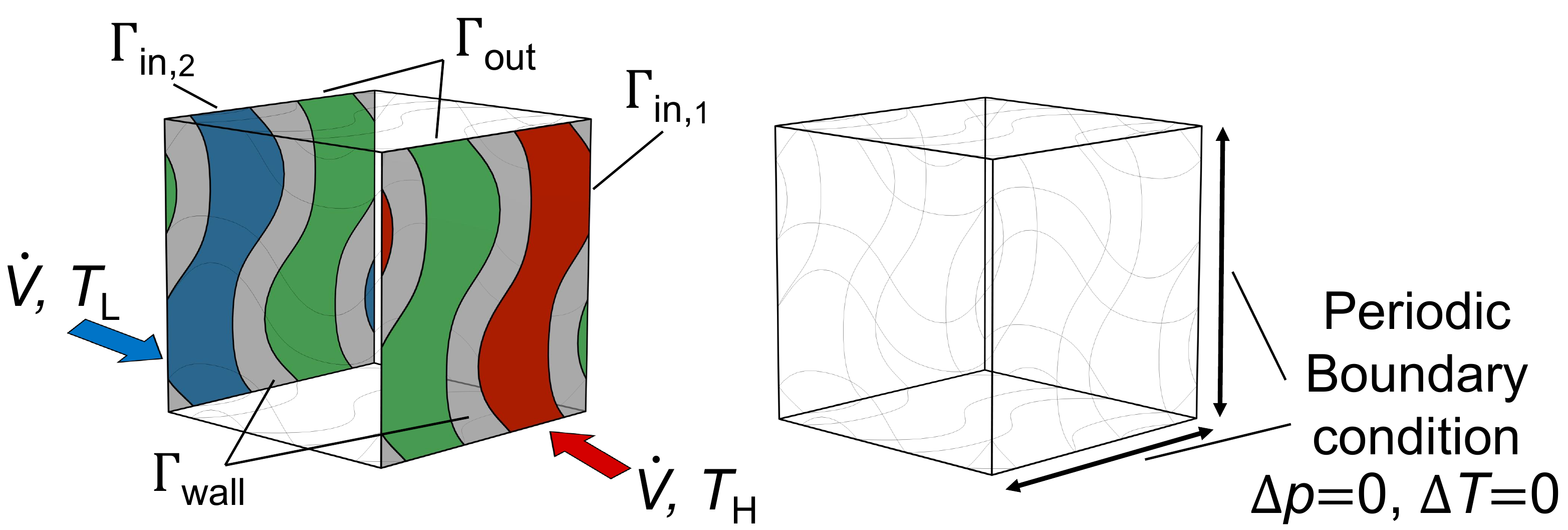}
        \caption{Boundary conditions applied to the RVE used for computing effective properties $\alpha$, $\beta$, and $h^*_i$}
        %\vspace{-0.3cm}
        \label{fig3}
    \end{center}%
    %\vspace{-0.2cm}
\end{figure}
%%%%%%%%%%%%%%%%%%%%%%%%%%%%%%%%%%
In this study, we focus on the periodicity of TPMS structures and quantitatively evaluate the corresponding effective properties through numerical homogenization based on the RVE method.
The RVE is defined as a region that statistically represents the entire material \cite{Hill1963}.
For periodic structures, macroscopic effective properties can be appropriately evaluated by using a unit cell with periodic boundary conditions.
Accordingly, the unit cell of the gyroid structure is adopted as the RVE in this study.

The gyroid structure defined by Eq.~\eqref{eq:gyroid} partitions the space into two geometrically identical regions with equal volume fractions, and the material properties of the hot and cold fluids are identical in this study.  
Therefore, the porosity $\varepsilon_i$, the wall-fluid $i$ interfacial area in the RVE $A_{\mathrm{RVE},i}$, the viscous resistance coefficient $\alpha_i$, the inertial resistance coefficient $\beta_i$, and the effective thermal conductivity of fluid $i$ $k^*_i$ are the same for fluid 1 and fluid 2. 
Accordingly, in this study, the subscript $i$ for these effective properties is omitted, and they are hereafter denoted as $\varepsilon$, $A_{\mathrm{RVE}}$, $\alpha$, $\beta$, $k^*$.
Regarding the effective heat transfer coefficient $h^*_i$, its value depends on the magnitude of the Darcy velocity $|\mathbf{U}_i|$, so $h^*_1$ and $h^*_2$ may differ at the same location. 
However, due to the symmetry of the gyroid structure and the identical material properties of the fluids, the interpolation function expressed in terms of $|\mathbf{U}_i|$ is identical for both fluids.

The porosity $\varepsilon$ is defined as the ratio of the volume of fluid $i$ to the total volume of the RVE:
\begin{equation}
    \varepsilon = \frac{V_{\mathrm{f},i}}{V_{\mathrm{RVE}}} = \frac{V_{\mathrm{f},i}}{L_{\mathrm{cell}}^3},
    \label{eq:porosity}
\end{equation}
where $V_{\mathrm{f},i}$ represents the volume of each fluid region within the RVE and $V_{\mathrm{RVE}}$ is the total volume of the RVE.
The contact surface area between the wall and each fluid in the RVE, $A_{\mathrm{RVE}}$, is directly measured from the geometry of the gyroid unit cell.

To compute the coefficients $\alpha$, $\beta$, and the effective heat transfer coefficient $h^*_i$, coupled thermofluid simulations are conducted on the RVE.
Fig.~\ref{fig3} illustrates the boundary conditions applied to the RVE for deriving the effective properties.
A counterflow configuration is assumed in this study: the hot fluid enters through the inlet surface $\Gamma_{\mathrm{in},1}$ with temperature $T = T_\mathrm{H}$ and volumetric flow rate $\dot{V}$, while the cold fluid enters through the inlet surface $\Gamma_{\mathrm{in},2}$ with temperature $T = T_\mathrm{L}$ and the same flow rate $\dot{V}$.
Both fluids exit through the outlet surface $\Gamma_{\mathrm{out}}$, where a pressure boundary condition $p = 0$ is imposed.
The wall surface $\Gamma_{\mathrm{wall}}$, located on the same plane as the inlet surfaces, is treated as adiabatic.
To capture the average flow behavior, fully developed flow conditions are imposed at the inlet and outlet, and periodic boundary conditions are applied to the lateral surfaces.
Simulations are performed for various gyroid unit cells with different wall thicknesses and a range of flow rates.
The wall thickness is controlled by varying the level-set parameter $c$, which defines the gyroid structure based on Eq.~\eqref{eq:gyroid}.

The pressure drop per unit length, $-\Delta p / L_{\mathrm{cell}}$, is described as a function of the Darcy velocity component $U$ aligned with the pressure drop direction, the viscous resistance coefficient $\alpha$, and the inertial resistance coefficient $\beta$.
Specifically, the coefficients $\alpha$ and $\beta$ are determined by fitting the relationship between the pressure gradient and $U$ using the following equation:
\begin{equation}\label{eq:Darcy_Forchheimer}
-\frac{\Delta p}{L_{\mathrm{cell}}} = -\alpha U - \beta U^2.
\end{equation}
This equation corresponds to the Darcy-Forchheimer model, which accounts for both viscous and inertial effects in porous media flow \cite{Nield2006}.

The effective heat transfer coefficient $h^*_i$ can be calculated using Eq.~\eqref{effective_heat_transfer_coefficient}, and its value depends on the wall thickness and the Darcy velocity.
As illustrated in Fig.~\ref{fig3}, each fluid in the RVE is subjected to the same flow rate, resulting in $h_i^* = h_1^* = h_2^*$ for the specified flow rate $\dot{V}$.
In this study, thermofluid simulations using the RVE are performed to calculate the effective heat transfer coefficient $h^*_i$ for each combination of wall thickness and flow rate.
The resulting $h^*_i$ values are then fitted as a function of wall thickness and the magnitude of the Darcy velocity to construct the interpolation function.

Next, we describe the method for computing the effective thermal conductivity.
To derive the effective thermal conductivity, a unit temperature difference $\Delta T = 1$ is imposed across a pair of opposing surfaces of the target region within the RVE, and the heat conduction equation is solved.
When computing the effective thermal conductivity of the fluid, the unit temperature difference is applied to only one of the fluid regions.
Isothermal boundary conditions are applied on the surfaces where the temperature difference is imposed.
Periodic boundary conditions are applied on the surrounding surfaces perpendicular to these surfaces, and adiabatic boundary conditions are applied on the remaining surfaces.
Under these boundary conditions, the effective thermal conductivities of the fluid and solid are calculated based on Fourier’s law $\mathbf{q} = -k \nabla T$ using the following equations:
\begin{align}
\frac{1}{A_{\mathrm{cell}}} \int_{\Psi_{\mathrm{f}}} \mathbf{q} \cdot \mathbf{n} \, dA &= -k_{\mathrm{f}}^* \frac{\Delta T}{L_{\mathrm{cell}}}, \\
\frac{1}{A_{\mathrm{cell}}} \int_{\Psi_{\mathrm{s}}} \mathbf{q} \cdot \mathbf{n} \, dA &= -k_{\mathrm{s}}^* \frac{\Delta T}{L_{\mathrm{cell}}},
\end{align}
where $ A_{\mathrm{cell}} $ is the cross-sectional area of the unit cell ($ A_{\mathrm{cell}} = L_{\mathrm{cell}}^2 $), $ \Psi_{\mathrm{f}} $ and $ \Psi_{\mathrm{s}} $ denote the surfaces on the high-temperature side where the unit temperature difference is applied in the fluid and solid regions, respectively, $ \mathbf{q} $ is the heat flux vector, $ \mathbf{n} $ is the outward unit normal vector on the surface.

\subsection{Full-Scale Performance Evaluation}
\label{sec2.4}

In this study, commonly used nondimensional performance parameters are introduced to quantitatively evaluate the heat transfer, flow resistance, and overall performance of multiple HX designs.
The actual velocity and temperature distributions within the full-scale two-fluid HX are obtained through thermofluid simulations based on the governing equations (Eq.~\eqref{eq:NS}).  
Using these simulation results, the velocity, temperature, and pressure distributions at the inlet and outlet are integrated to evaluate the heat transfer rate $Q$ and the pressure drop $\Delta P$ of the HX according to the following equations:
\begin{equation} \label{eq:heat_transfer}
\begin{aligned} 
    Q_i &= \left| \int_{\Gamma_{\mathrm{out},i}} \rho c_p T \mathbf{u} \cdot \mathbf{n} \, dA - \int_{\Gamma_{\mathrm{in},i}} \rho c_p T \mathbf{u} \cdot \mathbf{n} \, dA \right| , \\
    Q &= \frac{Q_1 + Q_2}{2},
\end{aligned}
\end{equation}
\begin{equation} \label{eq:pressure_drop}
\begin{aligned} 
   \Delta P_i & =   \frac{1}{A_{\mathrm{in},i}} \int_{\Gamma_{\mathrm{in},i}} p \, dA - \frac{1}{A_{\mathrm{out},i}} \int_{\Gamma_{\mathrm{out},i}} p \, dA, \\
    \Delta P &= \frac{\Delta P_1 + \Delta P_2}{2}.
\end{aligned}
\end{equation}
Based on these metrics, dimensionless performance parameters are constructed.

First, the Reynolds number $\mathrm{Re}$ is introduced as a dimensionless parameter representing the flow state:
\begin{equation}
    \mathrm{Re} = \frac{\rho u_{\mathrm{in}} D_\mathrm{h}}{\mu},
\end{equation}
where $ u_{\mathrm{in}} $ is the inlet velocity and $D_\mathrm{h} $ is the hydraulic diameter.
The hydraulic diameter $D_\mathrm{h} $ is calculated as follows \cite{Iyer2022,Barakat2024,Barakat2024new,Wang2025}:
\begin{equation}
   D_\mathrm{h} = \frac{4 V_{\mathrm{f,total}}}{A_{\mathrm{total}}},
\end{equation}
where $ V_{\mathrm{f,total}} $ is the total fluid volume in the core region of the HX and $ A_{\mathrm{total}} $ is the total heat exchange surface area.

The Fanning friction factor $ f $, representing the flow resistance due to viscous effects, is calculated from the pressure drop as follows:
\begin{equation}
    f = \frac{D_\mathrm{h} \Delta P}{2\ \rho u_{\mathrm{in}}^2 L_{\mathrm{HX}}},
\end{equation}
where $ L_{\mathrm{HX}} $ is the length of the HX along the flow direction.

For the heat transfer performance, the overall heat transfer coefficient $ U $ is calculated using the logarithmic mean temperature difference $ \Delta T_{\mathrm{LMTD}} $:
\begin{equation}
    \Delta T_{\mathrm{LMTD}} = \frac{(T_{h,\mathrm{in}} - T_{c,\mathrm{out}}) - (T_{h,\mathrm{out}} - T_{c,\mathrm{in}})}{\ln \left( \frac{T_{h,\mathrm{in}} - T_{c,\mathrm{out}}}{T_{h,\mathrm{out}} - T_{c,\mathrm{in}}} \right)},
\end{equation} 
\begin{equation} \label{eq:U}
    U = \frac{Q}{A_{\mathrm{total}}  \Delta T_{\mathrm{LMTD}}}.
\end{equation}
As a dimensionless parameter characterizing convective heat transfer, the Nusselt number $ \mathrm{Nu} $ is introduced:
\begin{equation}
    \mathrm{Nu} = \frac{U D_\mathrm{h}}{k_{\mathrm{f}}},
\end{equation}
where $ k_{\mathrm{f}} $ is the thermal conductivity of the fluid.
Furthermore, to evaluate the convective heat transfer performance normalized by flow conditions and fluid properties, the Colburn $ j $ factor is introduced:
\begin{equation}
    j = \frac{\mathrm{Nu}}{\mathrm{Re} \, \mathrm{Pr}^{1/3}},
\end{equation}
where $ \mathrm{Pr} $ is the Prandtl number.
It provides a standardized basis for comparing the heat transfer performance of different TPMS HXs.

In this study, we introduce the Performance Evaluation Criterion (PEC) \cite{Webb1972,Webb1981,Manglik1995} to evaluate the overall thermal-hydraulic performance of TPMS HXs, considering both heat transfer and flow resistance:
\begin{equation}
    \mathrm{PEC} = \frac{j/j_0}{f/f_0},
\end{equation}
where $j_0$ and $f_0$ are the reference values of the Colburn $j$ factor and the Fanning friction factor $f$, respectively. 
PEC provides a single nondimensional metric that balances heat transfer enhancement against pressure drop, enabling quantitative comparison among TPMS HX designs. 
We use PEC to identify the most efficient design by accounting for the trade-off between heat transfer rate and pressure drop.

\section{Optimization Problem}
\label{sec3}

This section describes the optimization method for TPMS HXs.
First, we present the method for deriving the optimized thickness distribution of TPMS structures by using TO.
Next, the governing equations for optimization and the formulation of the optimization problem are explained.

\subsection{Formulation Basics for Topology Optimization of TPMS Structures}
\label{sec3.1}

TO, originally proposed by Bends{\o}e and Kikuchi \cite{Bendsoe1988}, enables the optimization of material distribution in a fixed design domain under a high degree of design freedom. 
In TO, the material distribution within a fixed design domain $D$ is defined by a characteristic function $\chi_{\Omega} :D\rightarrow \{0,1\}$ representing the presence or absence of material.
This function generally takes discrete values: 0 for void regions and 1 for material regions.
In contrast, TPMS structures are characterized by level set functions, where local wall thickness is controlled by the level set parameter $c$.
However, restricting to binary values, such as in characteristic functions, does not allow for the representation of intermediate thicknesses.

To overcome this limitation, the density method \cite{Bendsoe2003}, which defines design variables as continuous values between 0 and 1, is effective.
By employing a continuous density function to control the wall thickness of TPMS structures, flexible thickness variations, including intermediate values, can be realized.
Furthermore, if the parameter controlling the wall thickness of TPMS structures is defined spatially discretely, it may cause abrupt thickness transitions and geometric discontinuities that degrade performance and complicate manufacturing.
To overcome this issue, applying filtering to the density function is effective.
Spatial smoothing of the design variables through filtering enables continuous and smooth control of the wall thickness distribution.
This filtering approach enables the realization of graded TPMS structures that are both physically feasible and manufacturable.

In this study, a convolutional filter \cite{Bourdin2001,Bruns2001} is applied to the density function $\gamma$, generating a spatially smooth design variable field $\hat{\gamma}$. 
This field does not directly define material presence, but serves as a wall thickness control parameter linked to the level set parameter $c$, which governs the local wall thickness of the TPMS structure.
Here, the design variable field and the level set parameter distribution $c$ are related by the following relationship:
\begin{equation}
\hat{\gamma} = \frac{c(\mathbf{x})-c_{\text{min}}}{c_{\text{max}}-c_{\text{min}}},
\end{equation}
where $c_{\text{min}}$ and $c_{\text{max}}$ denote the minimum and maximum values of $c$, respectively.

In this study, the design variable distribution $\hat{\gamma}$ is optimized based on the density method, and effective properties are defined as functions of the design variable $\hat{\gamma}$.
This framework enables the optimization of the wall thickness distribution in TPMS two-fluid HXs.
Gradients of the objective function with respect to the original design variable $\gamma$ are obtained by applying the chain rule.

\subsection{Governing equations incorporating design variable field}
\label{sec3.2}

In this study, the effective porous media model proposed in Section \ref{sec2.2} is adopted as the physical model for optimization.
The governing equations for velocity field are given by:
\begin{equation} \label{eq:flow}
\begin{aligned}
    &\nabla \cdot \mathbf{U}_i = 0, \\
    &\frac{\rho}{\varepsilon(\hat{\gamma})} \left( \mathbf{U}_i \cdot \nabla \right) \frac{\mathbf{U}_i}{\varepsilon(\hat{\gamma})} 
    = \nabla \cdot \left[ -p_i \mathbf{I} + \mathbf{K} \right] - \alpha(\hat{\gamma}) \mathbf{U}_i - \beta(\hat{\gamma}) |\mathbf{U}_i| \mathbf{U}_i, \\
    &\mathbf{K} = \frac{\mu}{\varepsilon(\hat{\gamma})} \left( \nabla \mathbf{U}_i + (\nabla \mathbf{U}_i)^{\mathrm{T}} \right) - \frac{2}{3} \frac{\mu}{\varepsilon(\hat{\gamma})} (\nabla \cdot \mathbf{U}_i) \mathbf{I}.
\end{aligned}
\end{equation}
The governing equations for temperature field  are given by:
\begin{equation} \label{eq:temperature}
\begin{aligned}
    &\rho c_p \mathbf{U}_i \cdot \nabla T_i + \nabla \cdot \left( - k_{\mathrm{f}}^*(\hat{\gamma}) \nabla T_i \right) = Q_i^*, \\
    &\nabla \cdot \left( - k_{\mathrm{s}}^*(\hat{\gamma}) \nabla T_{\mathrm{w}} \right) = -Q_1^* - Q_2^*, \\
    &Q_i^* = h^*_i(\hat{\gamma},\mathbf{U}_i) \left( T_{\mathrm{w}} - T_i \right) A_{\mathrm{RVE}} (\hat{\gamma}) /V_{\mathrm{RVE}}.
\end{aligned}
\end{equation}
Effective properties are defined as functions of the design variable $\hat{\gamma}$. 
In particular, the effective heat transfer coefficient $h^*_i$ is modeled as a function of both the design variable $\hat{\gamma}$ and the Darcy velocity $\mathbf{U}_i$. 
These effective properties are expressed as polynomial interpolation functions.

%%%%%%%%%%%%%%%%%%%%%%%%%%%%%%%%%%%%%%%%%
\begin{figure}[htb]
    \begin{center}
        \includegraphics[width=0.9\linewidth]{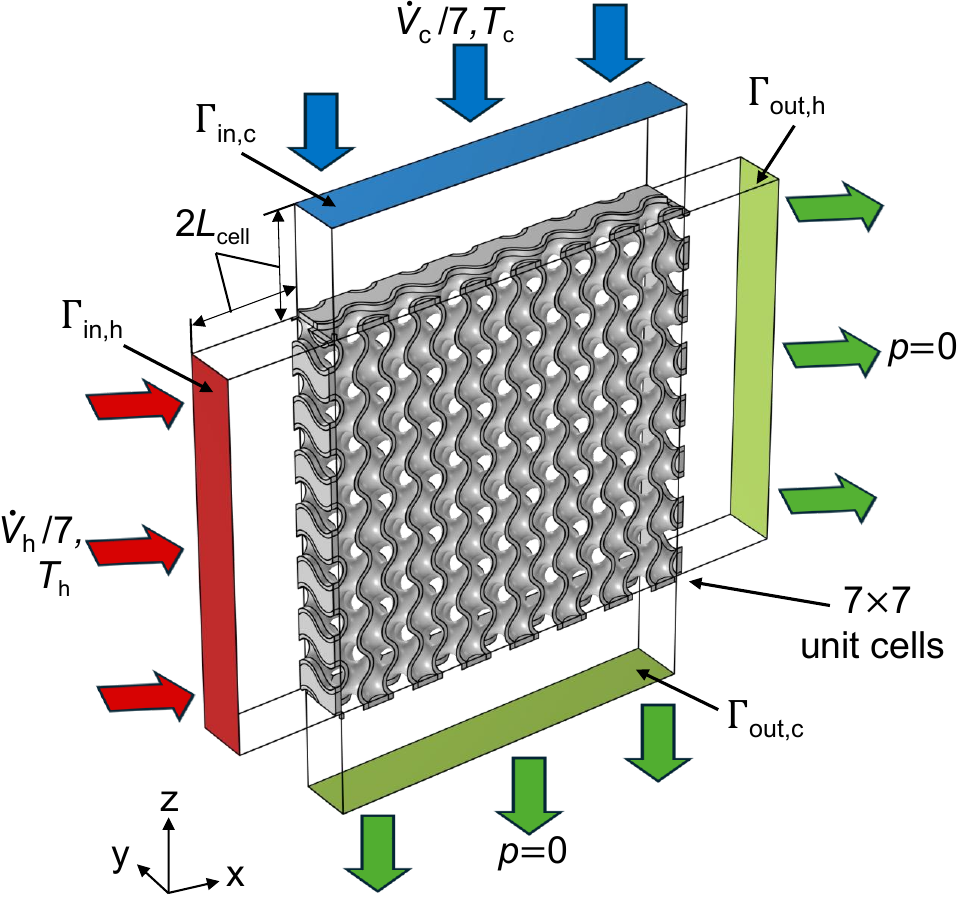}
        \caption{Computational domain and boundary conditions employed for numerical model validation}
        %\vspace{-0.3cm}
        \label{fig4}
    \end{center}%
    %\vspace{-0.2cm}
\end{figure}
%%%%%%%%%%%%%%%%%%%%%%%%%%%%%%%%%%
\subsection{Problem Formulation}
\label{sec3.3}

In HX design, it is crucial to consider both heat transfer performance and pressure drop.
In this study, an objective function that incorporates these two factors with a weighting factor $w$ is employed.
The optimization problem is formulated as follows:
\begin{equation} \label{eq:optimization_problem}
\begin{aligned}
     \text{Find}                              & \quad\gamma                                                           \\ 
        \text{minimize} & \quad J= -Q_{\text{ave}} + w \Delta p_{\text{ave}}        \\ 
        \text{subject to}                        & \quad 0 \leq \gamma(\mathbf{x}) \leq 1,\quad \forall \mathbf{x} \in D 
\end{aligned}
\end{equation}
The average heat transfer rate \( Q_{\text{ave}} \) is defined as the average of the absolute values of the heat transfer amounts calculated from the temperature differences between the inlet and outlet of fluids 1 and 2:
\begin{equation} 
\begin{aligned}
Q_{\text{ave}} = \frac{\rho c_p}{2} \Bigg( & \left| \int_{\Gamma_{\mathrm{out},1}} T_1 \mathbf{U}_1 \cdot \mathbf{n} \, dA 
- \int_{\Gamma_{\mathrm{in},1}} T_1 \mathbf{U}_1 \cdot \mathbf{n} \, dA \right| \\
& + \left| \int_{\Gamma_{\mathrm{out},2}} T_2 \mathbf{U}_2 \cdot \mathbf{n} \, dA 
- \int_{\Gamma_{\mathrm{in},2}} T_2 \mathbf{U}_2 \cdot \mathbf{n} \, dA \right| \Bigg),
\end{aligned}
\end{equation}
where \(\Gamma_{\mathrm{in},i}\) and \(\Gamma_{\mathrm{out},i}\) denote the inlet and outlet surfaces of fluid \(i\), respectively.
The average pressure drop \( \Delta p_{\text{ave}} \) is defined as the average of the pressure differences between the inlet and outlet of fluids 1 and 2:
\begin{equation}
\begin{aligned}
\Delta p_{\text{ave}} = \frac{1}{2} \Bigg( & \frac{1}{A_{\mathrm{in},1}} \int_{\Gamma_{\mathrm{in},1}} p_1 \, dA - \frac{1}{A_{\mathrm{out},1}} \int_{\Gamma_{\mathrm{out},1}} p_1 \, dA \\
& + \frac{1}{A_{\mathrm{in},2}} \int_{\Gamma_{\mathrm{in},2}} p_2 \, dA - \frac{1}{A_{\mathrm{out},2}} \int_{\Gamma_{\mathrm{out},2}} p_2 \, dA \Bigg).
\end{aligned}
\end{equation}
This study performs optimization with different weighting factors $w$ to obtain multiple design solutions reflecting the trade-off between heat transfer performance and pressure drop.

\section{Numerical Examples} 
\label{sec4}

This section demonstrates the effectiveness of the proposed approach through numerical examples.
First, the accuracy of the numerical solver is validated by comparing its results with experimental data.
Subsequently, the effective properties are calculated using a gyroid unit cell.
The accuracy of the effective porous media model is then assessed by comparing the pressure drop and heat transfer rate between the effective and full-scale models.
Next, optimization of wall thicknesses for TPMS-based two-fluid HXs is performed using the effective model.
The resulting optimized structure is dehomogenized, and full-scale numerical simulations are conducted to demonstrate the effectiveness of the proposed method.
All numerical analyses are performed using the finite element software COMSOL Multiphysics (version 6.2), and the optimization is carried out in MATLAB (version 2022b).

% ======================================================================
\begin{table}[t]
    \begin{center}
            \caption{Material properties used for validation}
        \label{table1}
        \scalebox{1}{
            \begin{tabular}{lll}
                \hline
                Parameter                            & Solid       & Fluid                   \\
                \hline
                Thermal conductivity   [$\mathrm{W}/(\mathrm{m}\,\mathrm{K})$]          & 0.18      & 0.631   \\
                Specific heat capacity  [$\mathrm{J}/(\mathrm{kg}\,\mathrm{K})$]      & 2450  & 4178    \\
                Density   [$\mathrm{kg}/\mathrm{m}^3 $] & 1181   &  992.2 \\
                Viscosity  [$\mathrm{Pa}\,\mathrm{s}$]   & - & $6.53 \times 10^{-4}$\\
                \hline
            \end{tabular}
        }
    \end{center}
    \vspace{-0.3cm}
\end{table}
% ======================================================================
\subsection{Validation of Numerical Simulations}
\label{sec4.1}

To ensure the reliability of the numerical simulation results, a validation study was conducted based on a reference experiment. 
Since this study focuses on laminar flow conditions and a gyroid structure, the work by Dixit et al.~\cite{Dixit2022} was selected for validation, as it reports detailed experimental data on a gyroid two-fluid HX under such conditions.
In their study, the thermal performance of a HX composed of a gyroid structure with a cell size of 4.6 mm and a porosity of 0.8 was experimentally evaluated. 
The geometry consisted of seven repeated unit cells in all three spatial directions.

% ======================================================================
\begin{table*}[htb]
    \begin{center}
            \caption{Mesh dependency study results}
        \label{table2}
        \scalebox{1}{
            \begin{tabular}{llllll}
                \hline
                Mesh No.                            & No. of meshes       & Heat transfer rate [W]    & Error [\%]    & Pressure loss [Pa]    & Error [\%]         \\
                \hline
                Mesh 1     & 3,166,058      & 142.02  & 3.50   & 56.77 &  1.88       \\
                Mesh 2     & 4,810,819  & 139.93     & 1.98   & 56.10  & 0.68       \\
                Mesh 3      & 7,130,430 & 138.66 & 1.06    & 55.91   & 0.34      \\
                Mesh 4      & 10,865,567 & 137.59   & 0.28   & 55.79  & 0.13       \\
                Mesh 5   & 14,462,935 & 137.24   & 0.02   & 55.73  & 0.02      \\
                Mesh 6        & 18,356,943  & 137.21  & standard   & 55.72 & standard     \\
                \hline
            \end{tabular}
        }
    \end{center}
    \vspace{-0.2cm}
\end{table*}
% =====================================================================
%%%%%%%%%%%%%%%%%%%%%%%%%%%%%%%%%%%%%%%%%
\begin{figure*}[h!]
    \centering
    \begin{minipage}[t]{0.48\textwidth}
        \centering
        \includegraphics[width=0.95\linewidth]{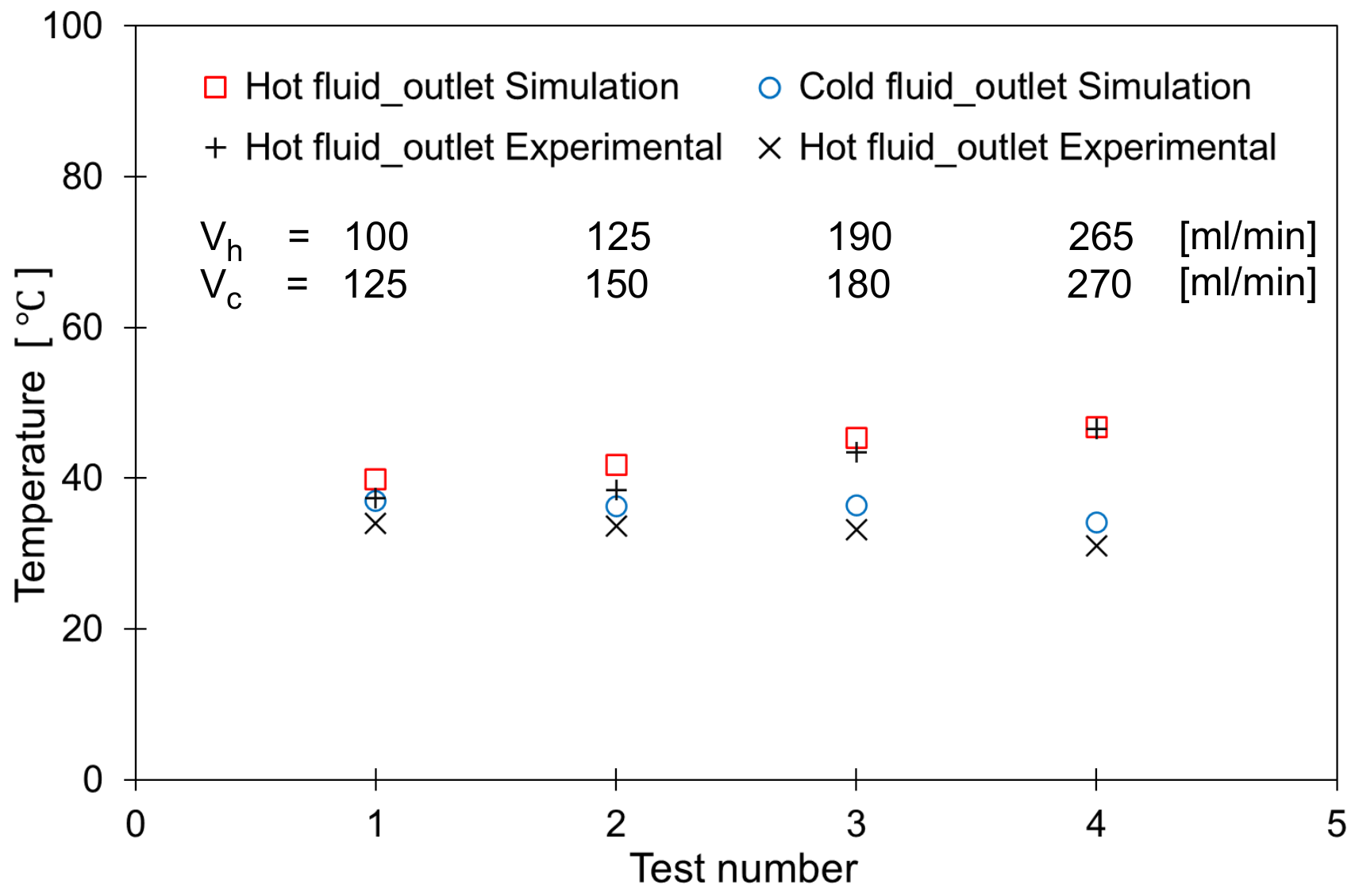}
        \subcaption{Outlet temperatures}
        \label{fig5a}
    \end{minipage} 
    \begin{minipage}[t]{0.48\textwidth}
        \centering
        \includegraphics[width=0.95\linewidth]{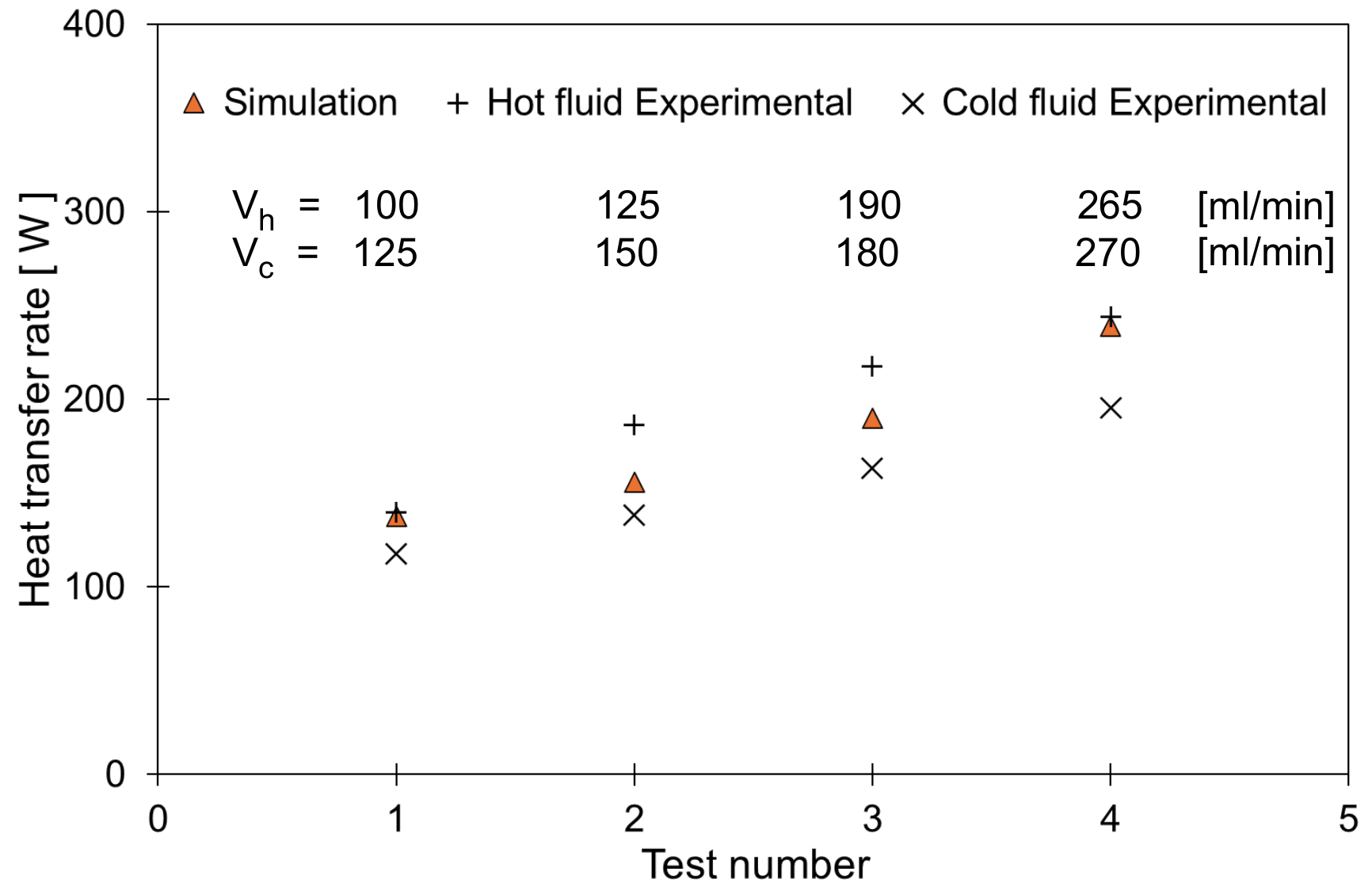}
        \subcaption{Heat transfer rates}
        \label{fig5b}
    \end{minipage}
    %\vspace{-0.2cm}
    \caption{Comparison of numerical and experimental results for validation}
    \vspace{-0.2cm}
    \label{fig5}
\end{figure*}
%%%%%%%%%%%%%%%%%%%%%%%%%%%%%%%%%%%%%%%%%
% ======================================================================
\begin{table}[h]
    \begin{center}
            \caption{Material properties used for optimization}
        \label{table3}
        \scalebox{1}{
            \begin{tabular}{lll}
                \hline
                Parameter                            & Solid       & Fluid                   \\
                \hline
                Thermal conductivity   [$\mathrm{W}/(\mathrm{m}\,\mathrm{K})$]          & 237     & 0.631   \\
                Specific heat capacity  [$\mathrm{J}/(\mathrm{kg}\,\mathrm{K})$]      & 900  & 4178    \\
                Density   [$\mathrm{kg}/\mathrm{m}^3 $] & 2700  &  992.2 \\
                Viscosity  [$\mathrm{Pa}\,\mathrm{s}$]   & - & $6.53 \times 10^{-4}$\\
                \hline
            \end{tabular}
        }
    \end{center}
    \vspace{-0.3cm}
\end{table}
% ======================================================================
In contrast, to reduce computational cost, the numerical model used for validation was simplified to a single unit cell in the y-direction, as shown in Fig.~\ref{fig4}. 
The gyroid geometry was adjusted to achieve a porosity of 0.8 and was constructed using Eq.~\eqref{eq:gyroid} with \(L_\mathrm{cell} = 4.6 \times 10^{-3}\,\mathrm{m}\) and \(c = 1.426 \times 10^{-3}\,\mathrm{m}\).
The solid domain was assumed to be composed of a polymer, and its material properties were referenced from Dixit's article.
Water was used as the working fluid. 
The reference temperature for the fluid properties was defined as the average of the inlet temperatures of the hot and cold streams, which was 313.15 K. The properties at this temperature were used.
The material properties employed in the validation are summarized in Table~\ref{table1}.
To enable a meaningful comparison, the inlet conditions were matched to those used in the Dixit's article.
However, to account for the fact that the thickness of the HX was reduced to one-seventh of the original geometry, the volumetric flow rates used in the simulation were scaled down accordingly to one-seventh of the experimental values.
The inlet temperatures of the hot and cold fluids were set to $T_\mathrm{h} = 333.15$  and $T_\mathrm{c} = 293.15$, respectively.
As boundary conditions, all surfaces except for the inlet and outlet faces were assumed to be adiabatic and subject to no-slip conditions.

Initially, a mesh independence study was conducted using the present simulation model.  
Several meshing schemes were evaluated by varying the minimum cell size and the number of boundary layer elements, and their effects on the pressure drop and heat transfer rate were assessed.  
The study was carried out under flow rate conditions of $\dot{V}_\mathrm{h} = 100$ and $\dot{V}_\mathrm{c} = 125 \ \mathrm{ml/min}$.  
The results are summarized in Table~\ref{table2}.  
As the number of elements increased, both the heat transfer rate and pressure drop showed a tendency to converge to constant values.
In this study, Mesh 4 is selected for the subsequent numerical analyses by considering the balance between computational accuracy and cost, as it results in less than 1\% error in both heat transfer rate and pressure drop.  
With this mesh configuration, the minimum element size and the number of boundary layer mesh layers are $9.6 \times 10^{-5}$ m and 4, respectively.  
All numerical results in this paper, except those based on the effective model, are obtained using this mesh configuration.

Subsequently, validation was conducted by comparing the numerical results with the experimental data reported by Dixit et al.
Fig.~\ref{fig5}(a) presents the comparison of outlet temperatures. 
While some discrepancies were observed, the overall trends reasonably agreed with the experimental data.  
Fig.~\ref{fig5}(b) shows the comparison of heat transfer rates. 
In the experimental results, a difference in heat transfer rate of approximately 22 to 55 W was observed between the hot and cold sides, which is due to environmental heat losses and other unavoidable experimental uncertainties.
The heat transfer rates obtained from the numerical simulation fell within this range, indicating reasonable consistency.
Although certain deviations between the experimental and numerical results were observed, these can be attributed to measurement errors in the experiments, unavoidable heat losses, and simplifications in the numerical model.  
Therefore, the present numerical simulation is considered to provide sufficient accuracy for the gyroid HXs under laminar flow conditions, and is thus applicable to the subsequent analyses.

\subsection{Computation of Effective Properties}
\label{sec4.2}

The effective properties are computed based on the method described in Section~\ref{sec2.3}.
In this study, optimization is conducted for a gyroid two-fluid HX composed of aluminum as the solid material and water as the working fluid. 
The material properties are defined at a reference temperature of 313.15 K, which corresponds to the average of the inlet temperatures of the hot and cold fluids. 
The material properties employed in the optimization are summarized in Table~\ref{table3}.

To derive the interpolation functions for the effective properties, 10 unit cells with varying wall thicknesses are used, corresponding to a unit cell length $L_\text{cell} = 4.6 \times 10^{-3} \, \mathrm{m}$ and level set parameters $c = 1.426,\, 1.75,\, 2.0,\, 2.25,\, 2.5,\, 2.75,\, 3.0,\, 3.25,\, 3.5,\, 3.75 \ (\times 10^{-3})$.
The effective properties of each unit cell are computed and expressed as functions of the design variable $\hat{\gamma}$ through polynomial interpolation.

%%%%%%%%%%%%%%%%%%%%%%%%%%%%%%%%%%%%%%%%%
\begin{figure*}[h]
    \centering
    \begin{minipage}[t]{0.45\textwidth}
        \centering
        \includegraphics[width=0.9\linewidth]{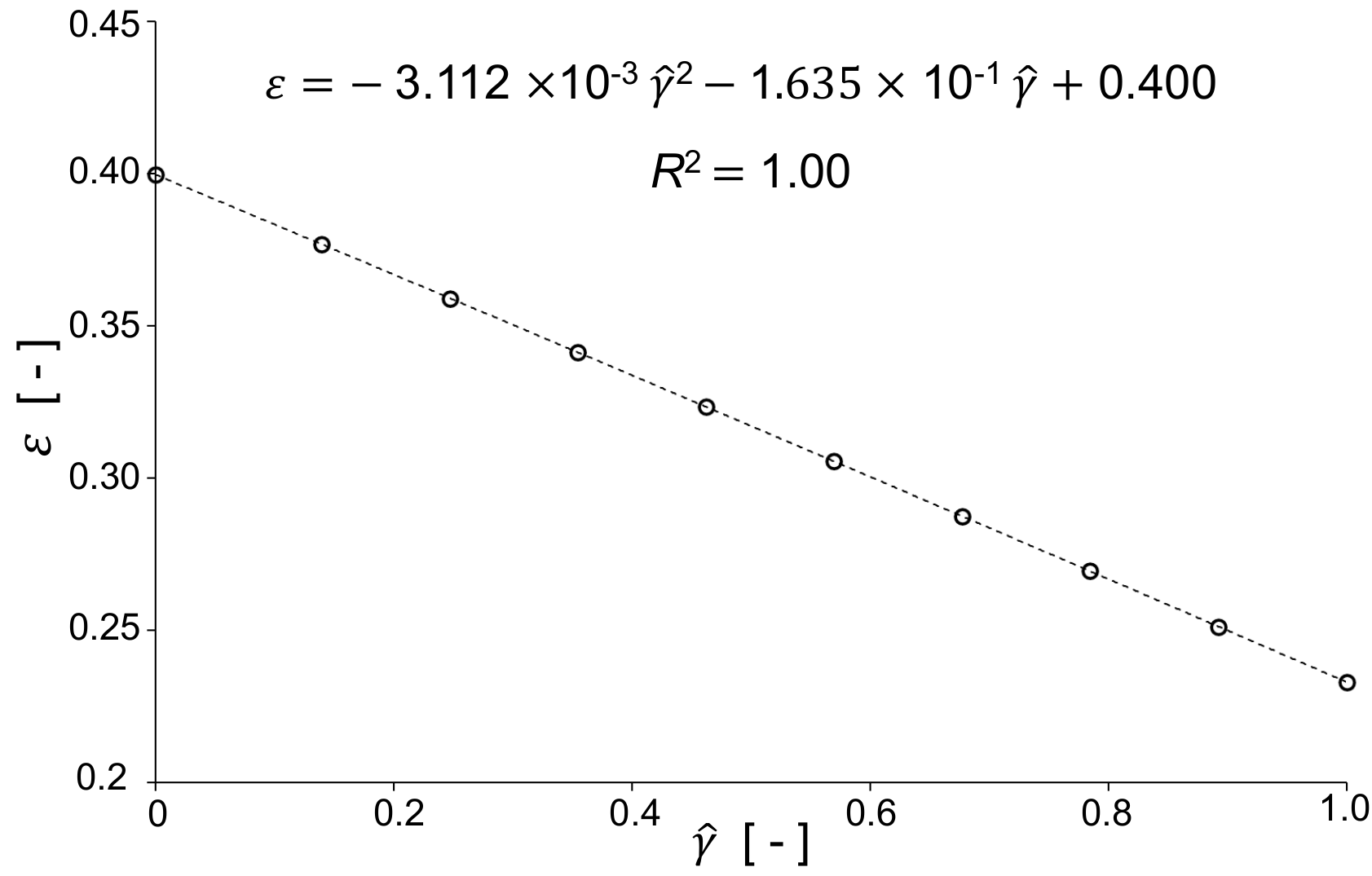}
        \subcaption{Interpolation result of the $\varepsilon$}
        \label{fig6a}
    \end{minipage} 
    \vspace{0.2cm}
    \begin{minipage}[t]{0.45\textwidth}
        \centering
        \includegraphics[width=0.9\linewidth]{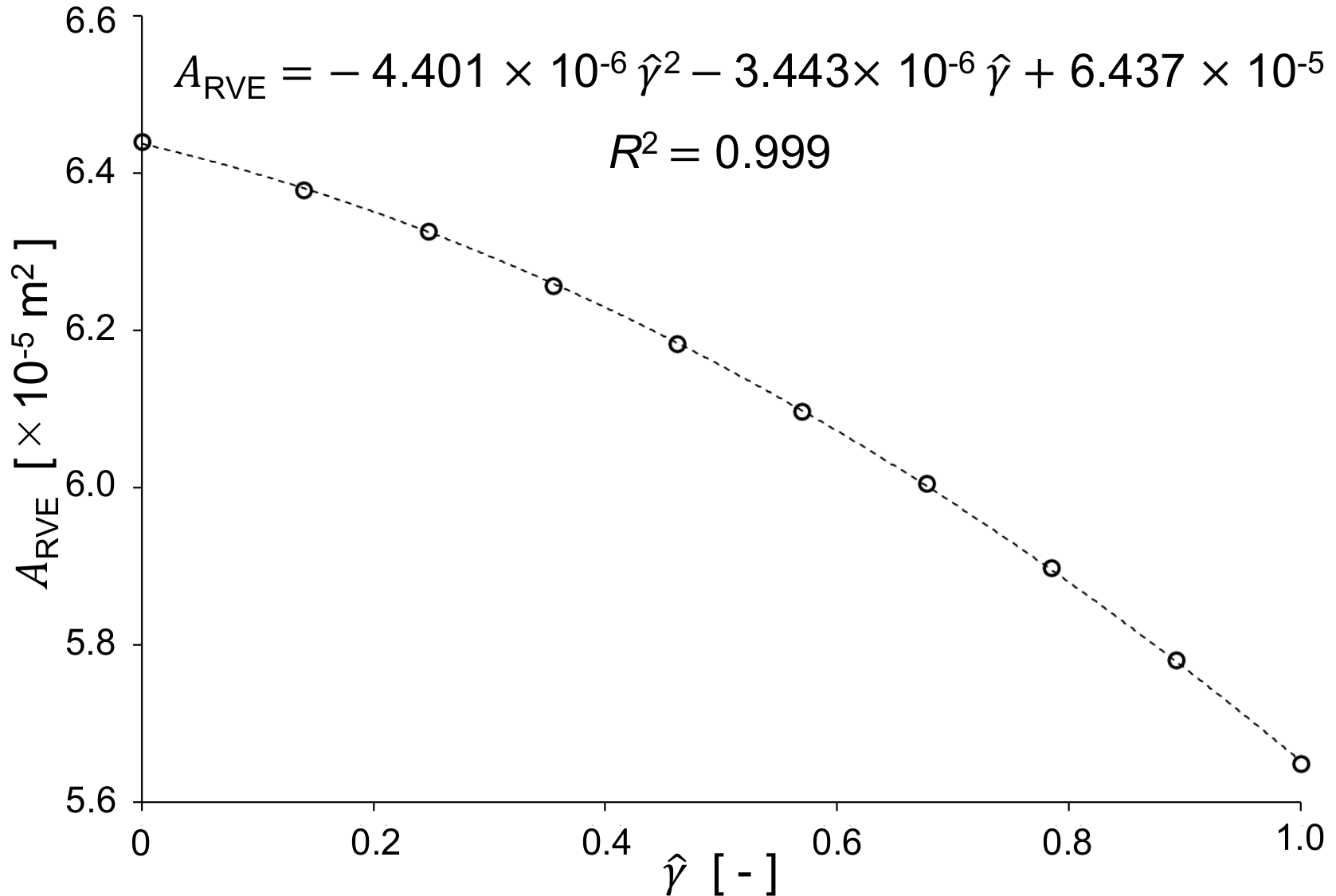}
        \subcaption{Interpolation result of the $A_{\mathrm{RVE}}$}
        \label{fig6b}
    \end{minipage}\\
        \begin{minipage}[t]{0.45\textwidth}
        \centering
        \includegraphics[width=0.9\linewidth]{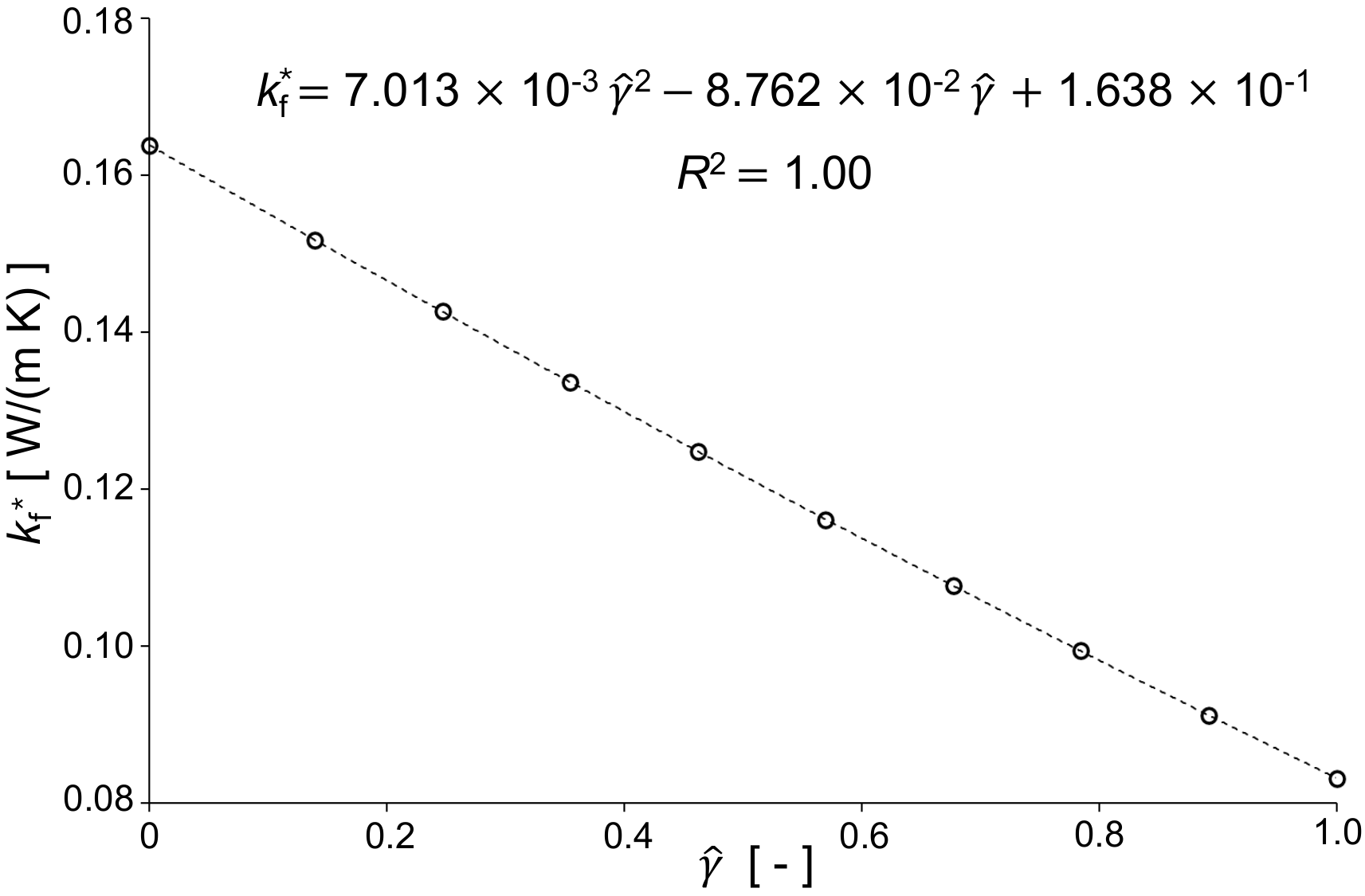}
        \subcaption{Interpolation result of the $k_f^*$}
        \label{fig6c}
    \end{minipage} 
    \vspace{0.2cm}
    \begin{minipage}[t]{0.45\textwidth}
        \centering
        \includegraphics[width=0.9\linewidth]{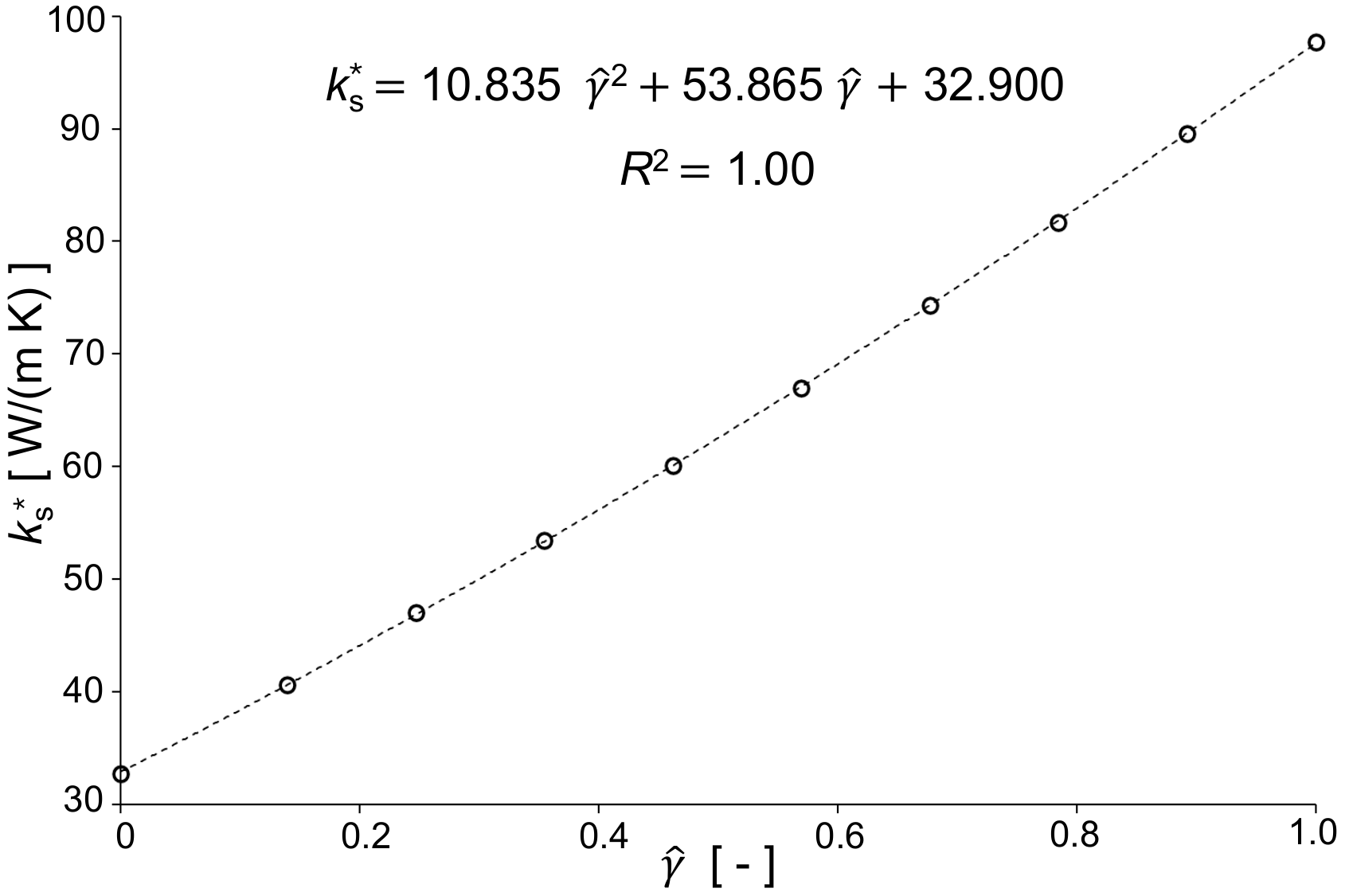}
        \subcaption{Interpolation result of the $k_s^*$}
        \label{fig6d}
    \end{minipage}
    %\vspace{-0.2cm}
    \caption{Interpolation functions of effective properties dependent on design variable}
    \vspace{-0.2cm}
    \label{fig6}
\end{figure*}
\begin{figure*}[htb]
    \centering
    \begin{minipage}[t]{0.45\textwidth}
        \centering
        \includegraphics[width=0.95\linewidth]{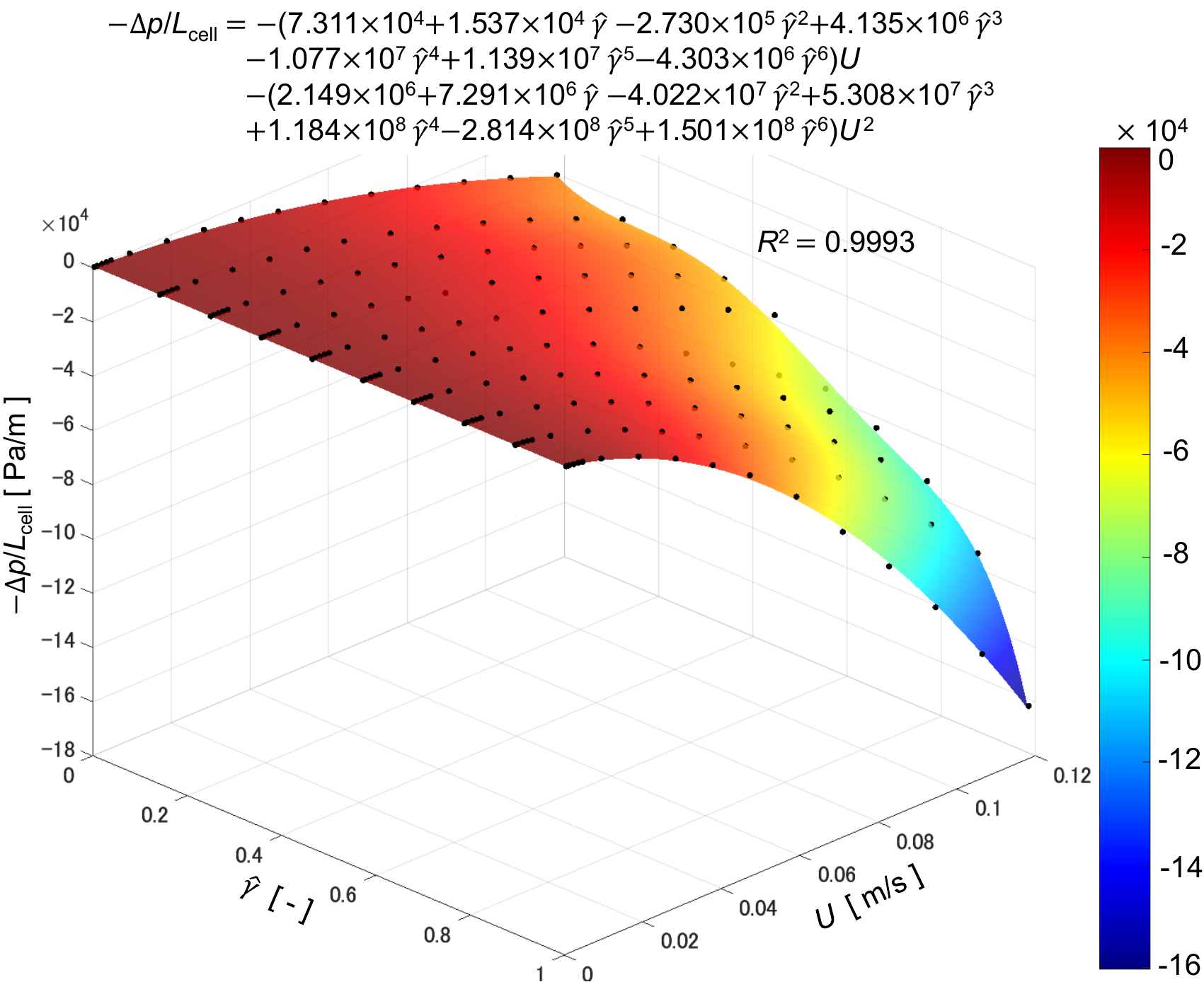}
        \subcaption{Interpolation results of the $\alpha$ and $\beta$}
        \label{fig7a}
    \end{minipage} 
    \begin{minipage}[t]{0.45\textwidth}
        \centering
        \includegraphics[width=0.95\linewidth]{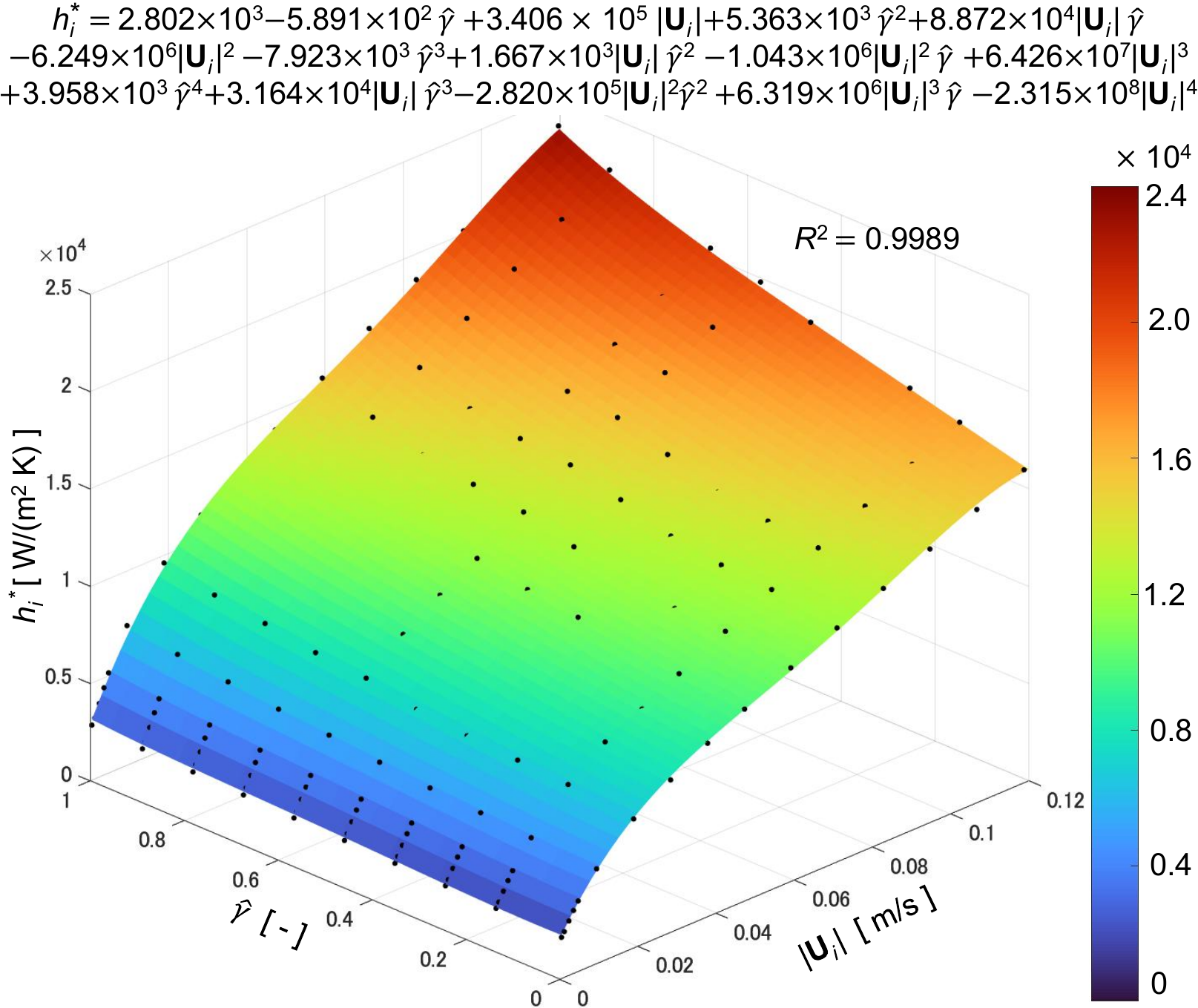}
        \subcaption{Interpolation result of the $h^*_i$}
        \label{fig7b}
    \end{minipage}
    %\vspace{-0.2cm}
    \caption{Interpolation functions of effective properties dependent on design variable and Darcy velocity}
    \vspace{-0.2cm}
    \label{fig7}
\end{figure*}

%%%%%%%%%%%%%%%%%%%%%%%%%%%%%%%%%%%%%%%%%
Figs.~\ref{fig6} and \ref{fig7} present the computed effective properties.
Fig.~\ref{fig6} illustrates the interpolation results for $\varepsilon$, $A_{\mathrm{RVE}}$, $k_f^*$ and $k_s^*$, each fitted with a quadratic function based on the data obtained from individual unit cells.
Fig.~\ref{fig7} shows the interpolation results for $\alpha$, $\beta$, and $h^*$. 
The inlet temperatures of the RVE were set to $T_\mathrm{H} = 314.15\,\mathrm{K}$ and $T_\mathrm{L} = 312.15\,\mathrm{K}$, and the volumetric flow rate $\dot{V}$ was varied from $1.0 \times 10^{-8}$ to $2.5 \times 10^{-6}\,\mathrm{m}^3/\mathrm{s}$.
Since $\alpha$, $\beta$, and $h^*_i$ depend not only on wall thickness but also on inlet flow rate, polynomial interpolation functions were constructed using all available data, covering 10 wall thicknesses and 20 flow rates for each thickness.
The interpolation functions shown in Figs.~\ref{fig7}(a) and \ref{fig7}(b) were derived using MATLAB.  
For Fig.~\ref{fig7}(a), since the pressure drop $\Delta p$ increases with both the design variable $\hat{\gamma}$ and the Darcy velocity $U$, the interpolation function was constructed to satisfy the following conditions:
\begin{equation} 
\frac{\partial (\Delta p/L_\mathrm{cell})}{\partial \hat{\gamma}} > 0, \quad \frac{\partial (\Delta p/L_\mathrm{cell})}{\partial U} > 0.
\end{equation}
Effective properties $\alpha$ and $\beta$ were determined based on the interpolation functions shown in Fig.~\ref{fig7}(a) and the relationship given by Eq.~\eqref{eq:Darcy_Forchheimer}.
The interpolation functions exhibit excellent agreement with the discrete data points, accurately capturing the variations in the effective properties over the entire range of the design variable.

\subsection{Verification of Effective Porous Medium Model}
\label{sec4.3}

\begin{figure}[htb]
    \begin{center}
        \includegraphics[width=0.975\linewidth]{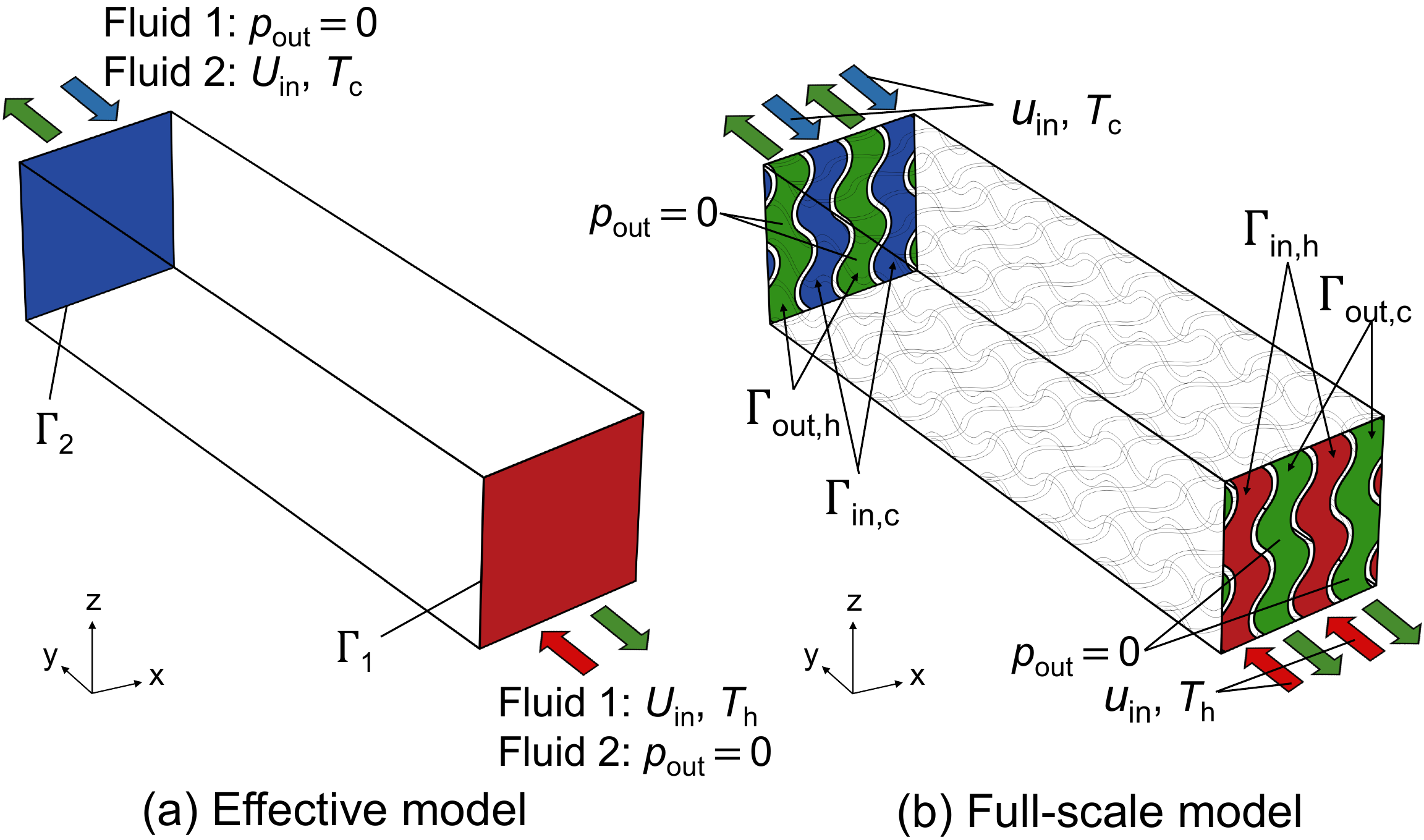}
        \vspace{-0.2cm}
        \caption{Analysis model used for verification of the effective porous media model}
        \label{fig8}
    \end{center}%
    %\vspace{-0.2cm}
\end{figure}

We verify the validity of the effective porous media model proposed in Section~\ref{sec2.2} and the interpolation functions of effective properties derived in Section~\ref{sec4.2}.  
Numerical simulations using both the effective porous media model and the full-scale model are performed over the same analysis domain.  
The validity of the constructed optimization model is assessed by comparing the pressure drop and heat transfer rate obtained from both models.  
The analysis domain used for verification is illustrated in Fig.~\ref{fig8}.  
In this model, a uniform distribution of the design variable $\hat{\gamma}$ is applied in the effective model, whereas the full-scale model is configured with a gyroid structure having a uniform wall thickness corresponding to the given design variable.  
Counter-flow conditions are assumed, with an inlet velocity of \(u_{\mathrm{in}} = 0.03\,\mathrm{m/s}\), a hot fluid inlet temperature of \(T_{\mathrm{h}} = 333.15\,\mathrm{K}\), and a cold fluid inlet temperature of \(T_{\mathrm{c}} = 293.15\,\mathrm{K}\).  
No-slip and adiabatic boundary conditions are applied to all surfaces except the inlet and outlet.  
In the full-scale model, fully developed flow conditions are imposed at the inlet and outlet boundaries.  
In the effective model, the Darcy velocity at the inlet is adjusted to match the volumetric flow rate of the full-scale model.

%%%%%%%%%%%%%%%%%%%%%%%%%%%%%%%%%%%%%%%%%
%%%%%%%%%%%%%%%%%%%%%%%%%%%%%%%%%%%%%%%%%
%%%%%%%%%%%%%%%%%%%%%%%%%%%%%%%%%%%%%%%%%
%%%%%%%%%%%%%%%%%%%%%%%%%%%%%%%%%%%%%%%%%
%%%%%%%%%%%%%%%%%%%%%%%%%%%%%%%%%%%%%%%%%
\begin{figure}[htb]
    \begin{center}
        \includegraphics[width=0.95\linewidth]{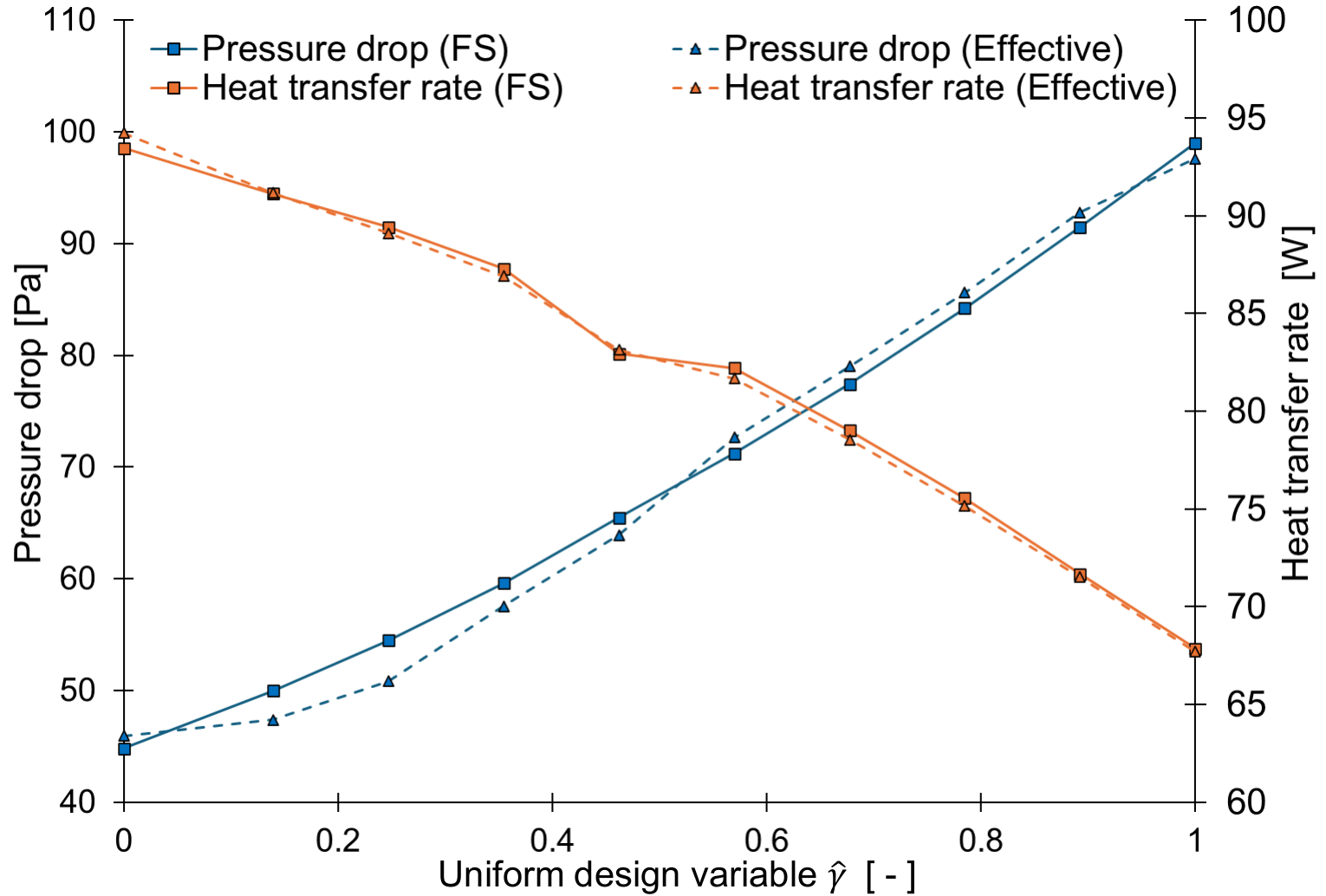}
        \caption{Comparison of simulation results between the effective porous media model and full-scale model}
        %\vspace{-0.3cm}
        \label{fig9}
    \end{center}%
    %\vspace{-0.2cm}
\end{figure}
%%%%%%%%%%%%%%%%%%%%%%%%%%%%%%%%%%%%%%%%%
Fig.~\ref{fig9} compares the simulation results obtained from both models.  
Dashed lines represent the results of the effective model, while solid lines correspond to the full-scale model.  
The maximum relative errors are 6.72\% for pressure drop and 0.83\% for heat transfer rate, demonstrating that the effective model accurately predicts the full-scale simulation.  
These findings confirm the validity of the proposed effective porous media model for TPMS two-fluid HXs and suggest that optimizing gyroid two-fluid HXs with the constructed effective model can improve performance.
The computational aspects of this approach, including the reduction in computational cost achieved by the effective model, the cost of deriving the interpolation functions, and the range of problems for which the constructed model is applicable, are discussed in the following section.

\subsection{Computational Efficiency and Applicability}
\label{sec4.4}

In this section, we discuss the computational cost, efficiency, and applicability of the proposed method.
The computation times reported in this section correspond to analyses performed using a 2.45 GHz AMD EPYC 7763 64-core processor with 512 GB RAM.

First, the computational cost reduction achieved by the effective model is evaluated.
Using the analysis model employed for verification in Section \ref{sec4.2} (Fig.~\ref{fig8}), the computational costs of the full-scale and effective models were compared.
For a uniform thickness distribution of $c=1.426\times10^{-3}$, the full-scale model contained 6,691,417 mesh elements and required 8,115 s for analysis.
In contrast, the effective model, with 256 mesh elements, required only 31 s, making it approximately 260 times faster than the full-scale simulation.
These results demonstrate that the proposed effective model can drastically reduce computational cost while maintaining reasonable accuracy.

%%%%%%%%%%%%%%%%%%%%%%%%%%%%%%%%%%%%%%%%%
\begin{figure*}[htb]
    \begin{center}
        \includegraphics[width=0.85\linewidth]{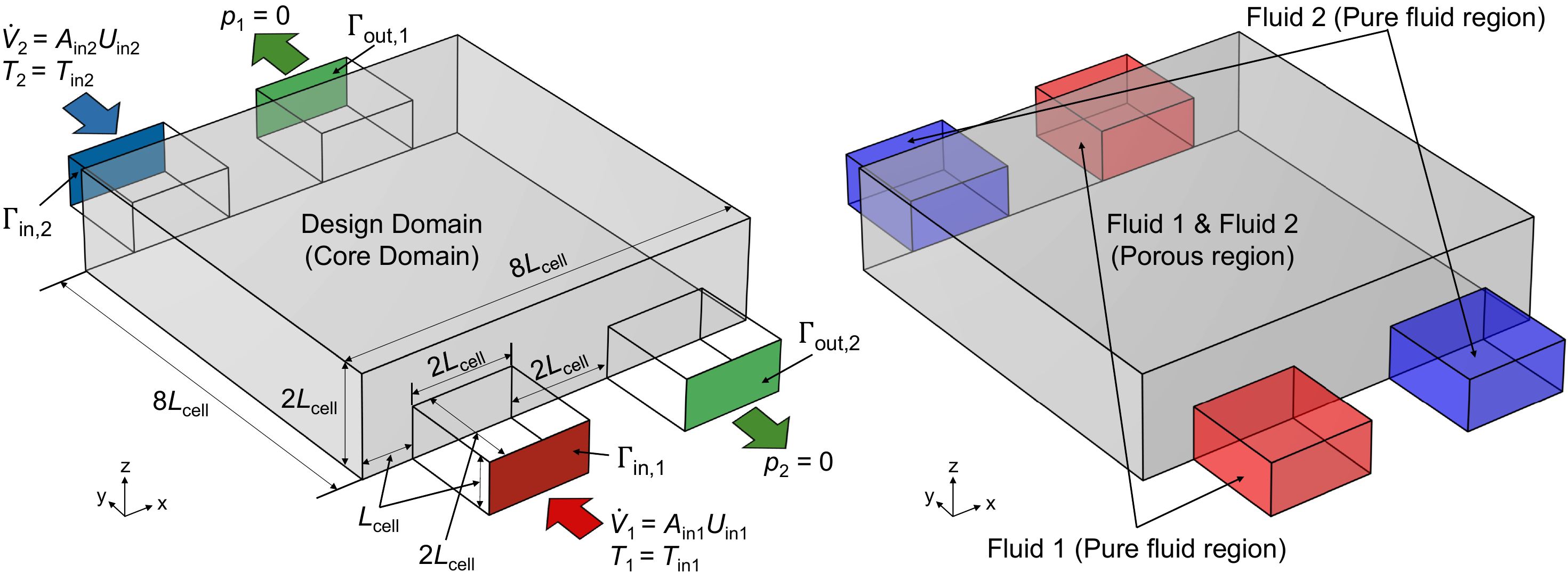}
        \caption{Problem setup for wall-thickness optimization of the gyroid two-fluid HX}
        %\vspace{-0.3cm}
        \label{fig10}
    \end{center}%
    %\vspace{-0.2cm}
\end{figure*}
%%%%%%%%%%%%%%%%%%%%%%%%%%%%%%%%%%%%%%%%%
%%%%%%%%%%%%%%%%%%%%%%%%%%%%%%%%%%%%%%%%%
\begin{figure}[htb]
    \begin{center}
        \includegraphics[width=0.975\linewidth]{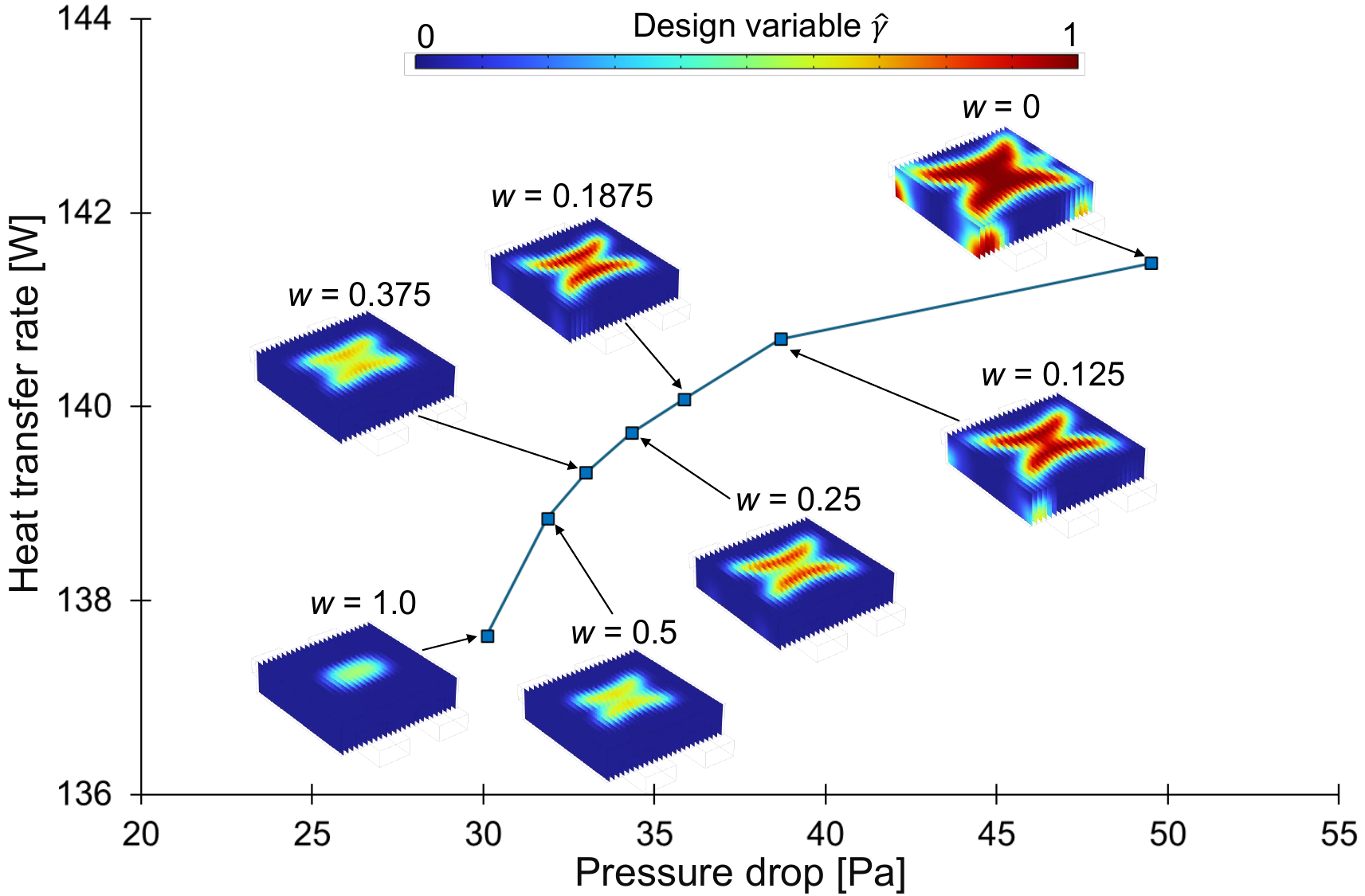}
        \caption{Optimization results for different values of the weighting factor $w$}
        %\vspace{-0.3cm}
        \label{fig11}
    \end{center}%
    %\vspace{-0.2cm}
\end{figure}
%%%%%%%%%%%%%%%%%%%%%%%%%%%%%%%%%%%%%%%%%
Next, we discuss the computational cost required to derive interpolation functions for the effective properties.
In this study, numerical analyses were performed for a sufficient number of RVEs to ensure that the resulting interpolation functions provide adequate predictive accuracy.
Specifically, ten RVEs with different thicknesses were prepared, and each was analyzed under twenty flow rate patterns.
For example, analyzing twenty flow rate patterns for a gyroid unit cell with $c=1.426\times10^{-3}$ and 227,471 mesh elements under the boundary conditions shown in Fig.~\ref{fig3} required 3,286 s.
Consequently, analyzing all ten RVEs to construct the interpolation functions would require approximately nine hours.
Efficient parallel computing can reduce this time.
While reducing the number of thickness variations or flow rate patterns can further lower computational cost, it may compromise the predictive accuracy of the effective model.
Since optimization involves multiple iterative calculations, performing optimization using the full-scale model would require an impractically high computational cost and take an excessively long time.
Considering this, the computational cost for constructing the interpolation functions of the effective properties is reasonable.

Finally, we discuss the applicability of the proposed method. Although the interpolation functions derived in Section \ref{sec4.2} are limited to the conditions considered in this study, the proposed method can be extended to other TPMS structures, materials, and unit cell sizes. 
While this study focused on counter-flow configurations, the method can also be applied to cross-flow or parallel-flow problems by appropriately modifying the RVE boundary conditions. Consequently, the proposed method can be widely applied to HXs with various TPMS types and flow configurations, providing a foundation for developing simplified models and optimization frameworks.

\subsection{Optimization Results}
\label{sec4.5}

Optimization of wall thicknesses distribution for the gyroid two-fluid HXs was conducted using the developed optimization model.  
The thickness distribution was determined by optimizing the spatial distribution of the design variable $\hat{\gamma}$.
The design variable was updated using the Method of Moving Asymptotes (MMA) \cite{Svanberg1987}.

In this study, a counterflow configuration was assumed, and the optimization was performed based on the problem setup shown in Fig.~\ref{fig10}.  
Fluid 1 was designated as the hot fluid and Fluid 2 as the cold fluid.  
The inlet temperatures were set to $T_{\mathrm{in}1} = 333.15~\mathrm{K}$ and $T_{\mathrm{in}2} = 293.15~\mathrm{K}$, respectively.  
The inlet flow rates $ \dot{V}_i $ for fluids 1 and 2 were set as boundary conditions.
These flow rates are given by $ \dot{V}_i = 2 L_{\mathrm{cell}}^2 U_{\mathrm{in}i}$, where $ L_{\mathrm{cell}} $ is the length of the unit cell, which is $ 4.6 \times 10^{-3}~\mathrm{m} $, and $ U_{\mathrm{in}i} $ denotes the average inlet velocity of fluid $ i $.
In this study, the flow rates $ \dot{V}_1 $ and $ \dot{V}_2 $ corresponding to $ U_{\mathrm{in}1} = U_{\mathrm{in}2} = 0.03~\mathrm{m/s} $.
Adiabatic and no-slip boundary conditions were applied to all surfaces except the inlet and outlet.  
The gyroid structures were embedded only within the core domain of the HX, while the other domains were treated as pure fluid regions. 
To obtain multiple optimized solutions, the optimization was performed multiple times with different values of the weighting factor $w$.
As the initial configuration, the design variable $\hat{\gamma}$ was uniformly set to 0.5, and the optimization was performed over 50 iterations.  
Fig.~\ref{fig11} presents the optimized results corresponding to each value of $w$.  
All values shown in Fig.~\ref{fig11} were evaluated using the effective porous media model.  
As the value of the weighting factor $w$ increased, the regions with larger thicknesses decreased, leading to a corresponding reduction in pressure drop. 

%%%%%%%%%%%%%%%%%%%%%%%%%%%%%%%%%%%%%%%%%
\begin{figure*}[h]
    \centering
    \begin{minipage}[t]{0.45\textwidth}
        \centering
        \includegraphics[width=0.9\linewidth]{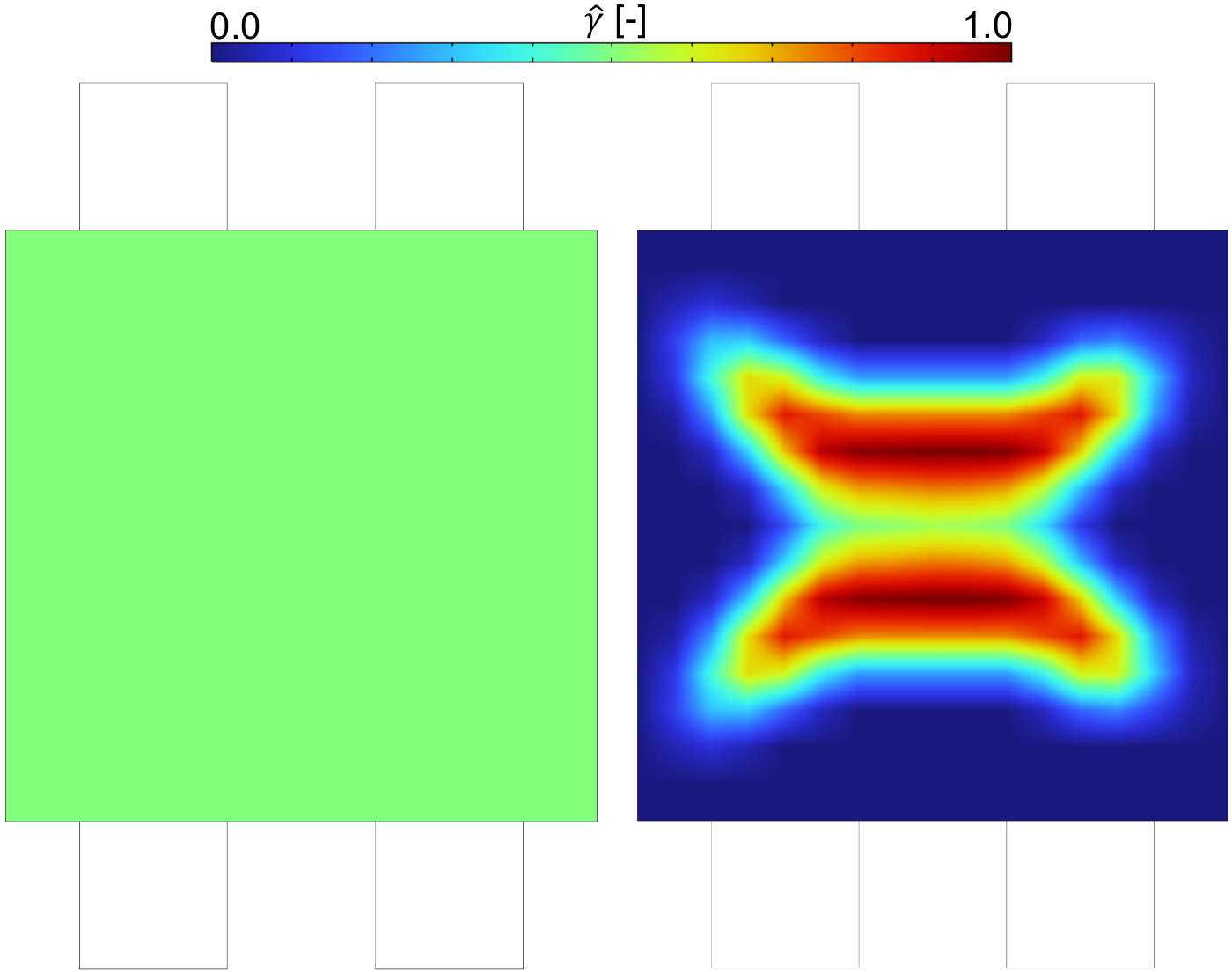}
        \subcaption{Design variable distribution}
        \label{fig12a}
    \end{minipage} 
    \vspace{0.2cm}
    \begin{minipage}[t]{0.45\textwidth}
        \centering
        \includegraphics[width=0.9\linewidth]{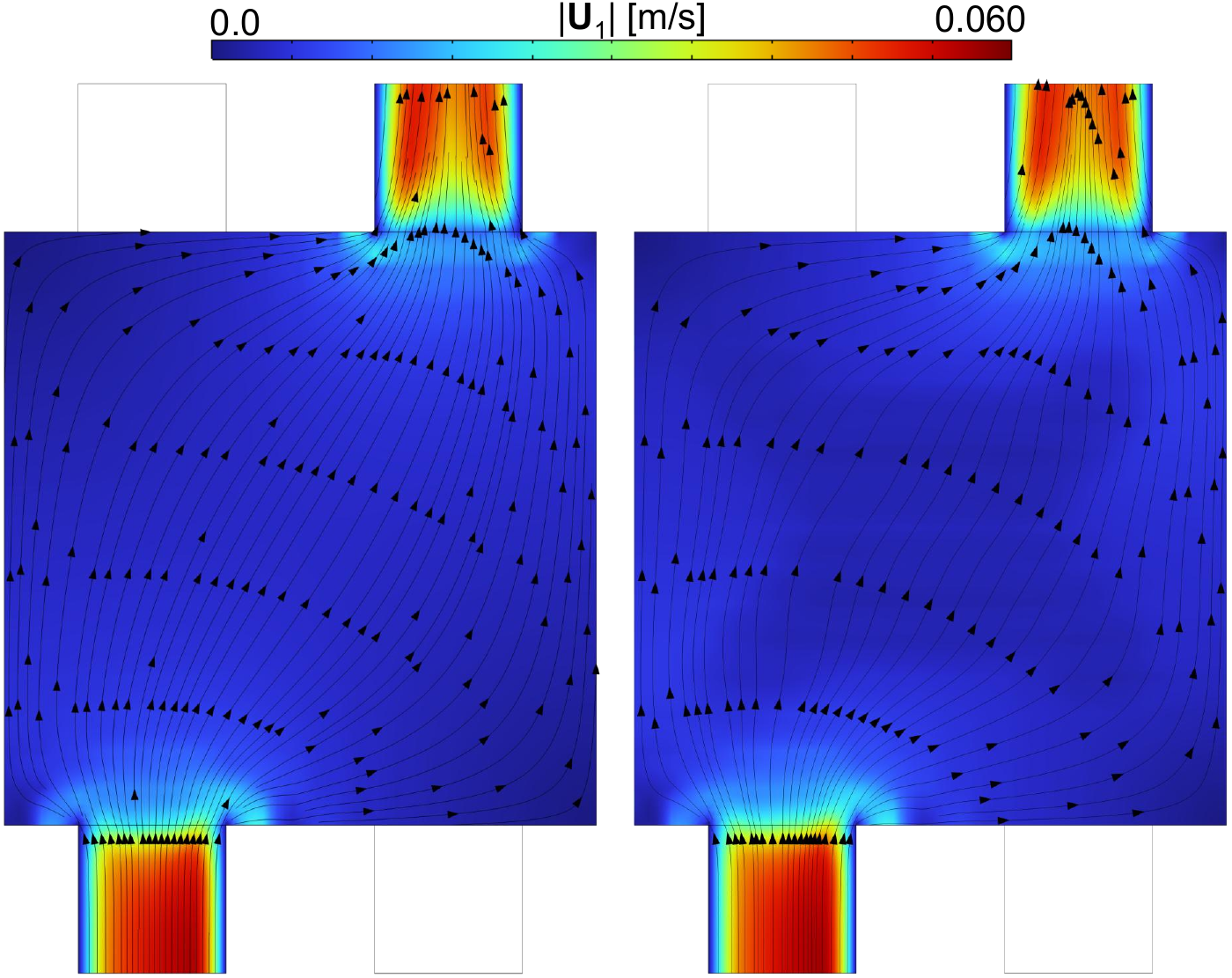}
        \subcaption{Darcy velocity distribution}
        \label{fig12b}
    \end{minipage} \\
    \vspace{0.2cm}
    \begin{minipage}[t]{0.45\textwidth}
        \centering
        \includegraphics[width=0.9\linewidth]{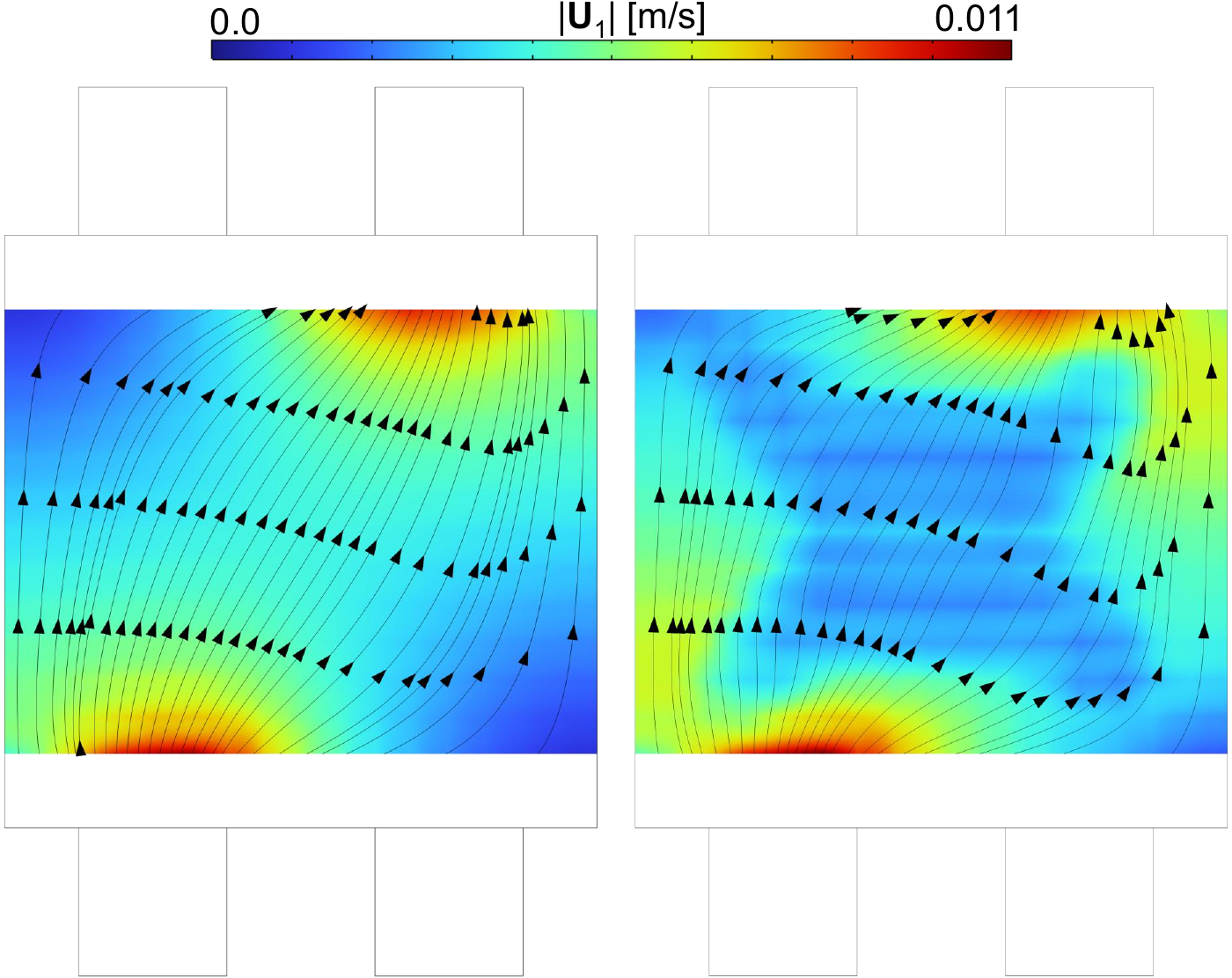}
        \subcaption{Darcy velocity distribution within the region \( L_{\mathrm{cell}} \leq  y \leq  7 L_{\mathrm{cell}} \)}
        \label{fig12c}
    \end{minipage}
    \vspace{0.2cm}
    \begin{minipage}[t]{0.45\textwidth}
        \centering
        \includegraphics[width=0.9\linewidth]{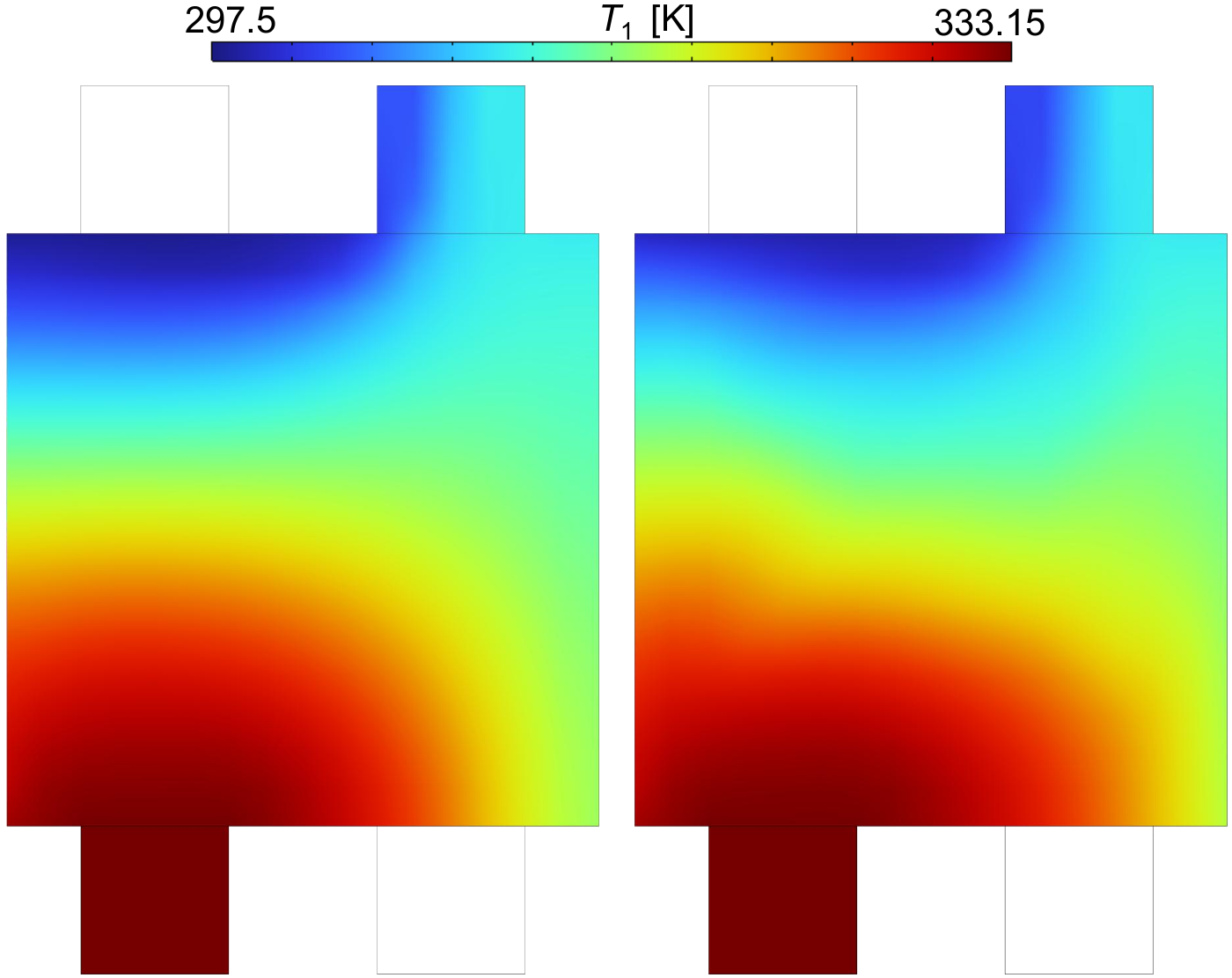}
        \subcaption{Temperature distribution}
        \label{fig12d}
    \end{minipage}
    \caption{Distributions at the cross-section \( z = L_{\mathrm{cell}}/2 \) focusing on the hot fluid for both the initial (left) and optimized (right) designs}
    \vspace{-0.2cm}
    \label{fig12}
\end{figure*}
%%%%%%%%%%%%%%%%%%%%%%%%%%%%%%%%%%%%%%%%%
Fig.~\ref{fig12} illustrates the distributions of the design variable, the Darcy velocity of the hot fluid, and the temperature distributions of the hot fluid at the cross-section $z = L_{\mathrm{cell}}/2$ for both the initial configuration ($\hat{\gamma} = 0.5$) and the optimized structure ($w = 0.25$).
In each panel, the left half corresponds to the uniform-thickness design, and the right half to the optimized-thickness design.
Fig.~\ref{fig12}(b) shows the overall Darcy velocity distribution, while Fig.~\ref{fig12}(c) presents the Darcy velocity distribution within the region $L_{\mathrm{cell}} \leq y \leq 7L_{\mathrm{cell}}$.
Fig.~\ref{fig12}(c) illustrates that the Darcy velocity distributions differ between the uniform-thickness and optimized-thickness designs.
The increased wall thickness in the central region raises flow resistance, leading to a local reduction in Darcy velocity.
As a result, whereas the flow in the uniform-thickness design proceeds more directly from the inlet to the outlet, the optimized-thickness design redirects more fluid toward the lateral regions of the core domain.

\begin{figure*}[h]
    \centering
    \begin{minipage}[t]{0.45\textwidth}
        \centering
        \includegraphics[width=0.9\linewidth]{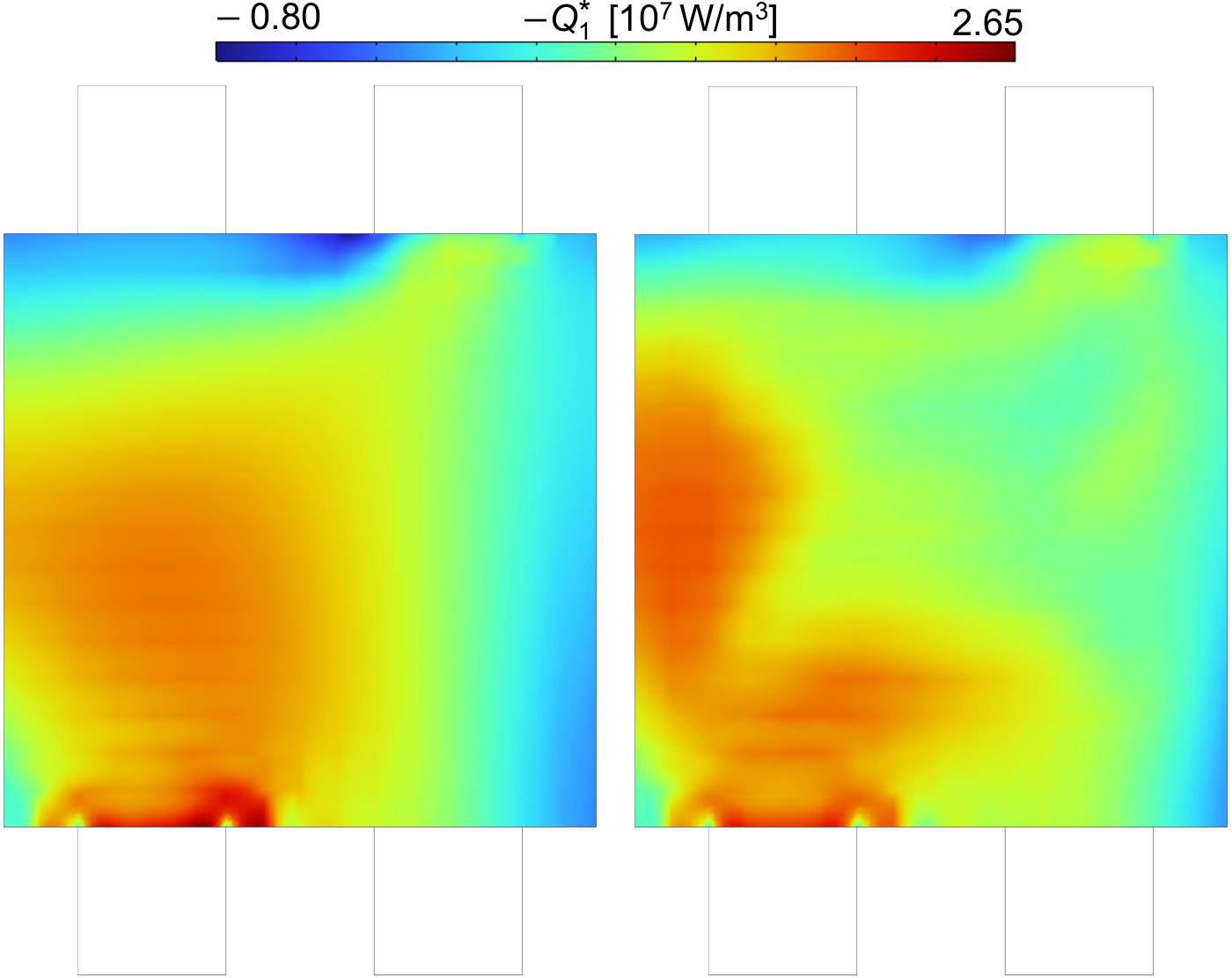}
        \subcaption{Distribution of the $-Q^*_1$}
        \label{fig13a}
    \end{minipage} 
    \begin{minipage}[t]{0.45\textwidth}
        \centering
        \includegraphics[width=0.9\linewidth]{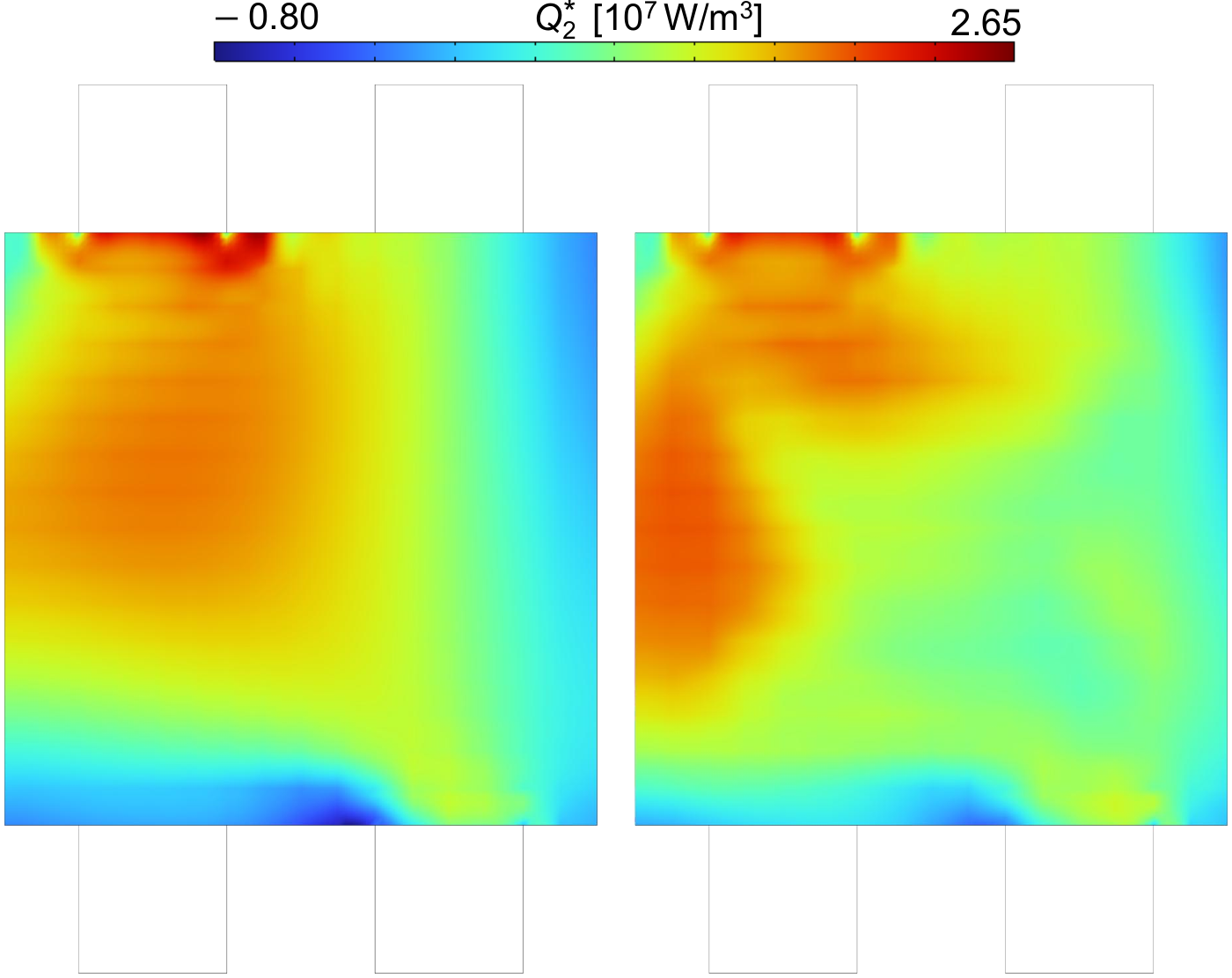}
        \subcaption{Distribution of the $Q^*_2$}
        \label{fig13b}
    \end{minipage} 
    \vspace{-0.2cm}
    \caption{Distributions of the volumetric heat source terms at the cross-section \( z = L_{\mathrm{cell}}/2 \) for both the initial (left) and optimized (right) designs}
    \vspace{-0.2cm}
    \label{fig13}
\end{figure*}

Fig.~\ref{fig13} presents the distributions of the volumetric heat source terms $Q^*_1$ and $Q^*_2$ for both the initial configuration ($\hat{\gamma} = 0.5$) and the optimized structure ($w = 0.25$).
Since the temperature of the hot fluid $T_1$ is generally higher than the wall temperature $T_\mathrm{w}$, the hot fluid transfers heat to the wall in most of the core domain.
In this study, the volumetric heat source term is defined as the heat received by the fluid per unit volume.
Therefore, $Q^*_1$ is generally negative, and $-Q^*_1$ serves as an indicator of the magnitude of heat exchange.
Negative values of $-Q^*_1$ or $Q^*_2$ appear in localized regions.
For the hot fluid, these negative values arise when its temperature falls below the wall temperature after passing near the inlet of the cold fluid.
Such occurrences are consistent with the physical behavior of the system.

The reasons for improved heat transfer in the optimized-thickness design are discussed.
Fig.~\ref{fig13} shows that the regions of enhanced heat transfer differ between the uniform-thickness and optimized-thickness designs.
In the uniform-thickness design, heat transfer is relatively high from the inlet toward the central region.
In contrast, the optimized-thickness design exhibits relatively high heat transfer in the left portion of the core domain.
Furthermore, the uniform-thickness design exhibits lower heat transfer in the right portion of the core.
As shown in Fig.~\ref{fig12}(c), in the optimized design, the increased central thickness redirects the flow toward the lateral regions, resulting in higher Darcy velocities in both the left and right portions of the core domain.
This local increase in Darcy velocity is considered to enhance heat transfer in these regions.
Consistent with these observations, the total heat transfer increased from 135.7 W in the uniform-thickness design with $\hat{\gamma}=0.5$ to 139.7 W in the optimized-thickness design with $w=0.25$.
These results indicate that the optimized-thickness design effectively utilizes the entire core domain for heat exchange, thereby improving the total heat transfer performance of the HX.

%%%%%%%%%%%%%%%%%%%%%%%%%%%%%%%%%%%%%%%%%
%%%%%%%%%%%%%%%%%%%%%%%%%%%%%%%%%%%%%%%%%
%%%%%%%%%%%%%%%%%%%%%%%%%%%%%%%%%%%%%%%%%
\begin{figure}[htb]
    \begin{center}
        \includegraphics[width=0.95\linewidth]{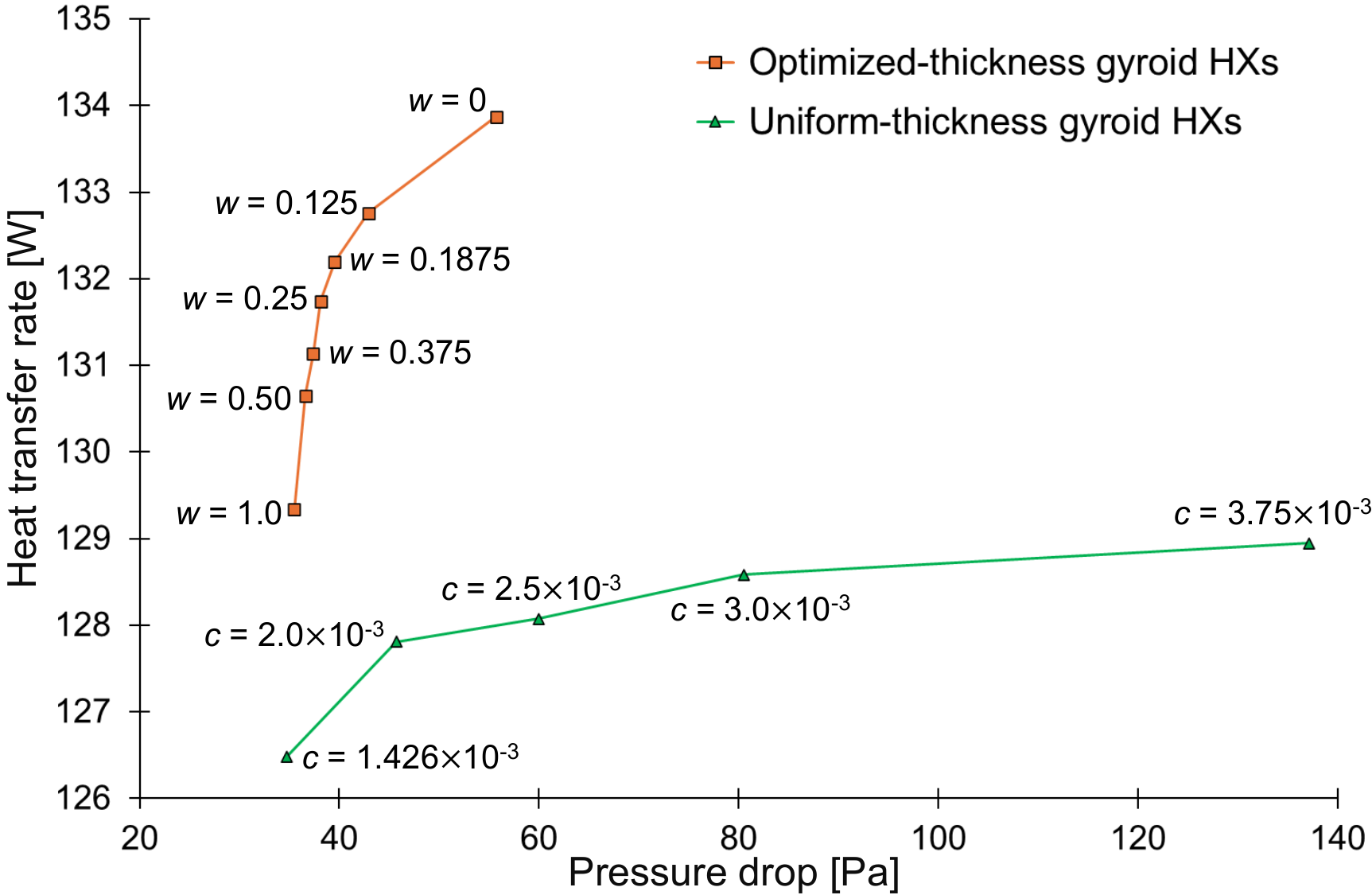}
        \caption{Full-scale simulation results for gyroid HXs with uniform and optimized thickness distribution}
        %\vspace{-0.3cm}
        \label{fig14}
    \end{center}%
    %\vspace{-0.2cm}
\end{figure}
%%%%%%%%%%%%%%%%%%%%%%%%%%%%%%%%%%%%%%%%%
\subsection{Full-Scale Performance Evaluation Results}
\label{sec4.6}

To further demonstrate the effectiveness of the optimized wall thickness distribution, we conducted dehomogenization and full-scale numerical analysis.
Thanks to filtering applied during the optimization, smoothly graded gyroid HXs were reconstructed from the optimized design variable distribution $\hat{\gamma}$.
To prevent mixing of the hot and cold fluids, 0.5 mm thick solid walls were appropriately placed at the interfaces between the fluid regions and the HX core.
In the full-scale analysis, the hot and cold fluids were given a uniform inlet velocity of 0.03 m/s, with inlet temperatures of 333.15 K and 293.15 K, respectively.
At the outlets, a pressure of 0 Pa was prescribed.  
All other boundaries were treated as adiabatic and no-slip boundaries.

Fig.~\ref{fig14} shows the results of the full-scale analysis. 
The gyroid HXs with optimized thickness distribution exhibited improved heat exchange performance compared to those with uniform thickness distribution.
This indicates that locally varying the gyroid HX thickness using the proposed method is more effective than a uniform increase in thickness for enhancing heat transfer.
Furthermore, as the weighting factor $w$ increased, the pressure drop decreased, demonstrating that the optimization using effective model appropriately captured the trade-off between heat transfer and pressure drop.

\begin{figure*}[htb]
    \centering
    \begin{minipage}[t]{0.45\textwidth}
        \centering
        \includegraphics[width=0.9\linewidth]{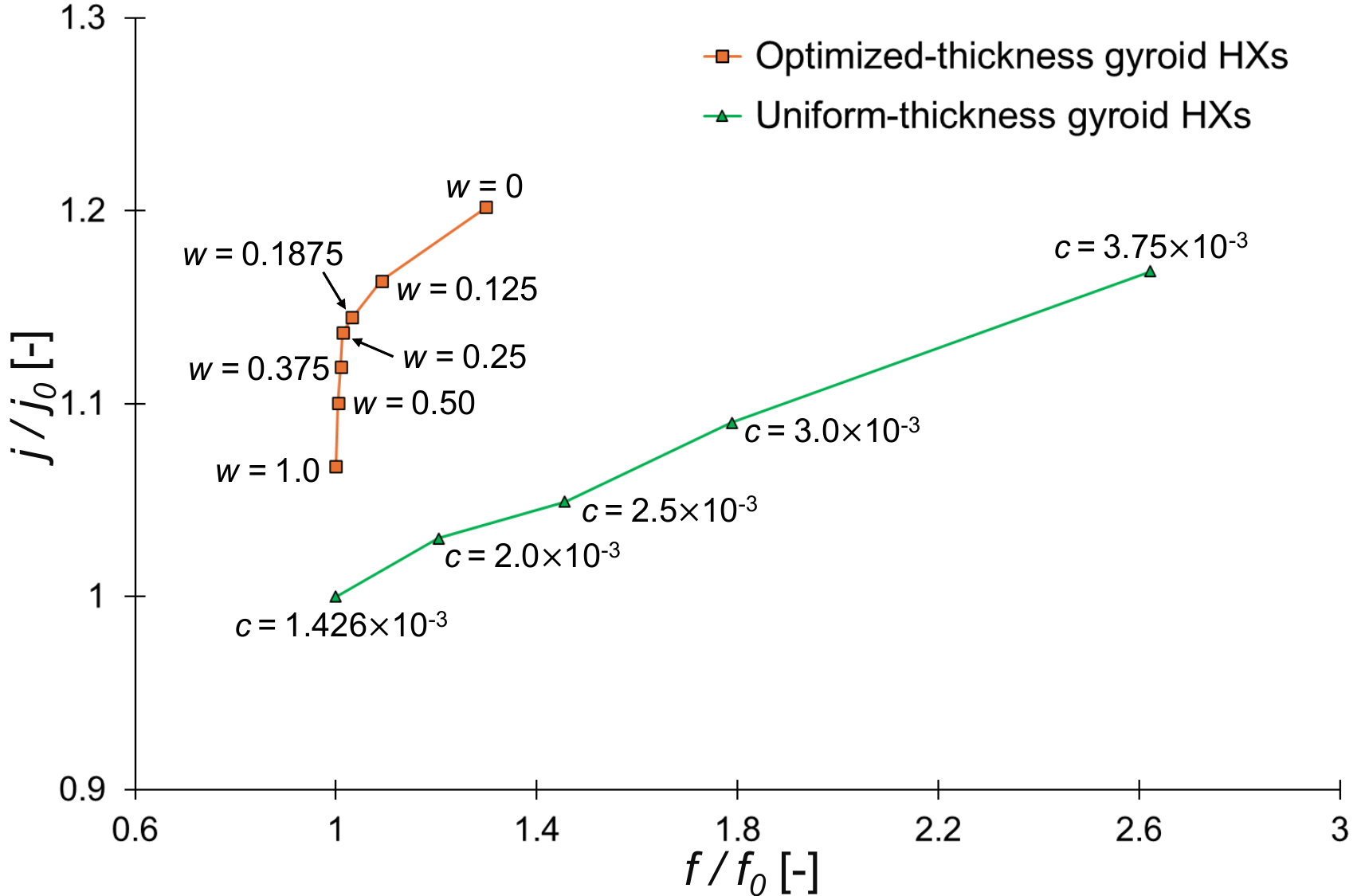}
        \subcaption{Normalized $j$ vs normalized $f$}
        \label{fig15a}
    \end{minipage}
    \begin{minipage}[t]{0.45\textwidth}
        \centering
        \includegraphics[width=0.9\linewidth]{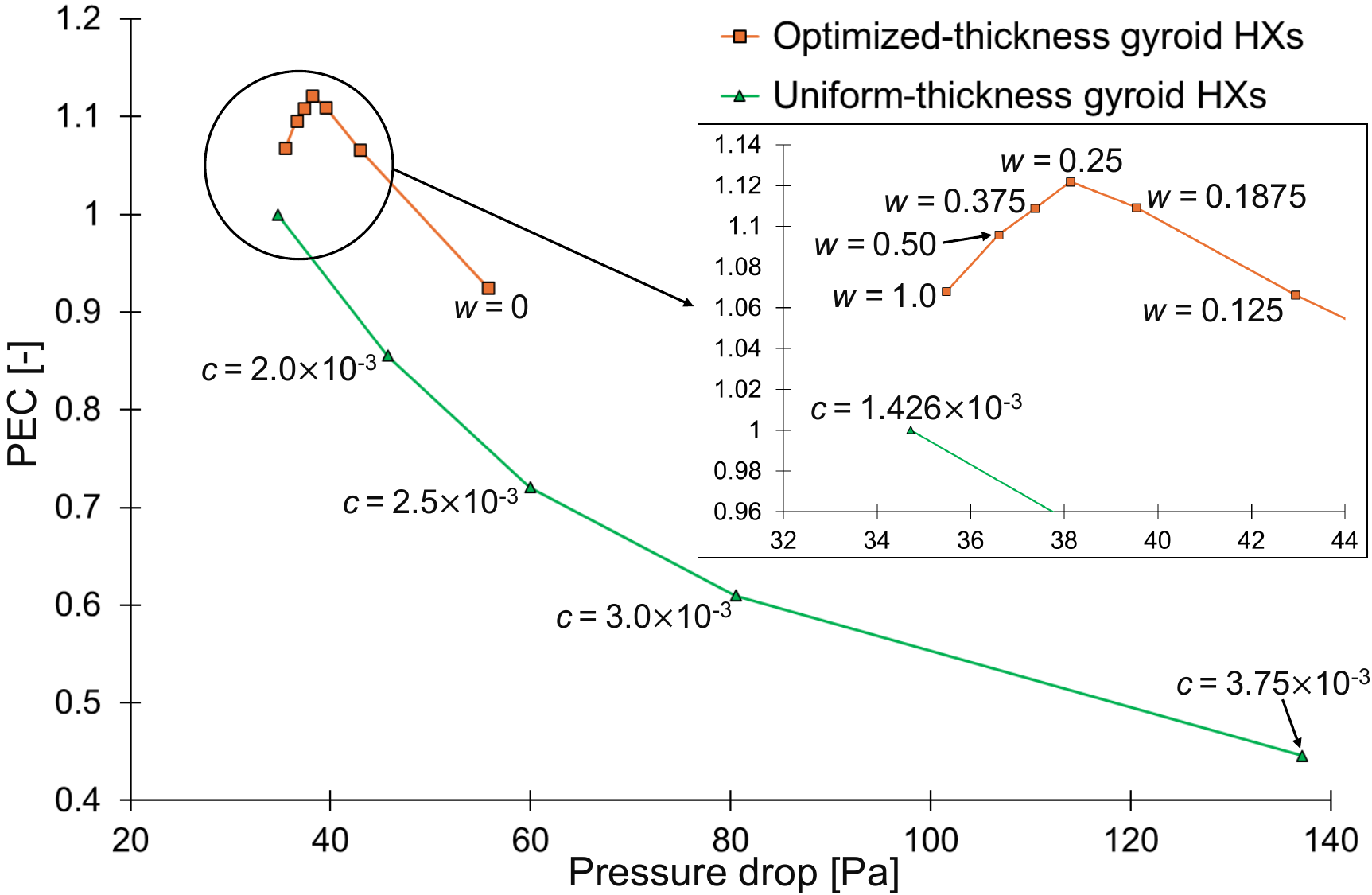}
        \subcaption{PEC values vs pressure drop}
        \label{fig15b}
    \end{minipage}
    \vspace{-0.2cm}
    \caption{Dimensionless performance factors for gyroid HXs with uniform and optimized wall thickness distribution}
    \label{fig15}
\end{figure*}

The trade-off between heat transfer and pressure drop was further evaluated.  
We used the PEC, which is a widely used performance metric for HXs that considers both heat transfer and pressure drop, to identify the most effective design.
Fig.~\ref{fig15} presents dimensionless performance factors for gyroid HXs with uniform-thickness and optimized-thickness distributions.
Fig.~\ref{fig15}(a) shows the scatter plot of normalized $j$ versus normalized $f$.
Fig.~\ref{fig15}(b) presents the corresponding PEC values for each design.
The baseline values of $j_0$ and $f_0$ were calculated for the uniform-thickness design ($c = 1.426 \times 10^{-3}$), which has the thinnest wall across the entire core domain.

For the uniform-thickness design, Fig.~\ref{fig15}(a) shows that increasing the wall thickness raises $j/j_0$ to a maximum of 1.17, while $f/f_0$ increases to 2.62.
This is because a uniform increase in wall thickness reduces the effective flow cross-sectional area, thereby increasing the mean flow velocity and enhancing the heat transfer coefficient, which leads to greater heat exchange.
However, the accompanying increase in flow resistance results in a significant rise in pressure loss.
Consequently, from the perspective of PEC, which reflects the trade-off between heat transfer and pressure drop, the HX with the thinnest wall ($c = 1.426 \times 10^{-3}$) was optimal among the uniform-thickness designs.
This result indicates that, in uniform-thickness designs, reducing the wall thickness effectively enhances HX performance by simultaneously decreasing both flow resistance and conductive thermal resistance.

%%%%%%%%%%%%%%%%%%%%%%%%%%%%%%%%%%%%%%%%%
%%%%%%%%%%%%%%%%%%%%%%%%%%%%%%%%%%%%%%%
%%%%%%%%%%%%%%%%%%%%%%%%%%%%%%%%%%%%%%%%%
\begin{figure*}[!t]
    \centering
    \begin{minipage}[t]{0.45\textwidth}
        \centering
        \includegraphics[width=0.84\linewidth]{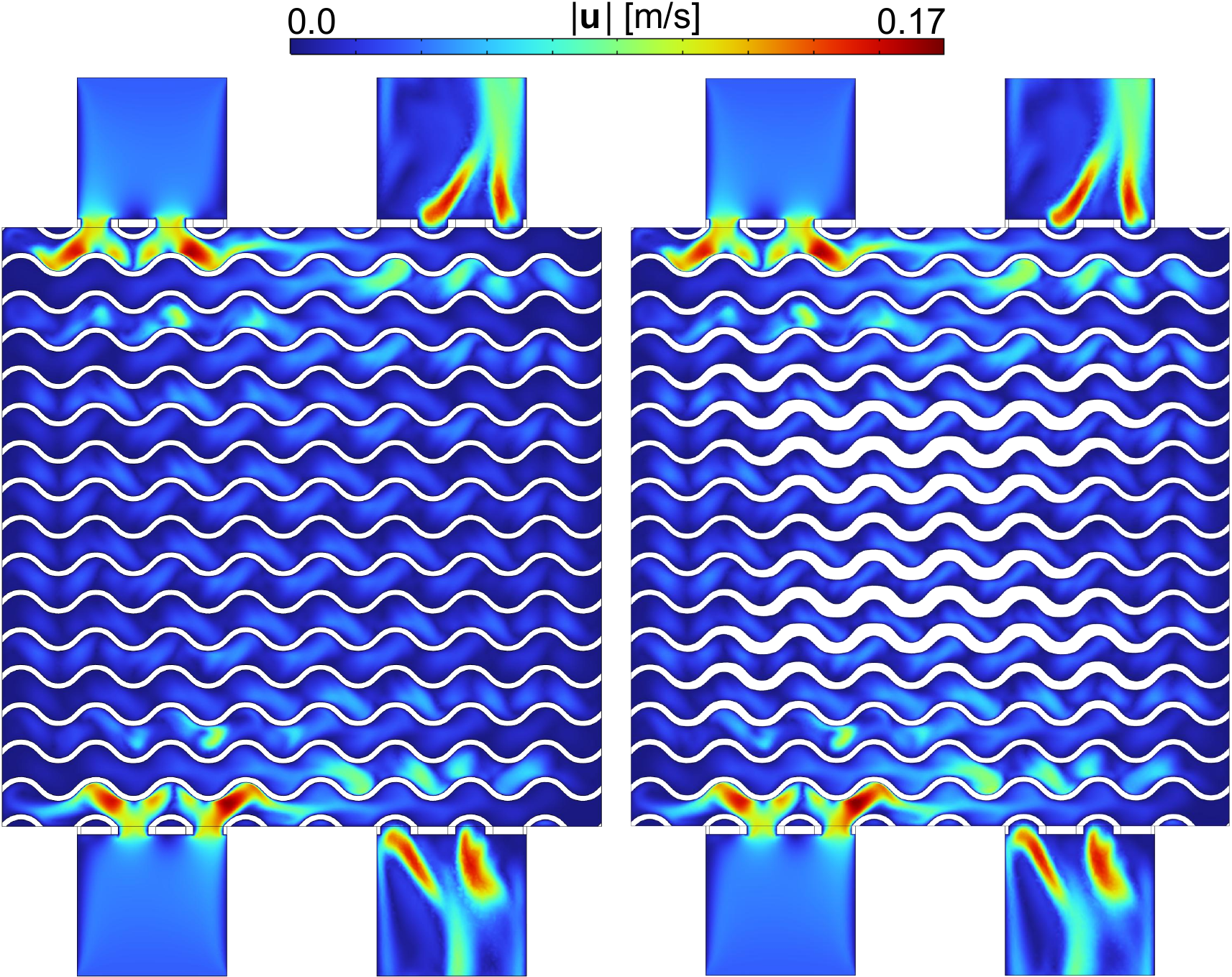}
        \subcaption{Velocity distribution}
        \label{fig16a}
    \end{minipage} 
    \vspace{0.2cm}
    \begin{minipage}[t]{0.45\textwidth}
        \centering
        \includegraphics[width=0.84\linewidth]{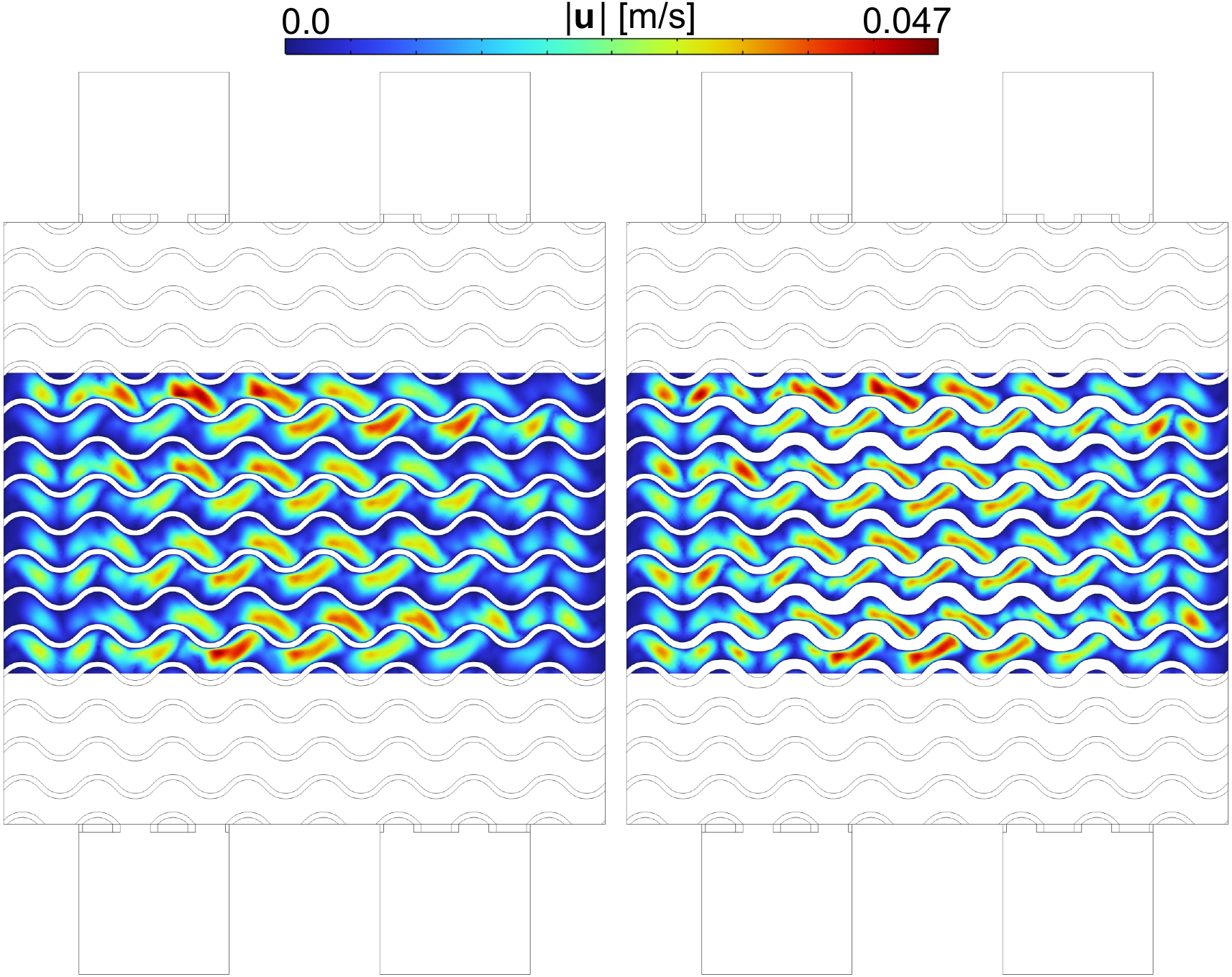}
        \subcaption{Velocity distribution within the region \( 2 L_{\mathrm{cell}} \leq  y \leq 6 L_{\mathrm{cell}} \)}
        \label{fig16b}
    \end{minipage} \\
    \vspace{0.2cm}
    \begin{minipage}[t]{0.45\textwidth}
        \centering
        \includegraphics[width=0.84\linewidth]{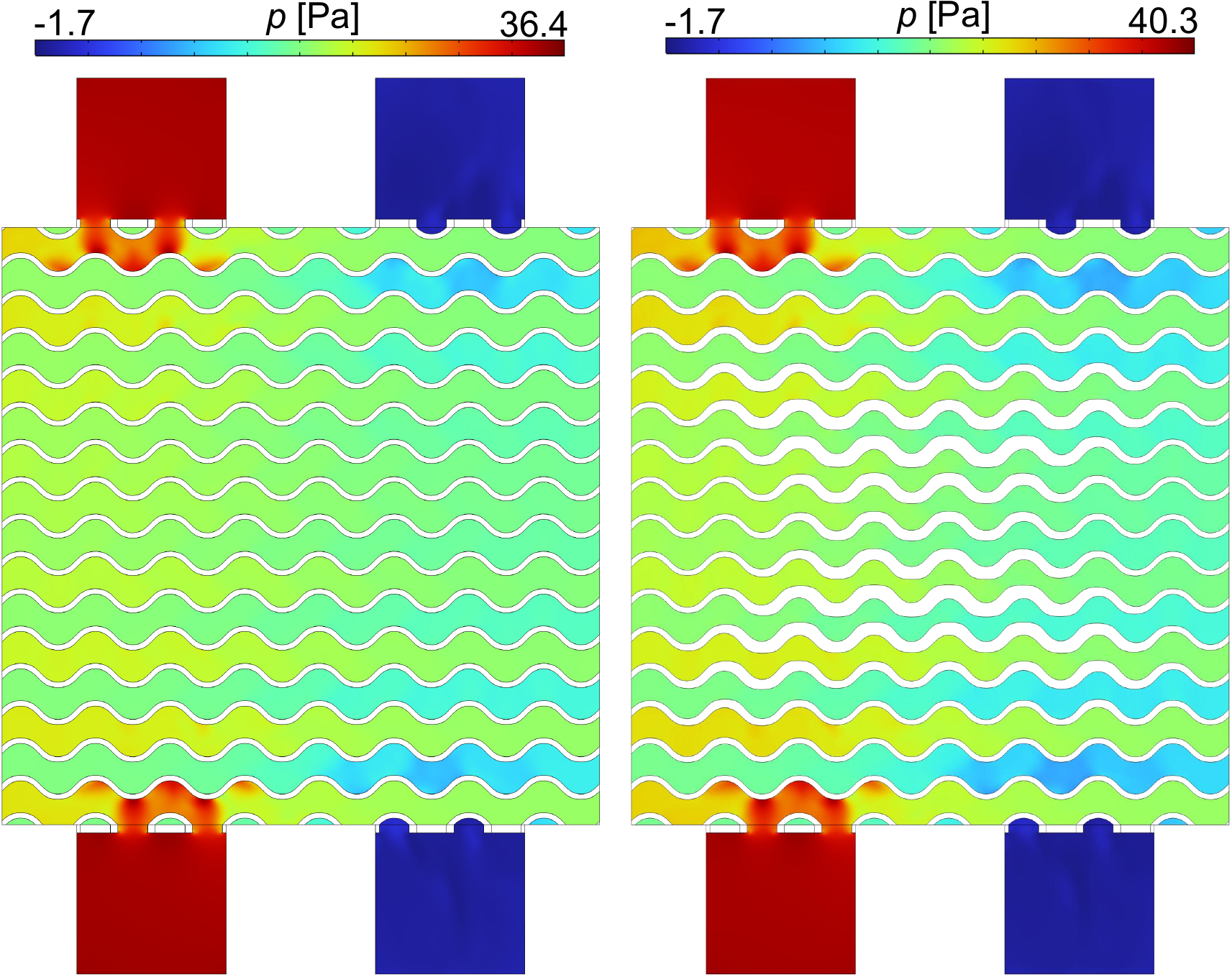}
        \subcaption{Pressure distribution}
        \label{fig16c}
    \end{minipage}
    \vspace{0.2cm}
    \begin{minipage}[t]{0.45\textwidth}
        \centering
        \includegraphics[width=0.84\linewidth]{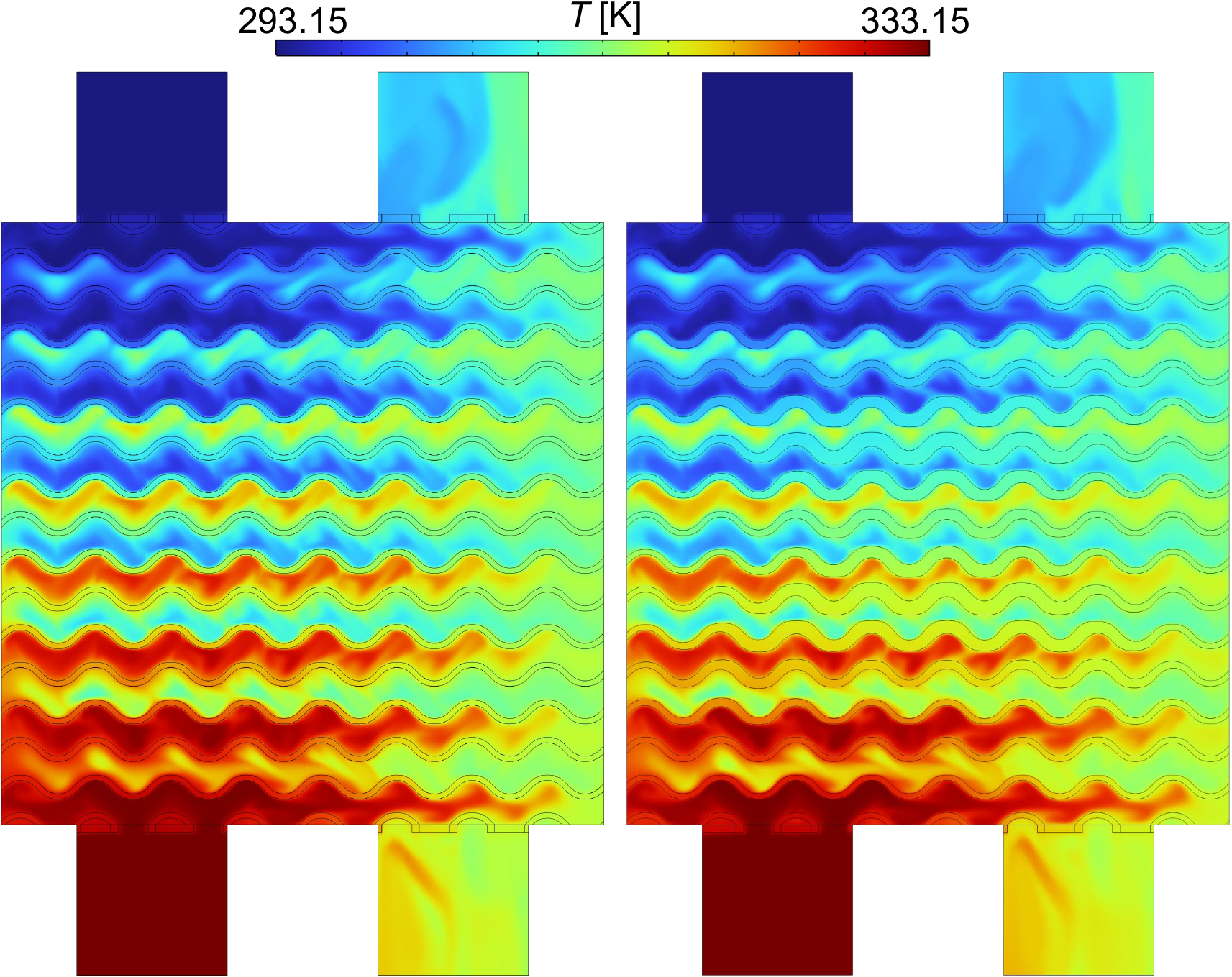}
        \subcaption{Temperature distribution}
        \label{fig16d}
    \end{minipage}\\
    \vspace{-0.2cm}
    \caption{Physical fields at the cross-section \( z = L_{\mathrm{cell}}/2 \) for the uniform-thickness (left, \( c = 1.426 \times 10^{-3} \)) and optimized-thickness (right, \( w = 0.25 \)) gyroid HXs}
    \vspace{-0.2cm}
    \label{fig16}
\end{figure*}
%%%%%%%%%%%%%%%%%%%%%%%%%%%%%%%%%%%%%%%%%
In contrast, as shown in Fig.~\ref{fig15}(b), the optimized-thickness designs achieved higher PEC values than the uniform-thickness design  for all weighting factors $w$ except $w=0$.
As an exception, for $w = 0$, the optimization did not account for pressure drop, resulting in a lower PEC value.
The highest PEC value, 1.122, was obtained for the design optimized with $w=0.25$.
This indicates that, compared to the uniform-thickness design with $c = 1.426 \times 10^{-3}$, which was the best-performing among the uniform-thickness designs, the optimized-thickness design achieves a 12.2\% increase in efficiency, considering the trade-off between heat transfer and pressure drop.
Heat transfer rates and pressure drops for the uniform-thickness ($c=1.426 \times 10^{-3}$) and optimized-thickness ($w=0.25$) designs were 126.5 W and 131.7 W, and 34.7 Pa and 38.1 Pa, respectively.
The volumetric heat densities for the designs with uniform and optimized thickness are $1.015 \times 10^7~\mathrm{W/m^3}$ and $1.057 \times 10^7~\mathrm{W/m^3}$, respectively, showing a 4.16\% improvement through optimization.

These results demonstrate that introducing the optimized non-uniform wall thickness distribution improves the efficiency of gyroid HXs.
Moreover, this improvement in heat transfer is attributed to mechanisms unique to non-uniform thickness distributions.
The following section examines the mechanisms underlying the performance enhancement in HXs with the optimized thickness distribution.

%%%%%%%%%%%%%%%%%%%%%%%%%%%%%%%%%%%%%%%%%
\begin{figure*}[!t]
    \begin{center}
        \centering
    \begin{minipage}[t]{0.8\textwidth}
        \centering
        \includegraphics[width=0.9\linewidth]{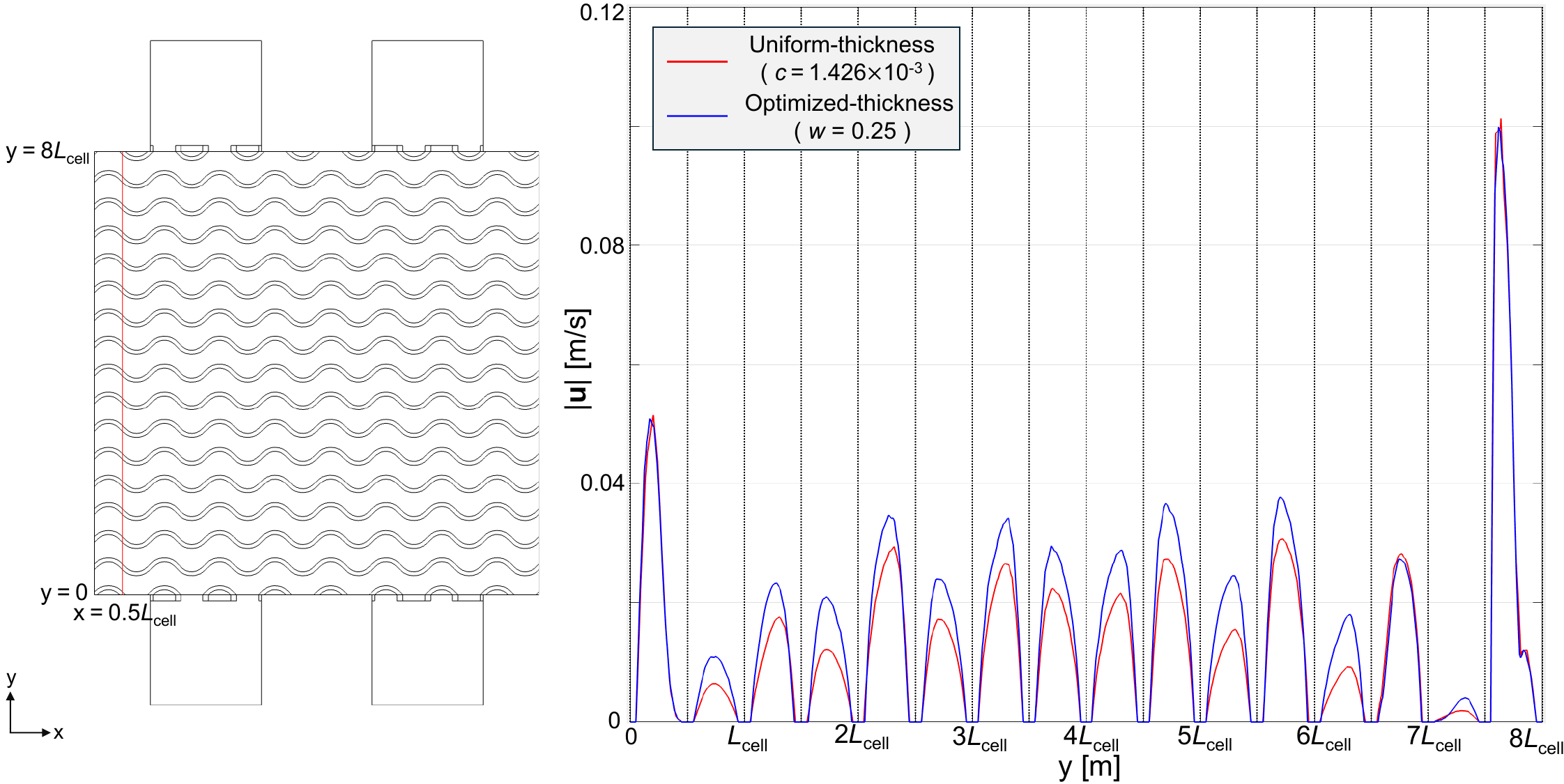}
        \subcaption{Velocity magnitude along $y$ at $x=0.5L_{\mathrm{cell}}$ in the cross-section.}
        \label{fig17a}
    \end{minipage} \\
    \vspace{0.2cm}
    \begin{minipage}[t]{0.8\textwidth}
        \centering
        \includegraphics[width=0.9\linewidth]{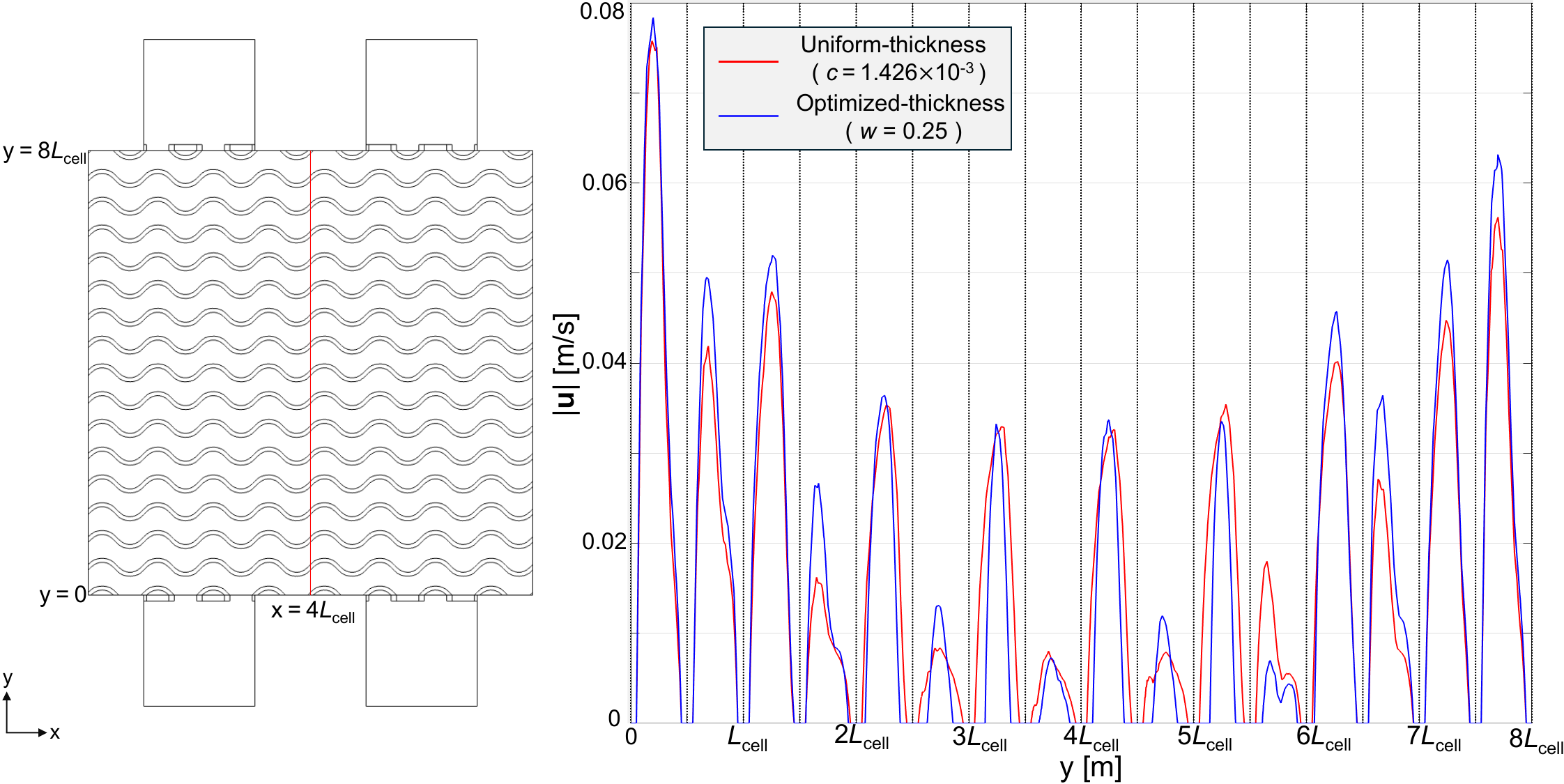}
        \subcaption{Velocity magnitude along $y$ at $x=4.0L_{\mathrm{cell}}$ in the cross-section.}
        \label{fig17b}
    \end{minipage} \\
    \vspace{0.2cm}
    \begin{minipage}[t]{0.8\textwidth}
        \centering
        \includegraphics[width=0.9\linewidth]{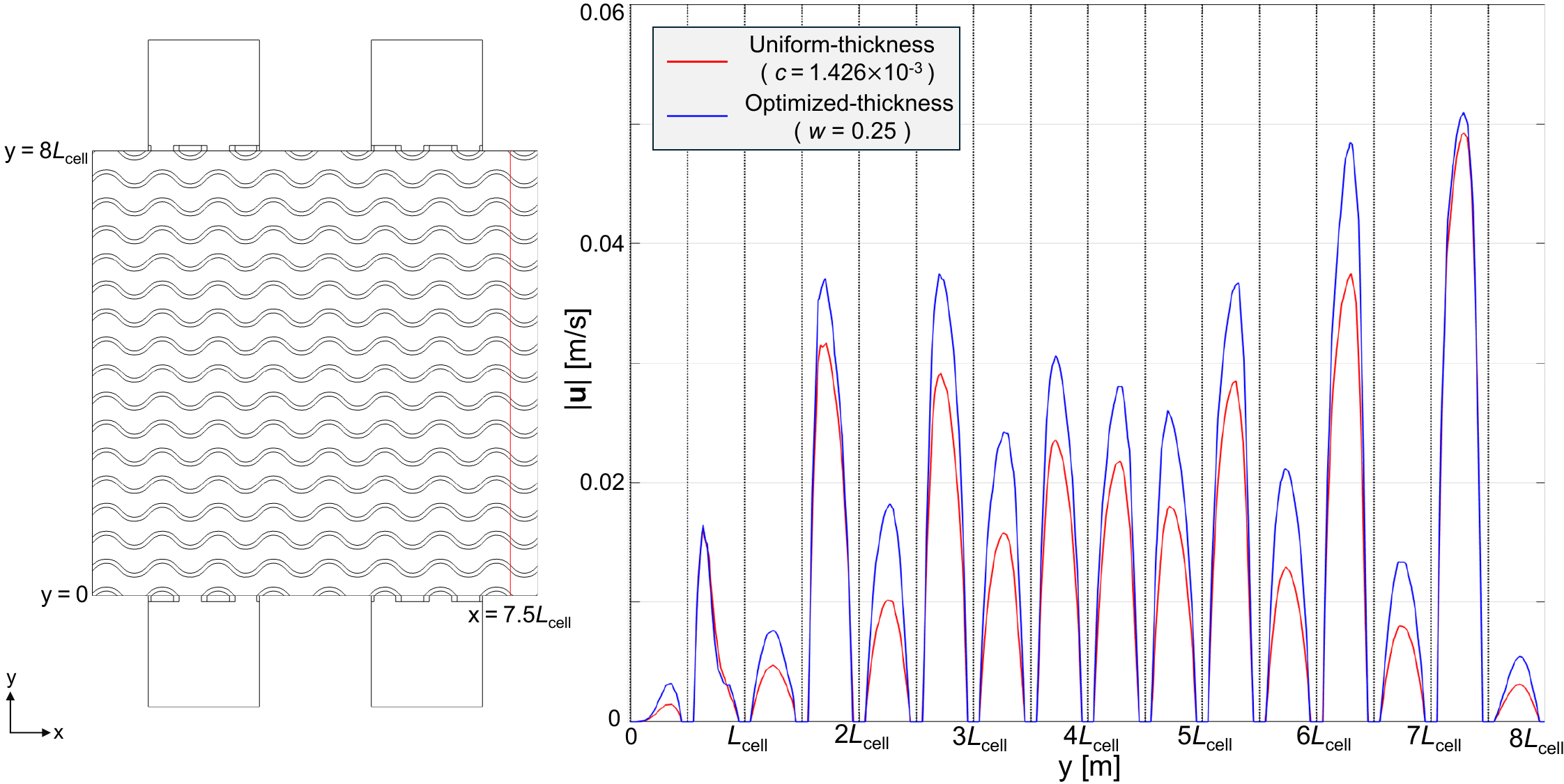}
        \subcaption{Velocity magnitude along $y$ at $x=7.5L_{\mathrm{cell}}$ in the cross-section.}
        \label{fig17c}
    \end{minipage} \\
        \caption{Velocity magnitude profiles at different positions in the cross-section}
        %\vspace{-0.3cm}
        \label{fig17}
    \end{center}%
\end{figure*}
%%%%%%%%%%%%%%%%%%%%%%%%%%%%%%%%%%%%%%%%%
%%%%%%%%%%%%%%%%%%%%%%%%%%%%%%%%%%%%%%%%%
\begin{figure}[htb]
    \begin{center}
        \includegraphics[width=0.9\linewidth]{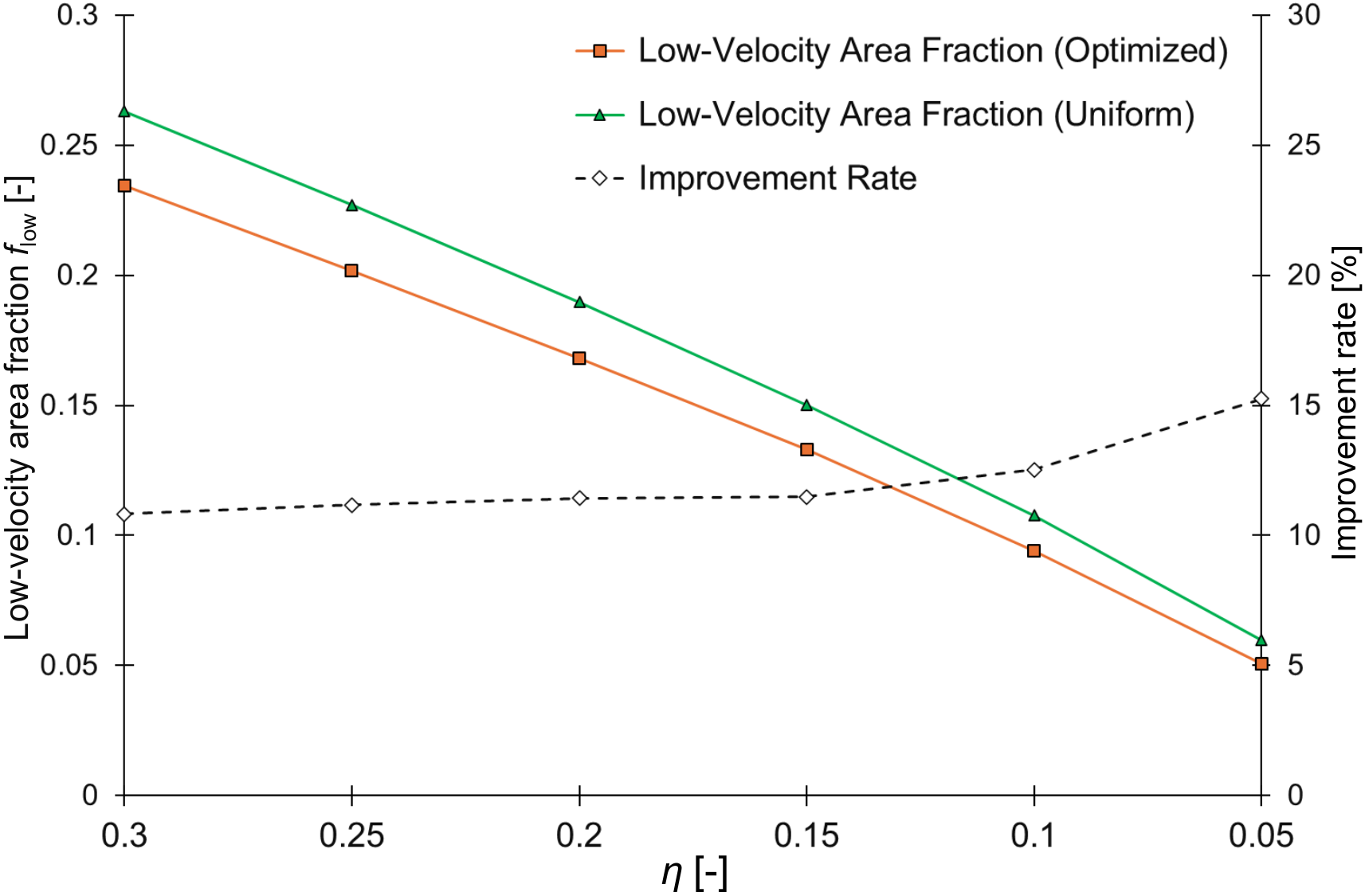}
        \caption{Low-velocity area fraction $f_{\mathrm{low}}$ for uniform and optimizedthickness designs at $z = L_{\rm cell}/2$, with improvement rate}
        %\vspace{-0.3cm}
        \label{fig18}
    \end{center}%
    %\vspace{-0.2cm}
\end{figure}
%%%%%%%%%%%%%%%%%%%%%%%%%%%%%%%%%%%%%%%%%
\subsection{Discussion}
\label{sec4.7}

To investigate the mechanisms responsible for the performance improvement with the optimized-thickness distribution, the physical fields obtained from full-scale simulations of both the optimized and baseline designs were compared.
Fig.~\ref{fig16} shows the physical fields at $z = L_{\mathrm{cell}}/2$ for both the uniform-thickness ($c = 1.426 \times 10^{-3}$, left) and optimized-thickness ($w = 0.25$, right) designs.
Fig.~\ref{fig16} shows the velocity distributions in (a) and (b), the pressure distribution in (c), and the temperature distributions in (d).

First, we focus on the velocity distribution. 
As shown in Fig.~\ref{fig16}(b), it can be seen that the optimized-thickness distribution leads to higher velocities near the left and right edges of the HX core region.
To clearly illustrate the differences in velocity distribution, Fig.~\ref{fig17} presents the velocity magnitude along the lines at $x = 0.5 L_{\mathrm{cell}}$, $x = 4.0 L_{\mathrm{cell}}$, and $x = 7.5 L_{\mathrm{cell}}$ within the cross-section.
From Fig.~\ref{fig17}(b), it can be observed that at the center of the HX ($x = 4.0 L_{\mathrm{cell}}$), there is no significant difference in velocity between the two designs.
In contrast, near the edges of the heat exchange region (Figs.~\ref{fig17}(a) and \ref{fig17}(c)), the optimized-thickness design exhibits higher velocity magnitudes than the uniform-thickness design.
These local increases in velocity are considered to contribute to the enhancement of heat transfer performance.
In the literature on conventional HXs, non-uniform velocity distributions have been shown to lead to performance degradation \cite{Ranganayakulu1996,Zhang2009,Singh2014}. 
Similarly, several studies on TPMS HXs have shown that improving velocity uniformity can enhance overall heat transfer efficiency, as it promotes more effective heat exchange throughout the entire HX \cite{Zhang2025conformal,Zhang2025gradient,Oh2023,Oh2025,Wang2024}.

Therefore, we quantitatively assessed velocity uniformity on the cross-section at $z = L_{\rm cell}/2$ using two metrics: the low-velocity area fraction $f_{\mathrm{low}}$ and the coefficient of variation (CV) of the velocity field. 
The $f_{\mathrm{low}}$ represents the portion of the fluid cross-sectional area where the flow velocity is below $\eta$ times the mean velocity $\bar{U}$, and is defined using the indicator function $I(\mathbf{u}(\mathbf{x}))$ as follows:
\begin{equation}
f_{\mathrm{low}} = \frac{\displaystyle \int_S I(\mathbf{u}(\mathbf{x})) \, dS}{\displaystyle \int_S \ dS},
\end{equation}
where
\begin{equation}
I(\mathbf{u}(\mathbf{x})) =
\begin{cases}
1, & \text{if } |\mathbf{u}(\mathbf{x})| \le \eta \bar{U},\\
0, & \text{otherwise,}
\end{cases}
\end{equation}
and $\bar{U}$ is the mean velocity over the cross-section:
\begin{equation}
\bar{U} = \frac{\displaystyle \int_S |\mathbf{u}(\mathbf{x})| \, dS}{\displaystyle \int_S dS}.
\end{equation}
The CV quantifies the overall velocity variability across the fluid region $S$ and is defined as:
\begin{equation}
\text{CV} = \frac{\sqrt{\displaystyle \frac{1}{\int_S dS} \int_S \big(|\mathbf{u}(\mathbf{x})| - \bar{U}\big)^2 \, dS}}{\bar{U}}.
\end{equation}
The CV, defined as the velocity standard deviation $\sigma_\mathbf{u}$ normalized by the mean velocity $\bar{U}$, offers a consistent measure of velocity variability independent of the mean flow magnitude.

By calculating both $f_{\mathrm{low}}$ and CV for the uniform-thickness and optimized-thickness designs, the reduction of low-velocity regions and the improvement in velocity uniformity can be quantitatively assessed.
The mean velocities $\bar{U}$ were 0.0225 m/s for the uniform-thickness design and 0.0250 m/s for the optimized-thickness design.
Despite the same inlet flow rate, the optimized design shows a higher mean velocity because the locally increased wall thickness reduces the effective flow cross-sectional area.
The calculated low-velocity area fraction $f_{\mathrm{low}}$ is shown in Fig.~\ref{fig18}.
Here, the values were computed for $\eta$ ranging from 0.3 to 0.05.
The corresponding improvement rates are also indicated in the figure and defined as:
\begin{equation}
\text{Improvement rate} =  \left( 1 - \frac{f_{\mathrm{low, opt}}}{f_{\mathrm{low, uni}}} \right)\times100\%,
\end{equation}
where $f_{\mathrm{low, opt}}$ and $f_{\mathrm{low, uni}}$ denote the low-velocity area fractions for the optimized-thickness and uniform-thickness designs, respectively.
As shown in Fig.~\ref{fig18}, $f_{\mathrm{low}}$ is smaller for the optimized-thickness design for all values of $\eta$, with improvement rates ranging from 10.8\% to 15.3\%. 
These results quantitatively demonstrate that the optimized thickness distribution effectively reduces low-velocity regions.
Additionally, the CV was calculated as 1.063 for the uniform-thickness design and 0.975 for the optimized-thickness design, representing an 8.3\% reduction in velocity variability.
Therefore, we conclude that the optimized-thickness distribution reduces low-velocity regions and enhances velocity uniformity, allowing more effective utilization of the entire HX core for heat transfer.

%%%%%%%%%%%%%%%%%%%%%%%%%%%%%%%%%%%%%%%%%
\begin{figure*}[htb]
    \begin{center}
        \includegraphics[width=0.95\linewidth]{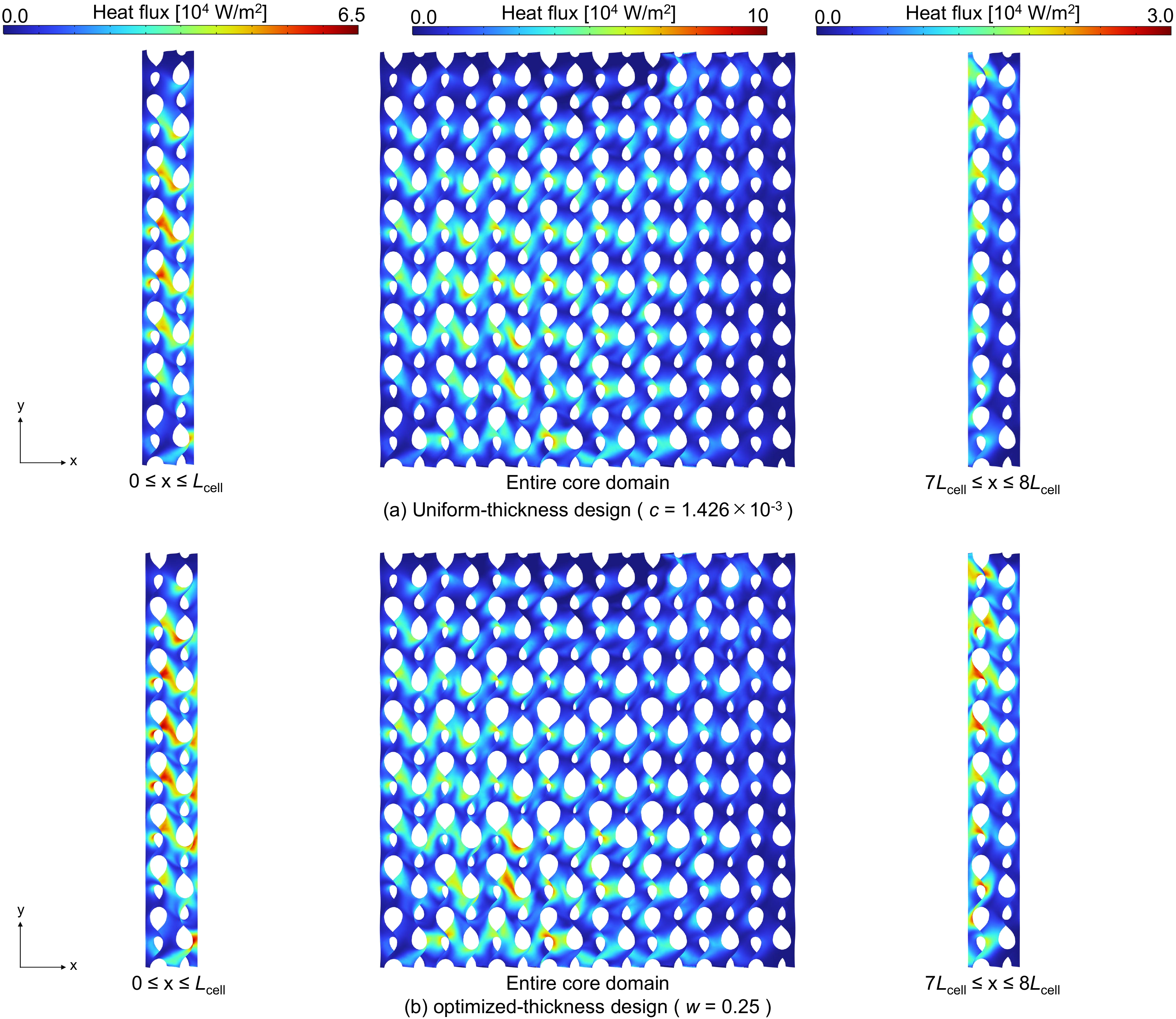}
        \caption{Heat flux distribution for uniform-thickness (top, \( c = 1.426 \times 10^{-3} \)) and optimized-thickness (bottom, \( w = 0.25 \)) gyroid HXs}
        %\vspace{-0.3cm}
        \label{fig19}
    \end{center}%
    %\vspace{-0.2cm}
\end{figure*}
%%%%%%%%%%%%%%%%%%%%%%%%%%%%%%%%%%%%%%%%%
To further demonstrate the advantage of the optimized-thickness distribution, the wall heat flux is visualized, clarifying how the improved velocity uniformity contributes to enhanced heat transfer performance.
Fig.~\ref{fig19} shows the heat flux distributions on the wall in contact with the hot fluid within the region $0.25 L_{\mathrm{cell}} \leq z \leq 0.75 L_{\mathrm{cell}}$ for both the uniform- and optimized-thickness designs.
As shown in Fig.~\ref{fig17}, the optimized-thickness distribution produces higher velocities near the edges of the core region.
Accordingly, Fig.~\ref{fig19} presents the heat flux distributions specifically in these edge regions ($0 \leq x \leq L_{\mathrm{cell}}$, $ 7 L_{\mathrm{cell}} \leq x \leq 8 L_{\mathrm{cell}}$).
A comparison between the uniform- and optimized-thickness designs confirms that the optimized distribution results in higher heat flux values at the edges.
This increase arises from the locally higher velocities, which strengthen convective heat transfer in these regions.
The overall heat transfer coefficient $U$, calculated from Eq.~\eqref{eq:U}, provides a quantitative measure of the heat transfer performance. For the optimized thickness distribution, $U$ reaches 577 W/(m$^2$ K), representing an approximately 14\% increase compared to $U=508$ W/(m$^2$ K) for the uniform-thickness design.
These results further indicate that the optimization of the thickness distribution, by enhancing velocity uniformity, effectively improves heat transfer performance.

In summary, the optimized-thickness design features a non-uniform thickness distribution with increased wall thickness near the center, which raises flow resistance in the central region and redirects more flow toward the left- and right-edge regions of the HX core.
This increases velocities in these regions, reducing low-velocity areas and improving velocity uniformity.
Consequently, as shown in Fig.~\ref{fig19}, heat transfer is enhanced at the core edges, allowing for more effective utilization of the entire core and improving heat exchange performance.

Finally, we discuss the discrepancies observed between the effective model and the full-scale simulation for the optimized design.  
Comparing Figs.~\ref{fig11} and \ref{fig14}, the full-scale simulation results in a higher pressure drop, whereas the effective model predicts a larger heat transfer rate.

We first discuss the discrepancy in pressure drop.
In the full-scale simulation, solid walls were introduced to prevent mixing of the hot and cold flow streams.
These walls narrow the flow paths near the inlets and outlets, forcing the fluid to converge into localized regions, as seen in the velocity distribution in Fig.~\ref{fig16}(a).
This results in an increase in pressure loss due to flow acceleration and redirection caused by fluid concentration.
The effective model cannot capture such geometrical and flow-concentration effects because the model assumes spatially averaged behavior and a more uniform flow field.
Furthermore, in the verification case of Section~\ref{sec4.3}, these solid walls were not included in the full-scale model, and the flow configuration closely resembled that assumed in the effective model.
As a result, the discrepancy in pressure drop between the effective and full-scale models was relatively small in verification.
Therefore, the presence of solid walls, not considered in the effective model, appears to be the primary cause of the pressure drop discrepancy between the effective and full-scale analyses of the optimized HX designs.

We next discuss the discrepancy in heat transfer rate.  
As described in Section~\ref{sec2.3}, the effective heat transfer coefficient is derived under the assumption of counter-flow within the RVE.  
Accordingly, the effective model assumes a counter-flow configuration throughout the entire core domain.  
This assumption remained valid in the verification case in Section~\ref{sec4.3}, where the full-scale geometry also maintained a counter-flow layout, leading to good agreement between the two models.  
In a counter-flow configuration, the hot fluid inlet aligns with the cold fluid outlet, maintaining a large temperature difference and ensuring efficient heat transfer.
However, within the optimization problem setup, there are regions where the flow is partially parallel-flow or cross-flow.
In these regions, the fluid temperatures tend to equalize more rapidly, reducing the local temperature gradient and thereby decreasing heat transfer efficiency.
Consequently, the full-scale simulation yields a lower heat transfer rate than the effective model.
The effective model overestimates heat exchange performance due to its assumption of ideal counter-flow.

These discrepancies arise from inherent modeling assumptions in the effective model and dehomogenization process, which are difficult to fully eliminate.
Nevertheless, the effective model remains highly valuable for design optimization, as it significantly reduces computational cost while generating optimized-thickness designs that exhibit higher performance than uniform-thickness designs.
This confirms that the TO framework proposed in this study is both practical and effective for the design of gyroid two-fluid HXs.

\section{Conclusion}
\label{sec5}

TPMS HXs exhibit excellent thermal performance and hold significant potential as compact, high-efficiency next-generation thermal devices.
In this study, we proposed an optimization method for TPMS two-fluid HXs to optimize the wall thickness distribution.
To enable optimization with reasonable computational cost, an effective porous media model was developed for TPMS HXs.
We introduced an effective heat transfer coefficient in the effective model to represent the heat exchange between the fluid and the solid wall, and implemented this heat exchange as a volumetric heat source term.
The proposed effective model was applied to gyroid two-fluid HXs, and the wall thickness distribution was optimized considering both pressure drop and heat transfer rate.
Furthermore, multiple optimized design variable fields were de-homogenized to reconstruct the corresponding gyroid HXs, and full-scale numerical analyses were conducted to demonstrate the effectiveness of the proposed method.
The main findings of this study are summarized as follows:
\begin{itemize}
    \item The proposed effective model for TPMS HXs was validated through comparison with full-scale model simulations.
    Under counterflow conditions, where the RVE used to calculate the effective properties assumes such flow, the maximum relative errors between the effective and full-scale models were 6.72\% for pressure drop and 0.83\% for heat transfer rate.
    These results indicate that the effective model can accurately predict the performance of the full-scale model, confirming its validity.
    \item Optimization of the wall thickness distribution of the gyroid HXs using the effective model yielded optimized designs with thicker walls near the central region of the core domain.
    Analysis of the effective model’s physical fields for the optimized designs indicates that the increased central thickness redirects the flow toward the lateral regions, improving velocity uniformity.
    This more uniform flow allows effective utilization of the entire core domain for heat exchange, thereby increasing the total heat transfer compared to the uniform-thickness design.
    \item Full-scale simulations of the gyroid HX with an optimized-thickness distribution demonstrated superior heat transfer performance compared to the uniform-thickness HX.
    Specifically, the optimized-thickness design at a weighting factor of $w = 0.25$ achieved a maximum PEC of 1.122, corresponding to a 12.2\% improvement in overall performance over the uniform-thickness gyroid HX when accounting for the trade-off between heat transfer rate and pressure drop.
    Furthermore, the volumetric heat density increased by 4.16\%, further confirming the enhanced thermal performance of the optimized-thickness design.
    \item Full-scale analysis of the physical fields showed that the optimized wall-thickness distribution suppresses low-velocity regions and improves velocity uniformity.
    In particular, the optimized design increases the flow velocity near the core edges compared with the uniform-thickness design, thereby enhancing local heat transfer.
    Consequently, the optimized configuration more effectively utilizes the entire HX core domain for heat exchange.
    These results highlight the effectiveness of wall-thickness optimization in improving the performance of gyroid HXs.
\end{itemize}

This paper demonstrated the effectiveness of the proposed optimization method in enhancing overall performance for gyroid HXs.
By calculating effective properties based on the TPMS unit cell and optimizing the wall thickness distribution using the proposed effective model, the framework can be extended to various TPMS structures.
Nonetheless, some error due to the porous media approximation is inevitable, and this study focused solely on wall thickness as the design variable, leaving room for further improvement.
Future work will focus on improving the predictive accuracy of the effective model and extending the optimization framework to additional design variables, such as unit cell size and multiple TPMS types, enabling more advanced optimizations and further enhancing the performance of TPMS HXs.
Furthermore, fabricating gyroid HXs with the optimized wall thickness and experimentally validating their performance remain important objectives for future research.

\section*{Acknowledgments}
This work was supported by JST SPRING (Grant Number JPMJSP2138) and JSPS KAKENHI (Grant No. 23K26018).

% \appendix
% \section{Example Appendix Section}

\bibliographystyle{elsarticle-num} % 数字順スタイル
\bibliography{references}         % references.bib を読み込む

\end{document}